\documentclass[reqno,12pt]{amsart}
\usepackage{a4wide}
\usepackage{amsmath}
\usepackage{amssymb}
\usepackage{amsthm}

\numberwithin{equation}{section}

\newtheorem{theo}{Theorem}
\newtheorem{conj}{Conjecture}
\newtheorem{coro}{Corollary}
\newtheorem{prop}{Proposition}

\newtheorem{lem}{Lemma}

\theoremstyle{remark}
\newtheorem{Remark}{Remark}
\newtheorem{Remarks}[Remark]{Remarks}

\def\al{\alpha}
\def\be{\beta}
\def\ga{\gamma}

\def\ep{\varepsilon}

\def\la{\lambda}

\def\ph{\varphi}

\def\om{\omega}
\def\Ga{\Gamma}
\def\De{\Delta}
\def\Th{\Theta}

\def\Om{\Omega}
\def\({\left(}
\def\){\right)}
\def\[{\left[}
\def\]{\right]}

\def\fl#1{\left\lfloor#1\right\rfloor}

\def\lcm{\operatorname{lcm}}

\def\dd{\textup{d}}

\begin{document}

\title[]{On the integrality of the Taylor coefficients of mirror maps}

\author[]{C. Krattenthaler$^\dagger$ and T. Rivoal}
\date{\today}

\address{C. Krattenthaler, Fakult\"at f\"ur Mathematik, Universit\"at Wien,
Nordbergstra{\ss}e~15, A-1090 Vienna, Austria.
WWW: \tt http://www.mat.univie.ac.at/\~{}kratt.}

\address{T. Rivoal,
Institut Fourier,
CNRS UMR 5582, Universit{\'e} Grenoble 1,
100 rue des Maths, BP~74,
38402 Saint-Martin d'H{\`e}res cedex,
France.\newline
WWW: \tt http://www-fourier.ujf-grenoble.fr/\~{}rivoal.}

\thanks{$^\dagger$Research partially supported 
by the Austrian Science Foundation FWF, grants Z130-N13 and S9607-N13,
the latter in the framework of the National Research Network
``Analytic Combinatorics and Probabilistic Number Theory"}



\begin{abstract}
We show that the Taylor coefficients of the series
${\bf q}(z)=z\exp({\bf G}(z)/{\bf F}(z))$ are integers, where ${\bf F}(z)$ and 
${\bf G}(z)+\log(z) {\bf F}(z)$ 
are specific solutions of certain hypergeometric differential
equations with maximal unipotent monodromy at $z=0$. 
We also address the question of finding the largest integer $u$
such that the Taylor coefficients of $(z ^{-1}{\bf q}(z))^{1/u}$ are still
integers.
As consequences, we are able to prove numerous integrality
results for the Taylor coefficients of mirror maps 
of Calabi--Yau complete intersections in weighted projective spaces, 
which improve and refine previous results by Lian and Yau, and by Zudilin.
In particular, we prove the general ``integrality'' conjecture of 
Zudilin about these mirror maps. 
A further outcome of the present study is the determination of the
Dwork--Kontsevich sequence $(u_N)_{N\ge1}$, where $u_N$ is the
largest integer such that $q(z)^{1/u_N}$ is a series with integer
coefficients, where $q(z)=\exp(F(z)/G(z))$, 
$F(z)=\sum _{m=0} ^{\infty} (Nm)!\,z^m/m!^N$ and
$G(z)=\sum _{m=1} ^{\infty} (H_{Nm}-H_m)(Nm)!\,z^m/m!^N$, with $H_n$ denoting
the $n$-th harmonic number, conditional on the conjecture that there
are no prime number $p$ and integer $N$ such that the $p$-adic 
valuation of $H_N-1$ is strictly greater than~$3$.
\end{abstract}

\maketitle

\section{Introduction and statement of results}

\subsection{Mirror maps}
{\em Mirror maps} have appeared quite recently in mathematics and
phys\-ics. Indeed, 
the term ``mirror map" was coined in the late  
1980s
by physicists whose research in string theory led them to discover
deep facts in algebraic geometry 
(e.g., pair of Calabi--Yau manifolds in ``mirror,'' enumeration of 
rational curves on one of these Calabi--Yau manifolds in dimension $3$).
In a sense, mirror maps can be viewed as higher dimensional generalisations of 
certain classical modular forms, which naturally appear in low dimensions 
(see some examples below).

The purpose of the present article is to prove rather sharp
integrality assertions for the Taylor coefficients of 
certain mirror maps coming from
hypergeometric differential equations, which are Picard--Fuchs equations of 
suitable one parameter families of Calabi--Yau 
complete intersections in weighted projective spaces.
The corresponding results
(see Theorems~\ref{thm:1}--\ref{thm:4}) encompass 
integrality results on these mirror maps which exist in the
literature, improving and refining them in numerous cases.

Before getting deeper into the subject, it is 
beneficial to define at this point a special case of 
a mirror map
(see just after~\eqref{eq:rajout6} below), which will be studied in great 
detail in the present paper. For a real number $x>0$ and 
an integer $m\ge 0$, set 
$$
H(x,m):=\sum_{n=0}^{m-1} \frac{1}{x+n},
$$
where, by convention, the empty sum is $0$. When $x=1$, $H(1,m)$ is simply 
the $m$-th harmonic number $H_m$. 
For integers $k\ge 1$ and $N\ge 1$, let us define the power 
series
\begin{equation*}
F_N(z):=\sum_{m=0}^{\infty} \frac{(Nm)!^k}{m!^{kN}} \,z^m
\end{equation*}
and
\begin{equation*}
G_N(z):=k\sum_{m=1}^{\infty} \frac{(Nm)!^k}{m!^{kN}}
\bigg(\sum_{j=1}^{N-1}H(j/N,m)-(N-1)H(1,m)\bigg)\,z^m,
\end{equation*}
which converge for $\vert z\vert <1/N^{kN}.$ 
(\footnote{Most of the time, the dependence on $k$ will not be 
explicit in our notation for the sake of better readability. 
Furthermore, for the results of the
present paper, it is sufficient to consider the series as 
{\it formal\/} power series.})
The functions $F_N(z)$ 
and $G_N(z)+\log(z)F_N(z)$ are solutions of the same hypergeometric 
differential equation, 
which is a special case of~\eqref{eq:equadiff} below. 
(It is of maximal unipotent monodromy, 
i.e., $F_N(N^{-kN}z)$ is a hypergeometric function with only $1$'s 
as lower parameters. 
Its other solutions around $z=0$ 
can then be obtained by Frobenius' method, see~\cite{yoshida}.)
The function 
\begin{equation}\label{eq:rajout6}
q_N(z):=z\exp(G_N(z)/F_N(z))\in z\mathbb{Q}[[z]]
\end{equation}
is usually called the {\em canonical coordinate}, and its compositional
inverse $z_N(q)$ is  
the prototype of a mirror map. In this paper, 
by abuse of terminology, we will 
also use the term ``mirror map" for any  
canonical coordinate.~(\footnote{Canonical coordinates and mirror
maps have distinct  
geometric meanings. However, 
in the number-theoretic study undertaken in the present paper,
they play strictly the same role, because $(z^{-1}q(z))^{1/\tau}\in 
1+z\mathbb{Z}[[z]]$ for some integer  
$\tau$ implies that $(q^{-1}z(q))^{1/\tau}\in 1+q\mathbb{Z}[[q]]$, and
conversely. (See~\cite[Introduction]{lianyau1}.)})  

\medskip

The case $N=1$ is trivial because $G_1(z)=0$ for any $k$. 
When $N=2$ and $k=1$, we have $F_2(z)=1/\sqrt{1-4z}$ and
$q_2(z)=4z/(1+\sqrt{1-4z})^2$. 
For small values of $N$ and $k$, 
modular forms rapidly enter the picture. Indeed, 
for $N=2$ and  
$k=2$, the compositional inverse of $q_2(z)$ is equal to
$\lambda(q)/16:=q\prod_{n=1}^{\infty}\big((1+q^{2n})/(1+q^{2n-1})\big)^8$,
which is 
a modular function of the variable $\tau$, with $\Im(\tau)>0$, defined
by $q=\exp(2i \pi \tau)$. Furthermore, when $N=2$ and   
$k=3$, the compositional inverse of $q_2(z)$ is equal to
$\lambda(q)(1-\lambda(q))/16$. See for example the 
discussion in~\cite[pp.~111--113]{andre} (and also for the importance
of such facts in Diophantine approximation),  
and see~\cite[Sec.~3]{lianyau} for a discussion of the modularity of the 
case $N=4, k=1$. In fact, Doran~\cite{doran} proved a result
that enables
one to describe all mirror maps of modular origin. 

The most famous non-modular example of
a mirror map is $q_5(z)$ (when $N=5, k=1$), which was used  
in the epoch-making paper 
by the physicists Candelas et al.~\cite{candelas}. 
Without going into details (see also~\cite{morrisonjams, pandha,VoisAA} 
and the references therein), let us give a short 
explanation of the importance of $q_5(z)$. 
Starting from the family ${\mathbf M}$ of quintic hypersurfaces 
in $\mathbb{P}^4(\mathbb{C})$ defined by $\sum_{k=1}^5x_k^5-5z
\prod_{k=1}^5x_k=0$ ($z$ being a complex parameter),  
Candelas et al.\ naturally associate another family $\mathbf{W}$ of
manifolds (the ``mirror of ${\mathbf M}$'') which turn out to be
Calabi--Yau. To $\mathbf{W}$, one can   
naturally associate a vector of periods (depending on $z$) which are
solutions of the same differential equation (namely, the Picard--Fuchs
equation of $\mathbf{W}$). This equation is simply the hypergeometric
differential equation satisfied by $F_5(z)$ and $G_5(z)+\log(z)F_5(z)$
(case $N=5$, $k=1$ above). Then they observed the non-trivial
property that the Taylor 
coefficients of $q_5(z)$ are integers.  
Furthermore, let us define the {\em Yukawa coupling}~(\footnote{The 
Yukawa coupling is a geometric object whose definition 
in a specific situation can be found
in~\cite[Definition 4.5.2]{batstrat}.  
In the present case, it can be computed as in~\eqref{eq:yukawa}.})
\begin{equation} \label{eq:yukawa}
K(q):= \frac{5}{1-5^5z_5(q)}\cdot \frac{1}{F_5(z_5(q))^2}\cdot
\bigg(\frac{qz_5'(q)}{z_5(q)}\bigg)^3 \in \mathbb{Q}[[q]], 
\end{equation}
where $z_5(q)$ is the compositional inverse of $q_5(z)$, and write it as 
$
K(q) = 5+ \sum_{d=1}^{\infty} k_d \,\frac{q^d}{1-q^d},
$
which is formally possible. Candelas et al.\ observed that the {\em
instanton number} $n_d:=k_d/d^3$ is an integer for all  
$d\ge 1$, which is already a non-trivial fact, but that furthermore
$n_d$ seems to be the number of rational curves of degree $d$ lying  
on the initial quintic ${\mathbf M}$, 
thereby providing an effective algorithm to compute these numbers. 
These striking observations generated much
interest amongst algebraic geometers, and this culminated 
in the work
of Givental~\cite{givental} and Lian et al.~\cite{lianliuyau} where it is proved that, 
if for a given $d$ the curves of degree $d$ are all rigid, then there are $n_d$ of them; 
see also the discussion in~\cite[p.~49]{VoisAA}. In fact, the coincidence was proved 
to be true for $d\le 9$, and the first difference occurs at $d=10$ (see~\cite{cott}).
(However, we do not address questions related to Yukawa couplings
in the present paper.)

\subsection{Integrality of mirror maps}
Deep results can also be proved without any appeal to algebraic geometry. 
Indeed, using $p$-adic methods, 
Lian and Yau were the first to prove that $q_5(z)$ has integral Taylor
coefficients. In fact,  
in~\cite[Sec.~5, Theorem~5.5]{lianyau}, they proved that
$q_N(z)\in z\mathbb{Z}[[z]]$ for $k=1$ and any $N$ which is a prime
number.
Zudilin~\cite[Theorem~3]{zud} extended their result to any $k\ge 1$
and any $N$ which is a prime power and made the following 
conjecture, which is also implicit in the cited articles of Lian and
Yau, for example at the end of the introduction of~\cite{lianyau2}. 
In fact, the conjecture probably belongs to the folklore 
of mirror symmetry theory, and 
it seems that it had been left open for $N$ not a prime power.

\begin{conj}
\label{conj:zudilin1}
For any integers $k\ge 1$ and $N\ge 1$, we have $q_N(z)\in z\mathbb{Z}[[z]]$.
\end{conj}

Such integrality questions for mirror maps and Yukawa couplings 
undoubtedly remain an important question for algebraic geometers, as
is witnessed by the very recent  
preprint~\cite{volog} (which is an elaborate version of
\cite{KoSVAA}) on this subject. 
The mirror maps that are considered in that paper comprise ours.
When both approaches apply simultaneously, our  
results in Theorems~\ref{thm:1}--\ref{thm:4} 
are stronger than Theorem~2 in \cite[Sec.~1.3]{volog}.
Indeed, we prove  
that certain mirror maps have integral Taylor
coefficients, while in \cite{volog} the weaker statement is proved that 
mirror maps have Taylor coefficients in $\mathbb{Z}[1/n]$, where
$n$ is an integer parameter of geometric origin which is at  
least $2$ (by assumption iii) just before Theorem~2 in \cite{volog}).
On the other hand, the range of applicability of \cite[Theorem~2]{volog}
is much wider than ours. 
It is interesting to note that our approach as well
as the one in \cite{volog} are heavily based on $p$-adic analysis, so that
both are clearly close in spirit, although the exact methods that are
applied are different.

\medskip

In~\cite{lianyau2}, 
Lian and Yau strengthened their result
from \cite{lianyau} by proving an observation made by
physicists:   
\begin{equation}\label{eq:raflianyau}
\big(z^{-1}q_N(z)\big)^{1/N}\in \mathbb{Z}[[z]]
\end{equation}
for $k=1$ and any prime $N$.

Our original goal was to settle 
Conjecture~\ref{conj:zudilin1} and to prove \eqref{eq:raflianyau} 
for arbitrary $k$ and $N$. 
In the present paper, we shall accomplish much more: 
we establish refinements of the above integrality assertions which enable
us to prove results that are even stronger and more general than 
Conjecture~\ref{conj:zudilin1} or \eqref{eq:raflianyau}.
In fact, our refinements will go into two different, only barely 
overlapping directions. 
One direction is inspired by a conjecture in  
Zudilin's paper \cite{zud}, while the other seems to be entirely new.
In the remainder of this introductory section, we describe these two
directions, and we present our results. Their proofs will then be
given in the subsequent sections.

\subsection{Refinements of \eqref{eq:raflianyau}, 
part~I} \label{subsec:1}
We describe the second direction of refinement of \eqref{eq:raflianyau}
first. Starting point is the observation that
$$
\sum_{j=1}^{N-1}H(j/N,m) -(N-1)H(1,m)= N H_{Nm}-NH_m.
$$
Hence, with 
$$
G_{L,N}(z):=\sum_{m=1}^{\infty} H_{Lm} \,\frac{(Nm)!^k}{m!^{kN}} \,z^m
$$
and $q_{L,N}(z):=\exp(G_{L,N}(z)/F_N(z))$, 
we have 
\begin{equation}\label{eq:truemap}
q_N(z)=zq_{N,N}(z)^{kN}q_{1,N}(z)^{-kN}.
\end{equation} 
We are ready to state
our first refinement of Conjecture~\ref{conj:zudilin1} and 
\eqref{eq:raflianyau}.

\begin{theo}\label{thm:1}
For any integers $k, N\ge 1$ and $L\in \{1, 2, \ldots, N\}$, we have 
$q_{L,N}(z) \in\mathbb{Z}[[z]].$
\end{theo}

Clearly, 
due to \eqref{eq:truemap},
Theorem~\ref{thm:1} implies Conjecture~\ref{conj:zudilin1}
and \eqref{eq:raflianyau}.

However, much more is true. For any
sequence $\mathbf{N}=(N_1, \ldots, N_k)$
of positive integers (the $N_j$'s are not necessarily distinct)  
and any integer $L\ge 1$, let us define the power series
\begin{align*}
F_{\mathbf{N}}(z)&:= \sum_{m=0}^{\infty} 
\bigg(
\prod_{j=1}^k
\frac{(N_jm)!}{m!^{N_j}}  
\bigg) 
z^m, \\
G_{L,\mathbf{N}}(z)&:= \sum_{m=1}^{\infty} H_{Lm}\bigg(\prod_{j=1}^k
\frac{(N_jm)!}{m!^{N_j}}  \bigg) z^m,
\end{align*}
and the function 
$q_{L,\mathbf{N}}(z):=\exp\big(G_{L,\mathbf{N}}(z)/F_{\mathbf{N}}(z)\big).$ 
The series $F_{\mathbf{N}}(z)$ and
$G_{\mathbf{N}}(z)+\break \log(z)F_{\mathbf{N}}(z)$, 
where $G_{\mathbf{N}}(z)$ is a suitable linear combination 
of the series $G_{L,\mathbf{N}}(z)$
(for different $L$'s), are solutions of a hypergeometric
differential equation with  
maximal unipotent monodromy (of type~\eqref{eq:equadiff} below). 
That differential equation is the Picard--Fuchs equation 
of a one parameter family of mirror manifolds 
$V'$ of a complete intersection $V$ of 
$k$ hypersurfaces $V_1, \ldots, V_k$ of degrees $N_1, \ldots, N_k$ in
$\mathbb{P}^{d+k}(\mathbb{C})$: $V$ is a family of Calabi--Yau manifolds  
if one chooses
$d$ equal to $\sum _{j=1} ^{k}N_j-k-1$. The mirrors $V'$ are
explicitly constructed in~\cite[Sec.~5.2]{batstrat}. 

Let us define $M_{\mathbf{N}}=\prod_{i=1}^k N_i!$ and, for $L\ge
1$, $V_{L,\mathbf{N}}$ as the largest integer 
such that $q_{L,\mathbf{N}}(z)^{1/V_{L,\mathbf{N}}} \in 
\mathbb{Z}[[z]]$.~(\footnote{Let $q(z)$ be a given power 
series in $\mathbb{Z}[[z]]$, and let
$V$ be the largest integer with the property that $q(z)^{1/V}\in \mathbb
Z[[z]]$. Then $V$ carries complete information about {\it all\/}
integers with that property: namely,
the set of integers $U$ with $q(z)^{1/U}\in \mathbb
Z[[z]]$ consists of all divisors of $V$. Indeed, it is clear that all
divisors of $V$ belong to this set. Moreover, 
if $U_1$ and $U_2$ belong to this set, then also $\lcm(U_1,U_2)$ does
(cf.\ \cite[Lemma~5]{HeRSAA} for a simple proof 
based on B\'ezout's lemma).})    
While we are not able to determine $V_{L,\mathbf{N}}$ precisely,
we shall prove the following result which, as 
should become clear from the discussion 
in Section~\ref{sec:DworkKont}, comes relatively close.

\begin{theo} \label{thm:2}
Let $\Th_L:=L!/\gcd(L!, L!\,H_L)$ be the denominator of $H_L$
when written as a reduced fraction.
Then, for any integers $N_1, \ldots, N_k\ge 1$ and 
$L\in \{1, 2, \ldots, \max(N_1, \ldots, N_k)\}$, we have 
$q_{L,\mathbf{N}}(z)^{\frac{\Th_L}{M_{\mathbf{N}}}} \in\mathbb{Z}[[z]].$
\end{theo}

\begin{Remarks} \label{rem:1}
(a)
For any integer $s \ge 1$, 
we have $q_{L,\mathbf{N}}(z)^{1/s}=1+s^{-1}M_{{\bf N}}H_L z + \mathcal{O}(z^2)$, 
and hence  Theorem~\ref{thm:2} is optimal when $L=1$. This is not necessarily the case for 
other values of $L$ and, in particular, Theorem~\ref{thm:2} can sometimes be improved 
when $L=N_1=\cdots =N_k$ (see Theorem~\ref{thm:3} below).

(b) Considering the case $k=1$ and $L=N$ in Theorem~\ref{thm:2},
the first few values of $M_{(N)}/\Th_N=\gcd(N!, N!\,H_N)$ (for 
$N\ge 1$) are
\begin{equation}\label{eq:numericalvalues}
1, 1, 1, 2, 2, 36, 36, 144, 144, 1440, 1440, 17280, 17280, 241920,
3628800, 29030400,\dots  
\end{equation}
In the On-Line Encyclopedia of Integer Sequences~\cite{oeis}, 
this sequence is entry {\tt A056612}.
\end{Remarks}

By 
forming a suitable product of the functions $q_{L,\mathbf{N}}(z)$,
Theorem~\ref{thm:2}
implies the integrality of Taylor 
coefficients of the corresponding mirror maps $q_{\mathbf{N}}(z)$ of
the mirror pair $(V,V')$.  
According to~\cite[Prop. 5.1.2]{batstrat}, the Yukawa 
coupling in this case is equal to 
$$
K(q):=\frac{N_1N_2\cdots N_k}{1-\lambda z_{\mathbf{N}}(q)}\cdot
\frac{1}{F_{\mathbf{N}}(z_{\mathbf{N}}(q))^2}\cdot 
\bigg(\frac{qz_{\mathbf{N}}'(q) }{z_{\mathbf{N}}(q)} \bigg)^{d},
$$
where $\lambda=\prod_{j=1}^k N_j^{N_j}$ and $z_{\mathbf{N}}(q)$ is 
the compositional inverse of $q_{\mathbf{N}}(z)$. When 
$d=3=\sum_{j=1} ^{k}N_j-k-1$, 
the formal expansion
$K(q)=K(0)+\sum_{m=1}^{\infty}k_m \frac{q^m}{1-q^m}$ enables
one to count rational curves on the Calabi--Yau threefold $V$, 
at least for small values of the degree of the curves.

An outline of the proof of Theorem~\ref{thm:2} is given in
Section~\ref{sec:2}, with details being filled in in later sections.
Since we shall refer to it below, 
we remark that an alternative way to define
the integer $\Th_L$ is via 
\begin{equation} \label{eq:ThL} 
\Th_L=
{\prod _{p\le L}^{}}p^{-\min\{0,v_p(H_L)\}},
\end{equation}
where $v_p(\alpha)$ denotes the $p$-adic valuation of $\alpha$. 
Here and in the sequel of the article, 
the letter $p$ will always represent a prime number.

Due to~\eqref{eq:truemap}, Theorem~\ref{thm:2} has the following
consequence for the original 
mirror map $q_N(z)$, thus improving significantly upon
\eqref{eq:raflianyau}.

\begin{coro} \label{coro:1}
For all integers $k\ge 1$ and $N\ge 1$, we have 
$$
\big(z^{-1}q_N(z)\big)^{\frac{\Th_N}{N!^k kN}} \in \mathbb{Z}[[z]],
$$
where $q_N(z)$ is the mirror map in~\eqref{eq:rajout6}.
\end{coro}

In particular, in the emblematic case of the mirror map $q_5(z)$ of
the quintic (case $N=5$, $k=1$), we obtain that  
$\big(z^{-1}q_5(z)\big)^{1/10} \in \mathbb{Z}[[z]]$, which improves
on~\eqref{eq:raflianyau} by a factor of $2$. 

\medskip
The next theorem presents the improvement of Theorem~\ref{thm:2} for
the case $L=N_1=\dots=N_k$ which was announced in 
Remark~\ref{rem:1}(a) above.

\begin{theo} \label{thm:3}
Let $N$ be a positive integer, 
$\mathbf N=(N,N,\dots,N)$, with $k$ occurrences of $N$,
and let $\Xi_1=1$, $\Xi_7=1/140$, and, for $N\notin\{1,7\}$,
\begin{equation} \label{eq:Xi} 
\Xi_N:=
{\prod _{p\le N}
^{}}p^{\min\{2+\xi(p,N),v_p(H_N)\}},
\end{equation}
where $\xi(p,N)=1$ if $p$ is a Wolstenholme prime {\em(}i.e.,
a prime $p$ for which $v_p(H_{p-1})\ge3$
{\em(}\footnote{Presently, only two such primes are known, namely 
$16843$ and $2124679$, and it is unknown whether there are infinitely
many Wolstenholme primes or not.}{\em))} or $N$ is divisible by
$p$, and $\xi(p,N)=0$ otherwise. 
Then 
$q_{N,\mathbf N}(z)^{\frac{1}{\Xi_{N}N!^k}} \in\mathbb{Z}[[z]].$
\end{theo}

\begin{Remarks} \label{rem:Xi7}
For better comprehension, we discuss the meaning of the
statement of Theorem~\ref{thm:3} and its implications;
in particular, we address some fine points of the definition of $\Xi_N$.

\smallskip
(a) The case of $N=1$ is trivial since 
$q_{1,(1,\dots,1)}(z)=1/(1-z)$. Furthermore, we have
$$\Xi_7=\frac {1} {140}=2^{v_2(H_7)}5^{v_5(H_7)}7^{v_7(H_7)},$$
which differs by a factor of $3$ from the right-hand side of
\eqref{eq:Xi} with $N=7$ 
(since $v_3(H_7)=v_3(\frac {363} {140})=1$).

\smallskip
(b) Since 
$q_{N,\mathbf N}(z)=1+H_N N!^kz+\mathcal O(z^2)$, 
it is clear that $q_{N,\mathbf
N}(z)^{1/(p^{v_p(H_N)+1}N!^k)}\notin \mathbb{Z}[[z]]$, so that the 
exponent of $p$ in the prime factorisation of
$V_{N,\mathbf N}$ can be at most $v_p(H_N N!^k)$.
In Theorem~\ref{thm:3}, this theoretically maximal exponent is
further cut down. First of all, the number $\Xi_N$ appearing there
contains no prime factor $p>N$. Moreover, for primes $p$ with $p\le N$
and $v_p(H_N)\ge3$, the definition of $\Xi_N$ cuts the theoretically
maximal exponent $v_p(H_NN!^k)$
of $p$ down to $2+v_p(N!^k)$ respectively
$3+v_p(N!^k)$, depending on whether $\xi(p,N)=0$ or $\xi(p,N)=1$.
In items~(c)--(e) below, we address the question of how serious this cut is
expected to be.

\smallskip
(c) Clearly, the minimum appearing in the exponent of $p$ in the definition
\eqref{eq:Xi} of $\Xi_N$ is $v_p(H_N)$ as long as $v_p(H_N)\le 2$. 
In other words, the exponent of $p$ in the prime factorisation of
$\Xi_N$ depends largely on the $p$-adic behaviour of $H_N$. 
An extensive discussion of this topic, with many interesting results,
can be found in \cite{boyd}. We have as well computed a table of
harmonic numbers $H_N$ up to 
$N=
1000000$.~(\footnote{The summary of the table is available at
{\tt
http://www.mat.univie.ac.at/\~{}kratt/artikel/H.html}.})
Indeed, the data suggest that
pairs $(p,N)$ with $p$ prime, $p\le N$, and $v_p(H_N)\ge 3$ 
are not very frequent. 
More precisely, so far only five examples are
known with $v_p(H_N)=3$: four for $p=11$, with $N=848,9338,10583$,
and $3546471722268916272$, and one for $p=83$
with
\begin{align} \notag
&\hbox{\small
$N=79781079199360090066989143814676572961528399477699516786377994370\backslash$}
\\
&\kern2cm
   \hbox{\small $78839681692157676915245857235055200779421409821643691818$}
\label{eq:boyd} 
\end{align}
(see \cite[p.~289]{boyd}; the value of $N$ in \eqref{eq:boyd}, not
printed in \cite{boyd}, was kindly communicated to us by David Boyd). 
There is no example known with $v_p(H_N)\ge4$.
It is, in fact, conjectured that no $p$ and $N$ exist with
$v_p(H_N)\ge4$. Some evidence for this conjecture (beyond mere
computation) can be found in \cite{boyd}.

\smallskip
(d) Since in all the five examples for which
$v_p(H_N)=3$ we neither have $p\mid N$ (the gigantic number in
\eqref{eq:boyd} is congruent to $42$ mod $83$) nor that the prime $p$ is a
Wolstenholme prime, the exponent of $p$ in the prime factorisation of
$\Xi_N$ in these cases is $2$ instead of $v_{p}(H_N)=3$. 

\smallskip
(e) On the other hand, should there be a prime $p$ and an integer $N$
with $p\le N$, $v_p(H_N)\ge3$, $p$ a Wolstenholme prime or $p\mid
N$, then the exponent of $p$ in the prime factorisation of
$\Xi_N$ would be $3$. However, no such examples are known. We
conjecture that there are no such pairs $(p,N)$. If this conjecture
should turn out to be true, then, given $N\notin\{1,7\}$,
the definition of $\Xi_N$ in \eqref{eq:Xi} could be simplified to
\begin{equation} \label{eq:Xisimpl} 
\Xi_N:=
{\prod _{p\le N}
^{}}p^{\min\{2,v_p(H_N)\}}.
\end{equation}
\end{Remarks}

In view of \eqref{eq:ThL} and the fact that $M_{\mathbf N}=N!^k$ for
the vector $\mathbf N$ in Theorem~\ref{thm:3}, this theorem
improves upon Theorem~\ref{thm:2} in the case $L=N$.
Namely, Theorem~\ref{thm:3} is {\it always} at least as strong as 
Theorem~\ref{thm:2}, and it is {\it strictly} stronger if $N\ne7$ and
$v_p(H_N)\ge 1$ for some prime $p$ less than or equal to $N$. 
Indeed, the smallest $N\ne 7$ with that property is $N=20$, in which case 
$v_5(H_{20})=1$.

We remark that strengthenings of
Theorem~\ref{thm:2} in the spirit of Theorem~\ref{thm:3} 
for more general choices of the parameters 
can also be obtained by our techniques but are omitted here. 

We outline the proof of Theorem~\ref{thm:3} in Section~\ref{sec:Xi},
with the details being filled in 
in the subsequent Section~\ref{sec:aux}.
As we explain in Section~\ref{sec:DworkKont}, 
we conjecture that Theorem~\ref{thm:3} cannot be improved if $k=1$,
that is, that for $k=1$ the largest integer $t_N$ such that
$q_{N,N}(z)^{1/t_N}\in\mathbb Z[[z]]$ is exactly $\Xi_N N!$. 
Propositions~\ref{prop:p>N} and \ref{prop:vp=3} 
in Section~\ref{sec:DworkKont} show that 
this conjecture would immediately follow if one could prove
the conjecture from Remarks~\ref{rem:Xi7}(c) above that there are no 
primes $p$ and integers $N$ with $v_p(H_N)\ge 4$. 

\medskip
Even if the series $q_{N,N}(z)$ appears in the identity
\eqref{eq:truemap}, relating our mirror map $q_N(z)$ to the series
$q_{L,N}(z)$ (with $L=1$ and $L=N$), Theorem~\ref{thm:3} does not
imply an improvement over Corollary~\ref{coro:1}, for the
coefficient of $z$ in $(z^{-1}q_N(z))^{\Th_N/(pN!^kkN)}$ is equal
to $\frac {\Th_N} {p}(H_N-1)$, and, thus,
it will not be integral for primes $p$
with $v_p(H_N)>0$. Still, there is an improvement of
Corollary~\ref{coro:1} in the spirit of Theorem~\ref{thm:3}. It
involves primes $p$ with $v_p(H_N-1)>0$ instead.

\begin{theo} \label{thm:3a}
Let $N$ be a positive integer with $N\ge2$, and let
\begin{equation} \label{eq:Om} 
\Om_N:=
{\prod _{p\le N}
^{}}p^{\min\{2+\om(p,N),v_p(H_N-1)\}},
\end{equation}
where $\om(p,N)=1$ if $p$ is a Wolstenholme prime 
or $N\equiv\pm1$~{\em mod}~$p$, and $\om(p,N)=0$ otherwise. 
Then 
$\big(z^{-1}q_{N}(z)\big)^{\frac{1}{\Om_{N}N!^kkN}} \in\mathbb{Z}[[z]].$
\end{theo}
\begin{Remarks} \label{rem:Om}
Also here, some remarks are in order to get a better understanding of
the above theorem.

\smallskip
(a)
If $v_p(H_N)< 0$, then $v_p(H_N)=v_p(H_N-1)$. Hence, 
differences in the prime factorisations of $\Xi_N$ and $\Om_N$ 
can only arise for primes $p$ with $v_p(H_N)\ge0$.

\smallskip
(b) Since $z^{-1}q_{N}(z)=1+(H_N-1)N!^kkNz+\mathcal O(z^2)$, it is clear that 
$$\big(z^{-1}q_{N}(z)\big)^{1/(p^{v_p(H_N-1)+1}N!^kkN)}\notin
\mathbb{Z}[[z]].$$
If $\tilde V_N$ denotes the largest integer such that 
$\big(z^{-1}q_N(z)\big)^{1/\tilde V_NkN}$ is a series with integer coefficients,
the exponent of $p$ in $\tilde V_N$ can be at most 
$v_p\big((H_N-1)N!^k\big)$.
In Theorem~\ref{thm:3a}, this theoretically maximal exponent is
further cut down. Namely, as is the case for $\Xi_N$, 
the number $\Om_N$ in \eqref{eq:Om}
contains no prime factor $p>N$. Moreover, for primes $p$ with $p\le N$
and $v_p(H_N-1)\ge3$, the definition of $\Om_N$ cuts the theoretically
maximal exponent $v_p\big((H_N-1)N!^k\big)$
of $p$ down to $2+v_p(N!^k)$ respectively
$3+v_p(N!^k)$, depending on whether $\om(p,N)=0$ or $\om(p,N)=1$.

\smallskip
(c)
Concerning the question whether there are at all primes $p$ and integers $N$
with high values of $v_p(H_N-1)$, 
we are not aware of any corresponding literature. 
Our table of harmonic numbers $H_N$ mentioned in Remarks~\ref{rem:Xi7}(c)
does not contain any pair $(p,N)$ with 
$v_p(H_N-1)\ge3$.~(\footnote{The summary of the 
corresponding table, containing pairs $(p,N)$ with $p\le N$ and
$v_p(H_N-1)>0$, is available at
{\tt
http://www.mat.univie.ac.at/\~{}kratt/artikel/H1.html}.})
In ``analogy" to the conjecture mentioned in
Remarks~\ref{rem:Xi7}(c), we conjecture that no $p$ and $N$ exist with
$v_p(H_N-1)\ge4$. It may even be true that there are no $p$ and
$N$ with $v_p(H_N-1)\ge3$, in which case 
the definition of $\Om_N$ in \eqref{eq:Om} could be simplified to
\begin{equation} \label{eq:Omsimpl} 
\Om_N:=
{\prod _{p\le N}
^{}}p^{\min\{2,v_p(H_N-1)\}}.
\end{equation}
\end{Remarks}

In view of \eqref{eq:ThL}, Theorem~\ref{thm:3a} 
improves upon Corollary~\ref{coro:1}.
Namely, Theorem~\ref{thm:3a} is {\it always} at least as strong as 
Corollary~\ref{coro:1}, and it is {\it strictly} stronger if 
$v_p(H_N-1)\ge 1$ for some prime $p$ less than or equal to $N$. 
The smallest $N$ with that property is $N=21$, in which case 
$v_5(H_{21}-1)=1$.

We sketch the proof of Theorem~\ref{thm:3a} in Section~\ref{sec:Om}.
We also explain in that section that
we conjecture that Theorem~\ref{thm:3a} with $k=1$ is optimal,
that is, that for $k=1$ the largest integer $u_N$ such that
$\big(z^{-1}q_{N}(z)\big)^{\frac{1}{Nu_N}} \in\mathbb{Z}[[z]]$ 
is exactly $\Om_N N!$. 
Propositions~\ref{prop:p>N2} and \ref{prop:vp=32} 
in Section~\ref{sec:Om} show that 
this conjecture would immediately follow if one could prove
the conjecture from Remarks~\ref{rem:Om}(c) above that there are no 
primes $p$ and integers $N$ with $v_p(H_N-1)\ge 4$. 
As a matter of fact, the sequence $(u_{2N})_{N\ge1}$ appears also in the 
On-Line Encyclopedia of Integer Sequences~\cite{oeis},
as sequence {\tt A007757}, contributed by R.~E.~Borcherds under 
the denomination ``Dwork--Kontsevich sequence'' around 1995,
without any reference or explicit formula for it, however.
(\footnote{In private communication, both, Borcherds and Kontsevich
could not remember where exactly this sequence and its denomination came 
from.}) 

\subsection{Refinements of \eqref{eq:raflianyau}, part~II}
We now move on to describe the other direction of refinement of 
Conjecture~\ref{conj:zudilin1} and \eqref{eq:raflianyau},
inspired by Zudilin's paper~\cite{zud}.
For any given integer $N\ge 1$, let $r_1, r_2, \ldots, r_d$ denote 
the integers in $\{1, 2,\ldots, N\}$ which are coprime to $N$. It is
well-known that $d=\varphi(N)$, Euler's totient function, which is given by  
$
\varphi(N) = N 
{\prod_{p\mid N}} \big(1-\frac{1}{p}\big).
$
Set
$
C_N := N^{\varphi(N)}
{\prod_{p\mid N}} p^{\varphi(N)/(p-1)},
$
which is an integer because $p-1$ divides $\varphi(N)$ for any prime $p$
dividing $N.$ 
Let us also define the Pochhammer symbol $(\al)_m$ for complex numbers 
$\al$ and non-negative integers $m$ by
$(\al)_m:=\al(\al+1)\cdots (\al+m-1)$ if $m\ge1$ and $(\al)_0:=1$.
It can be proved (see~\cite[Lemma~1]{zud}, 
or \eqref{eq:Bzudilin} together with Lemma~\ref{lem:10a}$(iii)$)  
that, for any integer $m\ge 0$,  
\begin{equation}\label{eq:definitionBNgras}
\mathbf{B}_N(m) := C_N^m \prod_{j=1}^{\varphi(N)} 
\frac{(r_j/N)_m}{m!}
\end{equation}
is an integer. 
We will as well use 
$\mathbf{B}_{\mathbf{N}}(m):=\prod_{j=1}^k \mathbf{B}_{N_j}(m)$ for  
vectors $\mathbf{N}=(N_1, \ldots, N_k)$ of positive integers.
Zudilin also 
established another representation for $\mathbf{B}_N(m)$ (see \cite[Lemma
4]{zud}). 
Namely, let $p_1, p_2, \dots, p_\ell$ denote the distinct prime factors   
of $N$, and let us define the vectors of 
integers
\begin{equation}
(\alpha_j)_{j=1, \ldots, \mu}  := \left( N, \frac{N}{p_{j_1}p_{j_2}},
\frac{N}{p_{j_1}p_{j_2}p_{j_3}p_{j_4}}, \ldots 
\right)_{1\le j_1<j_2<\dots \le \ell } \label{eq:aj}
\end{equation}
and
\begin{equation}
(\beta_j)_{j=1, \ldots, \eta} := \left(
\frac{N}{p_{j_1}}, \frac{N}{p_{j_1}p_{j_2}p_{j_3}}, \ldots,  1,1, \dots, 1 
\right)_{1\le j_1<j_2<\dots \le \ell }, \label{eq:bj}
\end{equation}
where $\al_1+\al_2+\dots + \al_\mu=\be_1+\be_2+\dots +\be_\eta.$ Then we 
have 
\begin{equation}\label{eq:Bzudilin}
\mathbf{B}_N(m) = \frac{\prod_{j=1}^\mu (\alpha_j
m)!}{\prod_{j=1}^\eta (\beta_j m)!}. 
\end{equation}
For example, we have 
$$
\mathbf{B}_4(m) = \frac{(4m)!}{(2m)!\,m!^2}, \; \mathbf{B}_6(m) =
\frac{(6m)!}{(3m)!\,(2m)!\,m!},  
\; \mathbf{B}_{30}(m) =
\frac{(30m)!\,(5m)!\,(3m)!\,(2m)!}{(15m)!\,(10m)!\,(6m)!\,m!^9}. 
$$

We now define the functions 
\begin{align} \notag
\mathbf{H}_N(m)&:= \sum_{j=1}^{\varphi(N)} H(r_j/N,m) - \varphi(N)H(1,m),\\
\label{eq:definitionFgrasNgras}
\mathbf{F}_{\mathbf{N}}(z) &:= \sum_{m=0}^{\infty} 
\bigg(
\prod_{j=1}^k
\mathbf{B}_{N_j}(m)
\bigg)
z^m,  
\end{align}
and
$$
\mathbf{G}_{\mathbf{N}}(z) := \sum_{m=1}^{\infty} \bigg(\sum_{j=1}^k
\mathbf{H}_{N_j}(m) \bigg)\bigg( \prod_{j=1}^k \mathbf{B}_{N_j}(m)\bigg)
z^m.  
$$

We can now state Zudilin's conjecture 
from \cite[p.~605]{zud}.
\begin{conj}[\sc Zudilin]\label{conj:zudilin2}
For any positive integers $N_1, N_2, \dots, N_k$, we have 
$
{\bf q}_{\mathbf{N}} (z):= z 
\exp(\mathbf{G}_{\mathbf{N}}(z)/\mathbf{F}_{\mathbf{N}}(z)) \in
z\mathbb{Z}[[z]]. 
$
\end{conj}
As explained by Zudilin, this conjecture implies 
Conjecture~\ref{conj:zudilin1} because any function 
$q_{N}(z)$ is equal to $\mathbf{q}_{\mathbf{N}}(z)$ for a suitable
$\mathbf{N}$. More precisely, for a given $N\ge 1$,  
the vector $\mathbf{N}$ 
to choose consists of all the positive divisors of $N$ but $1.$
For example, we have $q_9(z)={\mathbf q}_{(9,3)}(z)$, 
$q_{12}(z)={\mathbf q}_{(12,6,4,3,2)}(z)$
and $q_{35}(z)={\mathbf q}_{(7,5)}(z)$. However,
Conjecture~\ref{conj:zudilin2} does not 
imply any of the previous theorems. Zudilin proved that his conjecture holds 
under the condition that if a prime number divides $N_1N_2\cdots
N_k$ then it also divides each $N_j$.  
This applies to the function ${\mathbf q}_{(9,3)}(z)$, but neither to 
${\mathbf q}_{(12,6,4,3,2)}(z)$ nor to ${\mathbf q}_{(7,5)}(z)$,
for example. 

We claim that Conjecture~\ref{conj:zudilin2} follows from the theorem
below. For the statement of the theorem, 
for an integer $L\ge 1$, we need to define  
$$
\mathbf{G}_{L, \mathbf{N}}(z) := \sum_{m=1}^{\infty} H_{Lm}
\bigg(\prod_{j=1}^k \mathbf{B}_{N_j}(m)\bigg) z^m.  
$$
\begin{theo}\label{thm:4} For any integers $N_1,N_2, \ldots, N_k\ge 1$ and 
$L\in \{1, 2, \ldots, \max(N_1, \ldots, N_k)\}$, we have 
$
{\bf q}_{L,\mathbf{N}} (z) :=
\exp(\mathbf{G}_{L,\mathbf{N}}(z)/\mathbf{F}_{\mathbf{N}}(z) 
)
\in\mathbb{Z}[[z]].
$
\end{theo}
An outline of the proof of this theorem is given in
Section~\ref{sec:7}, with details being filled in in subsequent sections.

To see that Theorem~\ref{thm:4} implies Conjecture~\ref{conj:zudilin2}, 
we prove in Lemma~\ref{lem:rajout2} in Section~\ref{sec:8} that,
for a given $N$, we have 
\begin{equation}\label{eq:rajout2}
\mathbf{H}_N(m)
= \sum_{j=1}^{\mu} \alpha_j H_{\alpha_jm} -  \sum_{j=1}^{\eta} \beta_j
H_{\beta_jm}. 
\end{equation}
Clearly, the $\alpha_j$'s and $\beta_j$'s (defined in~\eqref{eq:aj}
and~\eqref{eq:bj}) are less than $N$. 
Therefore
$\sum_{j=1}^k \mathbf{H}_{N_j}(m)$ is a finite sum of terms of the form 
$\lambda H_{Lm}$, where $\la$ and $L$ are integers with 
$L\in \{1, 2,\dots, \max(N_1, \ldots, N_k)\}$,
and the claimed implication follows.

\medskip

It is proved in~\cite{villegas} 
that there are exactly fourteen mirror maps ${\bf q_{{\bf N}}}(z)$ associated to complete intersections 
Calabi--Yau threefolds in weighted projective spaces: 
${\bf q}_{(12)}(z)$, ${\bf q}_{(5)}(z)$,  ${\bf q}_{(8)}(z)$, 
${\bf q}_{(10)}(z)$, ${\bf q}_{(3,3)}(z)$, ${\bf q}_{(4,2,2)}(z)$, ${\bf q}_{(2,2,2,2)}(z)$, ${\bf q}_{(4,4)}(z)$, 
${\bf q}_{(6,6)}(z)$, 
${\bf q}_{(4,3)}(z)$, 
${\bf q}_{(6,2,2)}(z)$, 
${\bf q}_{(3,2,2)}(z)$, 
${\bf q}_{(3,6)}(z)$, 
${\bf q}_{(6,4)}(z)$.
All but the first one were already described in~\cite[Sec.~6.4]{batstrat}.  
Theorem~\ref{thm:4} implies that they all
have integral Taylor coefficients, 
which is a new result for the last five of them.

More generally, for any $\mathbf{N}$, the  
functions $\mathbf{F}_{\mathbf{N}}(z)$ and 
$\mathbf{G}_{\mathbf{N}}(z)+\log(z)\mathbf{F}_{\mathbf{N}}(z)$ 
satisfy a hypergeometric differential
equation ${\bf L}y=0$ with maximal unipotent monodromy at the origin, 
where the differential operator ${\bf L}$ is defined by 
\begin{equation}\label{eq:equadiff}
{\bf L}:=\left(z\frac{\dd}{\dd z}\right)^{\varphi(N_1)+\cdots +\varphi(N_k)} 
- z C_{\bf N} \prod_{j=1}^k\prod_{i=1}^{\varphi(N_j)} 
\left(z\frac{\dd}{\dd z} + \frac{r_{i,j}}{N_j}\right).
\end{equation}
Here, $C_{\bf N}=C_{N_1}C_{N_2}\cdots C_{N_k}$ 
and the $r_{i,j}\in \{1, 2, \ldots, N_j\}$ form the residue 
classes mod $N_j$ which are coprime to $N_j$.
This equation is the Picard--Fuchs equation 
of the mirror Calabi--Yau family of 
a one parameter family of Calabi--Yau complete intersections in a 
weighted projective space (see~\cite{corti} and~\cite[Sec.~3]{hosono}). 

\medskip
It is natural to expect 
refinements of Theorem~\ref{thm:4} in the spirit of Theorem~\ref{thm:2} 
and Theorem~\ref{thm:3a}. We did not make a systematic research in
this direction, but it could be  
interesting to do so. For example, in the case $k=1$ and $\mathbf{N}=(6)$, 
it seems that the following relations are best possible:
${\bf q}_{1,(6)}(z)^{1/60}$,  
${\bf q}_{2,(6)}(z)^{1/6}$,  
${\bf q}_{3,(6)}(z)^{1/2}$, 
${\bf q}_{4,(6)}(z)$, 
${\bf q}_{5,(6)}(z)$ and 
${\bf q}_{6,(6)}(z)$ are in $\mathbb{Z}[[z]]$. 

As a first step towards such refinements, 
we prove in Lemma~\ref{lem:diviBB} in Section~\ref{sec:8} that
$\mathbf{B}_{\mathbf{N}}(1)$ always divides  
$\mathbf{B}_{\mathbf{N}}(m)$ for any $m\ge 1$ and any $\mathbf{N}$.
Our techniques enable us to deduce that  
${\bf q}_{1,\mathbf{N}}(z)^{1/\mathbf{B}_{\mathbf{N}}(1)} \in
\mathbb{Z}[[z]]$, which proves the above assertion that 
${\bf q}_{1,(6)}(z)^{1/60}\in \mathbb{Z}[[z]]$.
In fact, this is optimal because ${\bf q}_{1,\mathbf{N}}(z)=
1+\mathbf{B}_{\mathbf{N}}(1)z+\mathcal{O}(z^2)$. 
It turns out that 
$\mathbf{B}_{\mathbf{N}}(1)$ is a natural generalisation of
the quantity $M_{\mathbf{N}}$ which
appears in Theorem~\ref{thm:2}. However, for larger values of the
parameter $L$, we do not know what 
the analogue of the quantity $M_{\mathbf{N}}/\Th_L$
appearing in Theorem~\ref{thm:2} would be.

\medskip

Let us consider the hypergeometric functions $F$ with rational parameters 
for which there exists a constant $C=C(F)>0$ such
that the Taylor coefficients of $F(Cz)$ are of the form 
\begin{equation}\label{eq:quotientfact}
\frac{\prod_{j=1}^r(a_j m)!}{\prod_{j=1}^s (b_j m)!},
\end{equation}
where $m$, the $a_j$'s and $b_j$'s are non-negative
integers, and $\sum_{j=1}^r a_j=\sum_{j=1}^s b_j$. 
A well-known theorem of Landau 
(cf.\ \cite[Proposition 1.1]{villegas}) determines the 
quotients of the form~\eqref{eq:quotientfact} which are integer-valued
for all integers $m\ge 0$.
The quantities $\mathbf{B}_{\mathbf{N}}(m)$ are such examples but do not 
describe all these quotients, 
as demonstrated by the example $(3m)!/(m!\,(2m)!)$.  

In our context, it is therefore natural to ask the following question: 
{\it what are the hypergeometric functions $F$ {\em(}with $z$ 
changed to $Cz$ for a suitable $C${\em)} such that 
\begin{itemize}
\item[$(i)$] their Taylor coefficients are integers and of the 
form~\eqref{eq:quotientfact},
\item[$(ii)$] the associated ``mirror map'' type functions
{\em(}that is, functions defined in a manner analogous to 
${\bf q}_{\bf N}(z)${\em)} have integral Taylor coefficients? 
\end{itemize}}
Of course, 
Condition~$(ii)$ has a sense only if the underlying differential 
equation has two solutions 
$F(z)$ and $G(z)+\log(z) F(z)$ with $F(z)$ and $G(z)$ 
holomorphic at $0$, 
which in turn implies that at least one of the lower parameters of $F$ 
is equal to $1.$ 
Rodriguez--Villegas proved in~\cite{villegas} a result which, once translated 
in our setting, says that the functions $F$ 
with maximal unipotent monodromy (i.e., with only $1$'s as lower parameters) 
and satisfying $(i)$ are 
exactly the
functions $\mathbf{F}_{\mathbf{N}}$ defined by~\eqref{eq:definitionFgrasNgras}. 
Theorem~\ref{thm:4} proves that such functions also satisfy $(ii)$. 
Numerical experiments seem to indicate that there are many examples of such
functions 
$F$ which are not of the form~\eqref{eq:definitionFgrasNgras}.

\subsection{Structure of the paper}
We now briefly review the organisation of 
the rest of the paper. 
Following the steps of previous authors, our approach
for proving Theorems~\ref{thm:1}--\ref{thm:4}
uses $p$-adic analysis. In particular, we make essential
use of Dwork's theory of formal congruences, 
which we recall in Section~\ref{sec:1}. 
Since the details of our proofs are involved, we provide brief
outlines of the proofs of Theorems~\ref{thm:2}--\ref{thm:4}
(with Theorem~\ref{thm:2} implying Theorem~\ref{thm:1}) in separate sections.
Namely, Section~\ref{sec:2} provides an outline of the proof of
Theorem~\ref{thm:2}, 
while Section~\ref{sec:7} provides an outline of the proof of
Theorem~\ref{thm:4}. In both cases, the proof is reduced to a certain
number of lemmas. The lemmas which are necessary for the proof of
Theorem~\ref{thm:2} are established in
Sections~\ref{sec:3}--\ref{sec:6}, while those necessary for the proof of
Theorem~\ref{thm:4} are established in
Sections~\ref{sec:8}--\ref{sec:10}.
Finally, Section~\ref{sec:Xi} contains an outline of the
proof of Theorem~\ref{thm:3}, and Section~\ref{sec:Om}
contains a sketch of the proof of Theorem~\ref{thm:3a}.
Both are largely based on arguments
already used in the proof of Theorem~\ref{thm:2}.
The outline of the proof of Theorem~\ref{thm:3} requires again several auxiliary lemmas, 
whose proofs are postponed to a separate section, Section~\ref{sec:aux}. 
The subsequent section, Section~\ref{sec:DworkKont}, reports on the
evidence to believe (or not to believe) that the value $t_N$
(defined in the next-to-last paragraph before Theorem~\ref{thm:3a})
is given by $t_N=\Xi_N N!$, $N=1,2,\dots$,
while Section~\ref{sec:DK} addresses 
the question of whether the Dwork--Kontsevich sequence $(u_N)_{N\ge1}$
(defined in the last paragraph in Section~\ref{subsec:1}) is (or is
not) given by $u_N=\Om_N N!$, $N=1,2,\dots$.

We draw the reader's attention to the fact that, while the general
line of our approach follows that of previous authors
(particularly \cite{lianyau2}), there
does arise a crucial difference (other than just technical
complications): the reduction and rearrangement of the sums
$C(a+Kp)$ in Sections~\ref{sec:2} and \ref{sec:Xi},
respectively $\mathbf C(a+Kp)$ in
Section~\ref{sec:7}, via the congruence \eqref{eq:J} 
requires a new reduction step, namely
Lemma~\ref{lem:12}, congruence~\eqref{eq:congrconj1U}, and
Lemma~\ref{lem:12a}, respectively. In fact, the
proofs of these two lemmas and the congruence 
form the most difficult parts of our
paper. (In previous work, the
use of \eqref{eq:J} sufficed because the corresponding authors
restricted themselves to $N$ being prime or a prime power.)

\medskip

We expect that (variations of) 
the detailed techniques that we present here can also
serve to prove integrality results for other types of mirror maps,
and in particular some of Observations~1--7 in
\cite[Sec.~4]{almkvist}. We also 
hope that our methods will
turn out to be useful in the context of the very general 
multivariable mirror maps coming from 
the Gelfand--Kapranov--Zelevinsky hypergeometric series: 
see~\cite[Sec.~7.1]{batstrat},~\cite{hosono} and~\cite[Sec.~8]{stienstra} 
for numerous examples related to Calabi--Yau manifolds which are 
complete intersections in products of weighted projective spaces.

\section{Dwork's theory of formal congruences} \label{sec:1}

Since the Taylor coefficients of $F_{\mathbf{N}}(z)$ and
$G_{L,\mathbf{N}}(z)$ (respectively $\mathbf F_{\mathbf{N}}(z)$ and
$\mathbf G_{L,\mathbf{N}}(z)$) are  
explicit, we could also obtain an explicit formula for the
coefficients of $q_{L,\mathbf{N}}(z)$ (respectively 
$\mathbf q_{L,\mathbf{N}}(z)$).  
Unfortunately, such a formula does not seem to be very useful to
prove the desired integrality properties.  
Instead, we will follow the authors of the 
articles cited in the Introduction, who all relied on Dwork's theory. 

First, consider a formal power series $S(z)\in \mathbb{Q}[[z]]$ and suppose 
that we want to prove that $S(z)\in\mathbb{Z}[[z]]$. 

\begin{lem} \label{lem:1} 
Let $S(z)$ be a power series in $\mathbb{Q}[[z]]$. If 
$S(z)\in\mathbb{Z}_p[[z]]$ for any prime number $p$, then
$S(z)\in\mathbb{Z}[[z]]$. 
\end{lem}
This is a consequence of the fact that, given $x\in \mathbb{Q}$, we have  
$x\in \mathbb{Z}$ if and only if 
$x\in \mathbb{Z}_p$ for all prime numbers $p$.
Hence we can work in $\mathbb{Q}_p$ for any fixed prime $p$.

\begin{lem}[\sc ``Dwork's Lemma''] \label{lem:2} 
Let $S(z)\in 1+z\mathbb{Q}_p[[z]]$. Then, 
we have $S(z)\in 1+z\mathbb{Z}_p[[z]]$ if and only if 
$$
\frac{S(z^p)}{S(z)^p }
\in 1+p z\mathbb{Z}_p[[z]].
$$
\end{lem}
\begin{proof} The proof is neither difficult nor long and can 
for example be found in the book of
Lang~\cite[Ch.~14, p.~76]{lang}. Lang attributes this lemma to 
Dieudonn\'e and Dwork.
\end{proof}

We now suppose that $S(z)=\exp(T(z)/\tau)$ for some $T(z)\in
z\mathbb{Q}[[z]]$ and some integer $\tau\ge 1$.
Dwork's Lemma implies the following result: 
$\tau$ being any fixed positive integer, 
we have $\exp\big(T(z)/ \tau \big) \in 1+z\mathbb{Z}_p[[z]]$ if 
and only if $T(z^p)-pT(z)\in p \tau z\mathbb{Z}_p[[z]]$.
(See~\cite[Corollary~6.7]{lianyau} for a proof.)
Since we will be interested in the case when $T(z)=g(z)/f(z)$ with
$f(z)\in 1+z\mathbb{Z}[[z]]$ and $g(z)\in z\mathbb{Q}[[z]]$, 
we state this result as follows.
\begin{lem}\label{lem:4} Given two formal power series 
$f(z)\in 1+z\mathbb{Z}[[z]]$ and $g(z)\in z\mathbb{Q}[[z]]$ and an
integer $\tau\ge 1$, we have 
$\exp\big(g(z)/(\tau f(z))\big) \in 1+z\mathbb{Z}_p[[z]]$ if and only if 
\begin{equation}\label{eq:uv=pvu}
f(z)g(z^p)-p\,f(z^p)g(z)\in p\tau z\mathbb{Z}_p[[z]].
\end{equation}
\end{lem}

Because of the special form of the functions which will play the role of 
$f(z)$ and $g(z)$, we 
will be able to deduce~\eqref{eq:uv=pvu} from the following crucial
result, also  
due to Dwork (see~\cite[Theorem~1.1]{dwork}). We state it in slightly 
bigger generality than 
necessary (i.e., our functions will be independent of $r$).

\begin{prop}[\sc ``Dwork's Formal Congruences Theorem''] 
\label{prop:dworkcongruence}
For all integers\break 
$r\ge 0$, let $A_r : \mathbb{Z}_{\ge 0}\to
\mathbb{Q}_p^{ \times}$, $g_r : \mathbb{Z}_{\ge 0}  
\to \mathbb{Z}_p\setminus\{0\}$ be mappings such that 

$(i)$ $\vert A_r(0)\vert_p =1$; 

$(ii)$ $A_r(m) \in g_r(m)\mathbb{Z}_p$;

$(iii)$ for all integers 
$u,v,n,r,s\ge 0$ such that $0\le u<p^s$ and 
$0 \le v<p$, 
we have 
\begin{equation} \label{eq:DworkA}
\frac{A_r(v+up +np^{s+1})}{A_r(v+up)}-\frac{A_{r+1}(u
+np^{s})}{A_{r+1}(u)} \in p^{s+1} \frac{g_{r+s+1}(n)}{g_r(v+up)}\,
\mathbb{Z}_p. 
\end{equation}
Furthermore, let $F(z) = \sum_{m=0}^{\infty} A_0(m) z^m$, 
$G(z) = \sum_{m=0}^{\infty} A_1(m) z^m$, and 
$$
F_{m,s}(z)=\sum_{j=mp^s}^{(m+1)p^s-1} A_0(j)z^j, \quad 
G_{m,s}(z)=\sum_{j=mp^s}^{(m+1)p^s-1} A_1(j)z^j.
$$
Then, for any integers  $m, s \ge 0$, we have 
\begin{equation}\label{eq:newcongruence}
G(z^p)F_{m,s+1}(z) - F(z)G_{m,s}(z^p)\in p^{s+1} g_s(m) \mathbb{Z}_p[[z]],
\end{equation}
or, equivalently, 
\begin{equation}\label{eq:newcongruence2}
\sum _{j=mp^s} ^{(m+1)p^s-1}\big(A_0(a+jp)A_1(K-j)-A_1(j)A_0(a+(K-j)p)
\big)\in p^{s+1} g_s(m) \mathbb{Z}_p
\end{equation}
for all $a$ and $K$ with $0\le a<p$ and $K\ge0$.
\end{prop}

\begin{Remarks} 
(a) Dwork proved his theorem with 
$A_r:\mathbb{Z}_{\ge 0}\to\mathbb{C}_p^{\times}$, $g_r:\mathbb{Z}_{\ge 0}\to\mathcal{O}_p\setminus\{0\}$ and 
 $A_r(m) \in g_r(m)\mathcal{O}_p$ 
(where $\mathcal{O}_p$ is the ring of integers in $\mathbb{C}_p$).
He obtained
a result similar to~\eqref{eq:newcongruence} and \eqref{eq:newcongruence2}, 
with $\mathcal{O}_p$ instead of
$\mathbb{Z}_p$. In our more restrictive setting,~\eqref{eq:newcongruence} 
and \eqref{eq:newcongruence2} 
hold because $\big(p^{s+1}g_s(m)\mathcal{O}_p\big) \cap \mathbb{Q}_p = 
p^{s+1}g_s(m)\mathbb{Z}_p.$

(b) For any integers  $a$ and $K$ with $0\le a<p$ and
$K\ge0$, the sum  
\begin{equation}\label{eq:rajout1}
\sum _{j=mp^s} ^{(m+1)p^s-1}\big(A_0(a+jp)A_1(K-j)-A_1(j)A_0(a+(K-j)p)
\big)
\end{equation}
is exactly the $(a+pK)$-th Taylor coefficient of 
$G(z^p)F_{m,s+1}(z) - F(z)G_{m,s}(z^p)$, 
which explains the equivalence between the formal
congruence~\eqref{eq:newcongruence} and the  
congruence~\eqref{eq:newcongruence2}. Note that in~\eqref{eq:rajout1} the 
value of $A_0$ and $A_1$ at negative integers 
must be taken as $0$.

(c) Most authors chose $g_s(m)=1$ or a constant in $m$ and $s$. We will
use instead  
$g_s(m)=A_s(m)$: this choice has already been made by Dwork
in~\cite[Sec.~2, p.~37]{dworkihes}.

(d) Dwork also applied his methods to the problems 
considered in the present paper. Indeed, he proved a 
result,
namely \cite[p.~311, Theorem~4.1]{dwork}, which implies that for any
prime $p$ that does not divide $N_1N_2\cdots N_k$, 
the mirror maps ${\mathbf q}_{{\mathbf N}}(z)$ have Taylor 
coefficients in $\mathbb{Z}_p$ (see~\cite[Proposition 2]{zud} for
details). This fact was used by the  
authors of~\cite{lianyau1, lianyau, zud} who 
focussed essentially 
on the remaining case when $p$ divides $N_1N_2\cdots N_k$. Our
approach is different, for we make no distinction of  
this kind between prime numbers.
\end{Remarks}

\section{Outline of the proof of Theorem~\ref{thm:2}}\label{sec:2}

In this section, we provide a brief outline of the proof of
Theorem~\ref{thm:2}, reducing it to
Lemmas~\ref{lem:12}--\ref{lem:11}. These lemmas are subsequently
proved in Sections~\ref{sec:4}--\ref{sec:6}, with two auxiliary lemmas
being the subject of the subsequent section.

By Dwork's Lemma (or rather its consequence given in
Lemma~\ref{lem:4}), we want to prove that  
$$F_{\mathbf N}(z)G_{L,\mathbf N}(z^p)-pF_{\mathbf N}(z^p)G_{L,\mathbf
N}(z) \in p
\frac{M_{\mathbf{N}}}{\Th_L} z \mathbb{Z}_p[[z]].$$ 
We follow the  
presentation in~\cite{lianyau2} and we let $0\le a<p$ and $K\ge 0$. 

The $(a+Kp)$-th Taylor coefficient of
$F_{\mathbf N}(z)G_{L,\mathbf N}(z^p)-pF_{\mathbf N}(z^p)G_{L,\mathbf
N}(z)$ is  
\begin{equation} \label{eq:C} 
C(a+Kp):=\sum_{j=0}^K B_{\mathbf N}(a+jp)B_{\mathbf N}(K-j) (H_{L(K-j)}-pH_{La+Ljp}),
\end{equation}
where $B_{\mathbf N}(m)=\prod_{j=1}^k B_{N_j}(m)$ with
$B_{N}(m):=\frac{(Nm)!}{m!^{N}}$  
(not to be confused with $\mathbf{B}_{\mathbf N}(m)$ and
$\mathbf{B}_{N}(m)$).   
In view of Lemma~\ref{lem:4}, proving Theorem~\ref{thm:2} is equivalent to
proving that 
\begin{equation} \label{eq:Ccong} 
C(a+Kp) \in p \frac {M_{\mathbf{N}}} {\Th_L}\mathbb{Z}_p
\end{equation}
for all primes $p$ and non-negative integers $a$ and $K$ with $0\le a<p$.

The following simple lemma will be frequently used in the sequel.

\begin{lem} \label{lem:multinomial/N!}
For all integers $m\ge 1$ and $N\ge 1$, we have 
$$
B_N(m) \in N! \,\mathbb{Z}.
$$
\end{lem}
\begin{proof} Set $U_m(N)=\frac{(Nm)!}{m!^N N!}$. For any $m, N\ge 1$,
we have the trivial relation  
$$
U_m(N+1) = \binom{Nm+m-1}{m-1} U_m(N).
$$
Therefore, since $U_m(1)=1$, the result follows by induction on $N$.
\end{proof}
We deduce in particular that $B_{\mathbf N}(m)\in M_{\mathbf{N}}\mathbb{Z}$ for
any $m\ge 1$. 

\medskip

Since 
$$
H_J = \sum_{j=1} ^{\fl{J/p}} \frac{1}{pj} + 
\underset{p\nmid j}{\sum_{j=1} ^{J}}\;
\frac{1}{j},
$$
we have
\begin{equation} \label{eq:J} 
pH_{J} \equiv H_{\fl{J/p}} \mod p\mathbb{Z}_p.
\end{equation}
Applying this to $J=La+Ljp$, we get
$$
pH_{La+Ljp}\equiv H_{\lfloor La/p\rfloor+Lj} \mod p\mathbb{Z}_p.
$$
By Lemma~\ref{lem:multinomial/N!}, we infer
\begin{equation}\label{eq:firstreduction}
C(a+Kp) \equiv \sum_{j=0}^K B_{\mathbf N}(a+jp)B_{\mathbf N}(K-j) (H_{L(K-j)}-H_{\lfloor
La/p\rfloor+Lj}) \mod pM_{\mathbf{N}}\mathbb{Z}_p.
\end{equation}
Indeed, if $K\ge 1$ or $a\ge1$, this is because 
$a+jp$ and $K-j$ cannot be simultaneously zero and therefore at least one of 
$B_{\mathbf N}(a+jp)$ or $B_{\mathbf N}(K-j)$ is divisible by $M_{\mathbf{N}}$ by 
Lemma~\ref{lem:multinomial/N!}. In the 
remaining case $K=a=j=0$, we note that the difference of harmonic numbers
in~\eqref{eq:firstreduction} 
is equal to $0$, and therefore the congruence~\eqref{eq:firstreduction} 
holds trivially because $C(0)=0$.
\medskip  

We now want to transform the sum on the right-hand side 
of~\eqref{eq:firstreduction} to a more manageable expression. 
In particular, we want to get rid of the floor function $\fl{La/p}$.
In order to achieve this, we will prove the following lemma
in Section~\ref{sec:4}.

\begin{lem} \label{lem:12}
For any prime $p$, non-negative integers $a$ and $j$ with $0\le a<p$, 
positive integers $N_1,N_2,\dots,N_k$, and $L\in\{1,2, \ldots, 
\max(N_1, \ldots, N_k)\}$, we have 
\begin{equation} \label{eq:congrconj1}
B_{\mathbf N}(a+pj)\left(H_{Lj+\lfloor La/p\rfloor} - H_{Lj}\right) \in p\frac
{M_{\mathbf{N}}} {\Th_L}\mathbb{Z}_p.
\end{equation}
\end{lem}
 
It follows from Eq.~\eqref{eq:firstreduction} and Lemma~\ref{lem:12} that 
 $$
C(a+Kp) \equiv \sum_{j=0}^K B_{\mathbf N}(a+jp)B_{\mathbf N}(K-j) \Big(H_{L(K-j)}-H_{Lj}\Big)
\mod p\frac {M_{\mathbf{N}}} {\Th_L}\mathbb{Z}_p\ ,
$$
which can be rewritten as 
\begin{equation}\label{eq:beforelemma4.2}
C(a+Kp) \equiv -\sum_{j=0}^K
H_{Lj}\Big(B_{\mathbf N}(a+jp)B_{\mathbf N}(K-j)-B_{\mathbf N}(j)B_{\mathbf N}(a+(K-j)p)\Big) \mod p\frac 
{M_{\mathbf{N}}} {\Th_L}\mathbb{Z}_p\ .
\end{equation}

We now use a combinatorial lemma due to Dwork (see~\cite[Lemma
4.2]{dwork}) which provides   
an alternative way to write the sum on the
right-hand side of~\eqref{eq:beforelemma4.2}: namely, we have 
\begin{equation} \label{eq:107}
\sum_{j=0}^K H_{Lj}\Big(B_{\mathbf N}(a+jp)B_{\mathbf N}(K-j)-B_{\mathbf N}(j)B_{\mathbf N}(a+(K-j)p)\Big)
=
\sum _{s=0} ^{r}
\sum _{m=0} ^{p^{r+1-s}-1}Y_{m,s},
\end{equation}
where $r$ is such that $K<p^r$, and 
\begin{equation*} 
Y_{m,s}:=\big(H_{Lmp^s}-H_{L\fl{m/p}p^{s+1}}\big)S(a,K,s,p,m),
\end{equation*}
the expression $S(a,K,s,p,m)$ being defined by
$$
S(a,K,s,p,m):=\sum _{j=mp^s} ^{(m+1)p^s-1}\big(B_{\mathbf N}(a+jp)B_{\mathbf N}(K-j)-B_{\mathbf N}(j)B_{\mathbf N}(a+(K-j)p)
\big).
$$
In this expression for $S(a,K,s,p,m)$, it is assumed that $B_{\mathbf
N}(n)=0$ for negative integers~$n$. 

\medskip

It would suffice to prove that 
\begin{equation}
\label{eq:yms}
Y_{m,s}\in p\frac {M_{\mathbf{N}}} {\Th_L}\mathbb Z_p 
\end{equation}
because~\eqref{eq:beforelemma4.2} and \eqref{eq:107}
would then imply that 
$
C(a+Kp) \in  p\frac
{M_{\mathbf{N}}} {\Th_L}\mathbb{Z}_p, 
$
as desired. 

We will prove \eqref{eq:yms} in the following manner.
The expression for $S(a,K,s,p,m)$ is of the form considered in 
Proposition~\ref{prop:dworkcongruence}. 
The proposition will enable us to prove the
following fact in Section~\ref{sec:5}.
\begin{lem} \label{lem:10}
For all primes $p$ and
non-negative integers $a,m,s,K$ with $0\le a<p$, we have 
\begin{equation}
\label{eq:congruenceS}
S(a,K,s,p,m)\in p^{s+1} B_{\mathbf N}(m) \mathbb{Z}_p.
\end{equation}
\end{lem}
Furthermore, in Section~\ref{sec:6} we shall prove the following lemma.
\begin{lem} \label{lem:11}
For all primes $p$, non-negative integers $m$, 
positive integers $N_1,N_2,\dots,N_k$, and $L\in\{1,2, \ldots, 
\max(N_1, \ldots, N_k)\}$, we have
\begin{equation} 
\label{eq:110}
B_{\mathbf N}(m)\big(H_{Lmp^s}-H_{L\fl{m/p}p^{s+1}}\big)\in \frac {M_{\mathbf{N}}}
{p^s\Th_L}\mathbb Z_p\ .
\end{equation}
\end{lem}

It is clear that~\eqref{eq:congruenceS} and~\eqref{eq:110}
imply~\eqref{eq:yms}. 
This completes the outline of the proof of Theorem~\ref{thm:2}.

\section{Two auxiliary lemmas} \label{sec:3}

The proof of Lemma~\ref{lem:10} in Section~\ref{sec:5} requires two further
auxiliary results, given in Lemmas~\ref{lem:6} and \ref{lem:gammap}
below, the proofs of which form the contents of this intermediary section.
In the proof of the first lemma, and also in later proofs, we shall
frequently make use of Legendre's formula for the $p$-adic valuation of $n!$,
where $n$ is a non-negative integer, which we recall here for
convenience:
\begin{equation} \label{eq:Leg}
v_p(n!) = \sum_{k=1}^{\infty} \left\lfloor \frac{n}{p^k} \right\rfloor.
\end{equation}

\begin{lem} \label{lem:6}
For all integers $N\ge 1$, $n\ge 0$, $s\ge 0$, $0\le u< p^s$, 
$p$ prime, 
we have
$$
\frac{B_N(u+np^s)}{B_N(u)}
\in B_N(n) \,\mathbb{Z}_p,
$$
where $B_N(m)$ is defined after \eqref{eq:C}.
\end{lem}
\begin{proof} 
We may rewrite the quotient in the assertion of the lemma as
\begin{align}
\frac{(Nu+Nnp^s)!}{(u+np^s)!^N}\cdot \frac{u!^N}{(Nu)!} 
&= \frac{(Nu+Nnp^s)!}{(Nu)!\,(Nnp^s)!}\cdot\frac{(Nnp^s)!}{(np^s)!^N} 
\cdot
\frac{(np^s)!^Nu!^N}{(u+np^s)!^N}
\notag
\\
&= \binom{Nu+Nnp^s}{Nu} \cdot \frac{(Nnp^s)!}{(np^s)!^N} \cdot
\binom{u+np^s}{u}^{-N}.
\label{eq:lem6}
\end{align}

The binomial coefficient $\binom{Nu+Nnp^s}{Nu}$ is an integer. 
For the multinomial coefficient
we observe that
\begin{align*}
v_p\big((Nnp^s)!/(np^s)!^N\big)
&=\sum _{k=1} ^{\infty}\left(\fl{\frac {Nnp^s} {p^k}}-
N\fl{\frac {np^s} {p^k}}\right)
\\ 
&=
\sum _{k=1} ^{\infty}\left(\fl{\frac {Nn} {p^k}}-
N\fl{\frac {n} {p^k}}\right)
\\
&=v_p\big((Nn)!/n!^N\big),
\end{align*}
which means that 
$$
\frac{(Nnp^s)!}{(np^s)!^N}  \in \frac{(Nn)!}{n!^N} \,\mathbb{Z}_p.
$$

Hence, to conclude the proof of the lemma, it suffices to show that 
\begin{equation} \label{eq:lem6v} 
v_p\bigg(\binom{u+np^s}{u}\bigg)= 0.
\end{equation}
In order to do so, we start with the formula
$$
v_p\bigg(\binom{u+np^s}{u}\bigg) = \sum_{k=1}^{\infty} 
\bigg(\left\lfloor\frac{u+np^s}{p^k}\right\rfloor-
\left\lfloor\frac{u}{p^k}\right\rfloor-\left\lfloor\frac{n
p^s}{p^k}\right\rfloor\bigg).
$$
We distinguish between two cases.

\medskip

(a)  $k\le s$. In this case, $(np^s)/p^k$ is an integer and therefore
$$
\left\lfloor\frac{u+np^s}{p^k}\right\rfloor-
\left\lfloor\frac{u}{p^k}\right\rfloor-\left\lfloor\frac{n
p^s}{p^k}\right\rfloor
 =
\left\lfloor\frac{u}{p^k}\right\rfloor +\frac{n p^s}{p^k}
-\left\lfloor\frac{u}{p^k}\right\rfloor- \frac{n p^s}{p^k} 
=0.
$$

(b)  $k> s$. In this case, we have $\lfloor u/p^k\rfloor=0$ because
$0\le u<p^s< p^k.$  

\medskip

Hence,  
$$
v_p\bigg(\binom{u+np^s}{u}\bigg) = \sum_{k=s+1}^{\infty} 
\bigg(\left\lfloor\frac{u+np^s}{p^k}\right\rfloor
-\left\lfloor\frac{n}{p^{k-s}}\right\rfloor\bigg).
$$
Since $0\le u/p^k<1$, we have 
$$
\left\lfloor\frac{u+np^s}{p^k}\right\rfloor 
-\left\lfloor\frac{n}{p^{k-s}}\right\rfloor \in \{0,1\},
$$
so that it suffices to show that the value $1$ cannot be attained. 
Arguing by contradiction, let us assume that,  
for some $n,u, p, s$, we have 
$$
\left\lfloor\frac{u+np^s}{p^k}\right\rfloor -
\left\lfloor\frac{n}{p^{k-s}}\right\rfloor=1.
$$
Writing, for simplicity,
$X=\Big\lfloor \dfrac{u+np^s}{p^k}\Big\rfloor$, our assumption
implies that  
\begin{equation}
\label{eq:inegalite}
\frac{n}{p^{k-s}}< X \le \frac{u+np^s}{p^k}.
\end{equation}
$X$ being an integer, we thus have 
\begin{equation}\label{eq:inegalite2}
\frac{n}{p^{k-s}} + \frac{1}{p^{k-s}} \le  X.
\end{equation}
Combining the right inequality of~\eqref{eq:inegalite} 
with~\eqref{eq:inegalite2}, we deduce that
$$
\frac{n}{p^{k-s}} + \frac{1}{p^{k-s}} \le \frac{u}{p^k} + \frac{np^s}{p^k},
$$
which means that $u\ge p^s$. This contradicts the hypothesis that
$u<p^s$ and finishes the proof. 
\end{proof}

During the proof of Lemma~\ref{lem:10} in Section~\ref{sec:5}, we will
also use certain properties of the $p$-adic gamma 
function $\Gamma_p$. This function is defined on integers $n\ge 1$ by 
$$
\Gamma_p(n) = (-1)^n \prod_{\stackrel{k=1}{(k,p)=1}}^{n-1} k.
$$
We will not consider its extension to $\mathbb{Q}_p.$
In the following lemma,
we collect the results on $\Gamma_p$ that we shall need later on.
\begin{lem}\label{lem:gammap}
$(i)$ For all integers $n\ge 1$, we have 
\begin{equation*}
\frac{(np)!}{n!} = (-1)^{np+1}p^n \Gamma_p(1+np).
\end{equation*}
$(ii)$ For all integers $k\ge 1,n\ge 1,s\ge 0$, we have 
\begin{equation*}
\Gamma_p(k+np^s) \equiv \Gamma_p(k) \mod p^s.
\end{equation*}
\end{lem}
\begin{proof} 
See~\cite[Lemma~7]{zud} for $(i)$ and~\cite[p.~71, Lemma~1.1]{lang} for $(ii)$.
\end{proof}

\section{Proof of Lemma~\ref{lem:12}} \label{sec:4}

The assertion is trivially true if $\lfloor La/p\rfloor=0$, that is,
if $0\le a<p/L$. We may hence assume that $p/L\le a<p$ from now on.
A further assumption upon which we agree without loss of generality
for the rest of the proof is that $N_k=\max(N_1, \ldots, N_k)$.

\subsection{First part: a weak version of Lemma~\ref{lem:12}}
In a first step, we prove that
\begin{equation} \label{eq:congrconj2}
B_{\mathbf N}(a+pj)\left(H_{Lj+\lfloor La/p\rfloor} - H_{Lj}\right) \in
p\mathbb{Z}_p\ .
\end{equation}
(The reader should note the absence of the term
$M_{\mathbf{N}}/\Th_L$ in comparison with \eqref{eq:congrconj1}.)

\medskip
For the proof of \eqref{eq:congrconj2}, we note that
the $p$-adic valuation of $B_{\mathbf N}(a+pj)$ is equal to
$$v_p\big(B_{\mathbf N}(a+pj)\big)=
\sum _{i=1} ^{k}\sum _{\ell =1} ^{\infty}\left(\fl{\frac {N_i(a+pj)}
{p^\ell }}- N_i\fl{\frac {a+pj} {p^\ell }}\right).$$
Obviously, all the summands in this sum are non-negative, whence, in
particular,
\begin{equation} \label{eq:111}
v_p\big(B_{\mathbf N}(a+pj)\big)\ge
\sum _{\ell =1} ^{\infty}\left(\fl{\frac {N_k(a+pj)} {p^\ell }}-
N_k\fl{\frac {a+pj} {p^\ell }}\right).
\end{equation}
On the other hand, by definition of the harmonic numbers, we have
$$
H_{Lj+\lfloor La/p\rfloor} - H_{Lj}=\frac {1} {Lj+1}+\frac {1} {Lj+2}+\dots+
\frac {1} {Lj+\lfloor La/p\rfloor}.
$$
It therefore suffices to show that
\begin{equation} \label{eq:100} 
v_p\big(B_{\mathbf N}(a+pj)\big)\ge 1+ 
\max_{1\le\ep\le \lfloor La/p\rfloor}v_p(Lj+\ep).
\end{equation}
The lower bound on the right-hand side of \eqref{eq:111}
can, in fact, be simplified since $0\le a<p$; namely, we have
\begin{equation} \label{eq:104} 
\fl{\frac {a+pj} {p^\ell }}=\fl{\frac {j} {p^{\ell -1}}}.
\end{equation}

For a given integer 
$\ep$ with $1\le\ep\le \lfloor La/p\rfloor$, let
$Lj+\ep=p^d\be$, where $d=v_p(Lj+\ep)$. If we use this notation in 
\eqref{eq:111}, together with \eqref{eq:104}, we obtain 
\begin{equation} \label{eq:103} 
v_p\big(B_{\mathbf N}(a+pj)\big)\ge
\sum _{\ell =1} ^{\infty}\left(\fl{\frac {N_ka} {p^\ell }-\frac
{N_k\ep}
{Lp^{\ell -1}}+ \frac {N_k\be} {L} p^{d+1-\ell }}-
N_k\fl{-\frac {\ep}
{Lp^{\ell -1}}+ \frac {\be} {L} p^{d+1-\ell }}\right).
\end{equation}
Since $\ep\le \lfloor La/p\rfloor$, 
we have $\frac {N_ka} {p^\ell }-\frac {N_k\ep}
{Lp^{\ell -1}}\ge0$, whence
\begin{equation} \label{eq:101} 
\fl{\frac {N_ka} {p^\ell }-\frac {N_k\ep}
{Lp^{\ell -1}}+ \frac {N_k\be} {L} p^{d+1-\ell }}\ge
\fl{ \frac {N_k\be} {L} p^{d+1-\ell }}.
\end{equation}
Clearly, we also have
\begin{equation} \label{eq:102}
 \fl{-\frac {\ep}
{Lp^{\ell -1}}+ \frac {\be} {L} p^{d+1-\ell }}\le
\fl{ \frac {\be} {L} p^{d+1-\ell }}.
\end{equation}
If we use \eqref{eq:101} and \eqref{eq:102} 
in \eqref{eq:103}, then we obtain
\begin{equation} \label{eq:105} 
v_p\big(B_{\mathbf N}(a+pj)\big)\ge
\sum _{\ell =1} ^{\infty}\left(\fl{
 \frac {N_k\be} {L} p^{d+1-\ell }}-
N_k\fl{ \frac {\be} {L} p^{d+1-\ell }}\right).
\end{equation}
Now we claim that $\be p^{d+1-\ell }/L$ cannot be an integer.
Indeed, if it were, then $L\ga p^{\ell -1}=\be p^d=Lj+\ep$ for a suitable
integer $\ga$. It would follow that $L$ divides $\ep$, contradicting
$1\le \ep\le La/p<L$.
However, the fact that $\be p^{d+1-\ell }/L$ is not an integer entails that
$$
 \frac {\be} {L} p^{d+1-\ell }-
\fl{ \frac {\be} {L} p^{d+1-\ell }}\ge\frac {1} {L},$$
as long as $\ell \le d+1$. Multiplication of both sides of this inequality by
$N_k$ leads to the chain of inequalities
$$
 \frac {N_k\be} {L} p^{d+1-\ell }-
N_k\fl{ \frac {\be} {L} p^{d+1-\ell }}\ge\frac {N_k} L\ge
1
$$
(it is here where we use the assumption $L\le
N_k=\max(N_1,\dots,N_k)$),
whence
$$\fl{\frac {N_k\be} {L} p^{d+1-\ell }}-N_k\fl{ \frac {\be} {L}
p^{d+1-\ell }}\ge 1, 
$$
provided $\ell \le d+1$.
Use of this estimation in \eqref{eq:105} gives
$$v_p\big(B_{\mathbf N}(a+pj)\big)\ge d+1=1+v_p(Lj+\ep).$$
This completes the proof of \eqref{eq:100}, and, hence, of
\eqref{eq:congrconj2}.

\medskip
For later use, we record that we have in particular shown that
for any 
$$D\le 1+\max_{1\le\ep\le \lfloor La/p\rfloor}v_p(Lj+\ep)$$ 
we have
\begin{equation} \label{eq:sumest} 
\sum _{\ell =2} ^{D}\left(\fl{\frac {N_k(a+pj)} {p^\ell }}-
N_k\fl{\frac {a+pj} {p^\ell }}\right)\ge D-1.
\end{equation}

\medskip
We now embark on the proof of \eqref{eq:congrconj1} itself.

\subsection{Second part: the case $j=0$}
In this case, we want to prove that
\begin{equation} \label{eq:congrconj4} 
B_{\mathbf N}(a)H_{\lfloor La/p\rfloor} \in p\frac {M_{\mathbf{N}}}
{\Th_L}\mathbb{Z}_p\ ,
\end{equation}
or, using \eqref{eq:J} (in the other direction), equivalently
\begin{equation} \label{eq:congrconj5} 
B_{\mathbf N}(a)H_{La} \in \frac {M_{\mathbf{N}}}
{\Th_L}\mathbb{Z}_p\ ,
\end{equation}
The reader should keep in mind that we still assume
that $p/L\le a<p$, so that, in particular, $a>0$. 

If $p>N_k=\max(N_1,\dots,N_k)$, then our claim, in the form
\eqref{eq:congrconj4}, reduces to\break 
$B_{\mathbf N}(a)H_{\lfloor La/p\rfloor} \in
p\mathbb{Z}_p$, 
which is indeed true because of \eqref{eq:congrconj2} with $j=0$.

Now let $p\le N_k$. Evidently, our claim, this time in the form
\eqref{eq:congrconj5}, holds for $a=1$. So, let $a\ge2$ from now on.

In a similar way as we did for the expression in
\eqref{eq:congrconj2}, we bound the $p$-adic valuation of the
expression in \eqref{eq:congrconj5} from below:
\begin{align} \notag
v_p\big(B_{\mathbf N}(a)H_{La}\big)
&=
\sum _{i=1} ^{k}\sum _{\ell =1} ^{\infty}\left(\fl{\frac {N_ia} {p^\ell }}-
N_i\fl{\frac {a} {p^\ell }}\right)+
v_p(H_{La})\notag\\
&\ge\sum _{i=1} ^{k}\sum _{\ell =1} ^{\infty}\fl{\frac {N_ia} {p^\ell }}
-\fl{\log_p La}\notag\\
&\ge\sum _{i=1} ^{k}\sum _{\ell =1} ^{\infty}\fl{\frac {2N_i} {p^\ell }}
-\fl{\log_p Lp}\notag\\
&\ge\fl{\frac {2N_k} {p}}+
\sum _{\ell =2} ^{\infty}\fl{\frac {2N_k} {p^\ell }}+
\sum _{i=1} ^{k-1}\sum _{\ell =1} ^{\infty}\fl{\frac {2N_i} {p^\ell }}
-\fl{\log_p L}-1
\notag\\
&\ge\fl{\frac {N_k} {p}}+
\sum _{i=1} ^{k}\sum _{\ell =1} ^{\infty}\fl{\frac {N_i} {p^\ell }}
-\fl{\log_p L}-1
\label{eq:unglA0}\\
&\ge\max\left\{1,\fl{L/p}\right\}+
\sum _{i=1} ^{k}v_p(N_i!)
-\fl{\log_p L}-1.
\label{eq:unglA}
\end{align}
If $p=2$, then we can continue the estimation \eqref{eq:unglA} as
\begin{equation} \label{eq:unglB} 
v_2\big(B_{\mathbf N}(a)H_{La}\big)\ge 
\sum _{i=1} ^{k}v_2(N_i!)
-\fl{\log_2 L}=
v_2\big(M_{\mathbf{N}}/\Th_L\big),
\end{equation}
where we used Lemma~\ref{lem:H_L} to obtain the equality.
(In fact, at this point it was not necessary to consider the case
$p=2$ because $a<p$ and because we assumed $a\ge2$. However, we shall
re-use the present estimations later in the third part of the current
proof, in a context where $a=1$ is allowed.)

From now on let $p\ge3$.
We use the fact that
\begin{equation} \label{eq:log}
k\ge\fl{\log_pk}+2 
\end{equation}
for all integers $k\ge2$ and primes $p\ge3$.
Thus, in the case that $L\ge2p$, 
the estimation \eqref{eq:unglA} can be continued as
$$
v_p\big(B_{\mathbf N}(a)H_{La}\big)
\ge 
1+\fl{\log_p\fl{L/p}}+
\sum _{i=1} ^{k}v_p(N_i!)
-\fl{\log_p L}
\ge 
\sum _{i=1} ^{k}v_p(N_i!)=v_p( M_{\mathbf{N}}),
$$
implying \eqref{eq:unglB} in this case.
If $p\le L<2p$, 
then the estimation \eqref{eq:unglA} can be continued as
\begin{align*}
v_p\big(B_{\mathbf N}(a)H_{La}\big)
&\ge 
1+
\sum _{i=1} ^{k}v_p(N_i!)
-2=v_p(M_{\mathbf{N}}/\Th_L),
\end{align*}
implying \eqref{eq:unglB} in this case also.
Finally, if $L<p$, it follows from \eqref{eq:unglA} that 
\begin{align*}
v_p\big(B_{\mathbf N}(a)H_{La}\big)
&\ge 
1+
\sum _{i=1} ^{k}v_p(N_i!)
-1=v_p(M_{\mathbf{N}}),
\end{align*}
implying \eqref{eq:unglB} also in this final case.
Everything combined, \eqref{eq:congrconj5} follows immediately.

\subsection{Third part: the case $j>0$}
Now let $j>0$.
If $p>N_k=\max(N_1,\dots,N_k)$, then \eqref{eq:congrconj1} reduces to 
\begin{equation*}
B_{\mathbf N}(a+pj)\left(H_{Lj+\lfloor La/p\rfloor} - H_{Lj}\right) \in
p\mathbb{Z}_p\ ,
\end{equation*}
which is 
again true because of \eqref{eq:congrconj2}.

Now let $p\le N_k$. The reader should keep in mind that we still assume
that $p/L\le a<p$, so that, in particular, $a>0$. 
In a similar way as we did for the expression in
\eqref{eq:congrconj2}, we bound the $p$-adic valuation of the
expression in \eqref{eq:congrconj1} from below. For the sake of
convenience, we write $T_1$ for $\max_{1\le\ep\le \lfloor
La/p\rfloor}v_p(Lj+\ep)$ and $T_2$ for $\fl{\log_p (a+pj)}$. 
Since it is somewhat hidden where our assumption $j>0$ enters the
subsequent considerations, we point out to the reader that
$j>0$ implies that $T_2\ge1$; without this property the split of
the sum over $\ell$ into subsums in the 
chain of inequalities below would be impossible.
So, using the above notation, we have (the detailed explanations for the
various steps are given immediately after the 
following chain of estimations)
{\allowdisplaybreaks
\begin{align} \notag
v_p\Big(&B_{\mathbf N}(a+pj)\left(H_{Lj+\lfloor La/p\rfloor} - H_{Lj}\right)\Big)\\
\notag
&=
\sum _{i=1} ^{k}\sum _{\ell =1} ^{\infty}\left(\fl{\frac {N_i(a+pj)}
{p^\ell }}- 
N_i\fl{\frac {a+pj} {p^\ell }}\right)+
v_p\big(H_{Lj+\lfloor La/p\rfloor} - H_{Lj}\big)\\
&=
\fl{\frac {N_k(a+pj)} {p }}-
N_k\fl{\frac {a+pj} {p }}
+
\sum _{\ell =2} ^{\min\{1+T_1,T_2\}}
\left(\fl{\frac {N_k(a+pj)} {p^\ell }}-
N_k\fl{\frac {a+pj} {p^\ell }}\right)
\notag\\
\notag
&\kern1cm
+
\sum _{\ell =\min\{1+T_1,T_2\}+1} ^{\infty}
\left(\fl{\frac {N_k(a+pj)} {p^\ell }}-
N_k\fl{\frac {a+pj} {p^\ell }}\right)\\
\notag
&\kern1cm
+
\sum _{i=1} ^{k-1}\sum _{\ell =1} ^{\infty}\left(\fl{\frac {N_i(a+pj)}
{p^\ell }}- 
N_i\fl{\frac {a+pj} {p^\ell }}\right)
+v_p\big(H_{Lj+\lfloor La/p\rfloor} - H_{Lj}\big)
\notag
\\
\notag
&\ge
\fl{\frac {N_ka} {p }}+\min\{1+T_1,T_2\}-1
+
\sum _{i=1} ^{k}\sum _{\ell =T_2+1} ^{\infty}
\left(\fl{\frac {N_i(a+pj)} {p^\ell }}-
N_i\fl{\frac {a+pj} {p^\ell }}\right)\\
&\kern1cm
+v_p\big(H_{Lj+\lfloor La/p\rfloor} - H_{Lj}\big)
\label{eq:ungl1}
\\
&\ge
\fl{\frac {N_ka} {p }}
+T_1+v_p\big(H_{Lj+\lfloor La/p\rfloor} - H_{Lj}\big)
+\min\{0,T_2-T_1-1\}
\notag\\
&\kern1cm
+
\sum _{i=1} ^{k}\sum _{\ell =\fl{\log_p(a+pj)}+1}
^{\infty}\left(\fl{\frac {N_i(a+pj)} {p^\ell }}- 
N_i\fl{\frac {a+pj} {p^\ell }}\right)
\label{eq:ungl2}\\
&\ge 
\max\left\{1,\fl{L/p}\right\}
+\min\{0,T_2-T_1-1\}
+
\sum _{i=1} ^{k}\sum _{\ell =1} ^{\infty}\fl{\frac {N_i} {p^\ell
}\cdot\frac {a+pj} {p^{\fl{\log_p(a+pj)}}}}
\label{eq:ungl3}\\
&\ge 
\max\left\{1,\fl{L/p}\right\}+\fl{\log_p (a+pj)}
-\fl{\log_p\big(Lj+\lfloor La/p\rfloor\big)}-1
\notag\\
&\kern1cm
+
\sum _{i=1} ^{k}\sum _{\ell =1} ^{\infty}\fl{\frac {N_i} {p^\ell
}\cdot\frac {a+pj} {p^{\fl{\log_p(a+pj)}}}}
\label{eq:ungl4}
\\
&\ge 
\max\left\{1,\fl{L/p}\right\}+\fl{\log_p j}
-\fl{\log_p\big(Lj+\lfloor La/p\rfloor\big)}
+
\sum _{i=1} ^{k}\sum _{\ell =1} ^{\infty}\fl{\frac {N_i} {p^\ell }}
\label{eq:ungl5}
\\
&\ge 
\max\left\{1,\fl{L/p}\right\}+\fl{\log_p j}
-\fl{\log_p L}-\fl{\log_p\left(j+\frac {1} {L}\lfloor La/p\rfloor\right)}-1
\notag\\
&\kern1cm
+
\sum _{i=1} ^{k}v_p(N_i!)
\label{eq:ungl6}\\
&\ge 
\max\left\{1,\fl{L/p}\right\}
-\fl{\log_p L}-1+v_p(M_{\mathbf{N}}).
\label{eq:ungl7}
\end{align}
}%
Here, we used \eqref{eq:sumest} in order to get \eqref{eq:ungl1}.
To get \eqref{eq:ungl3}, we used the inequalities
\begin{equation} \label{eq:ungl6a} 
\fl{\frac {N_ka} {p}}\ge\fl{\frac {N_k} {p}}\ge
\max\left\{1,\fl{L/p}\right\}
\end{equation}
and 
\begin{equation} \label{eq:ungl100} 
T_1+v_p\big(H_{Lj+\lfloor La/p\rfloor} - H_{Lj}\big)\ge0.
\end{equation}
To get \eqref{eq:ungl4}, we used that
$$T_2-T_1-1\ge \fl{\log_p (a+pj)}
-\fl{\log_p\big(Lj+\lfloor La/p\rfloor\big)}-1$$
and
$$
\fl{\log_p (a+pj)}
-\fl{\log_p\big(Lj+\lfloor La/p\rfloor\big)}-1=
\fl{\log_p j}-\fl{\log_p\big(Lj+\lfloor La/p\rfloor\big)}
\le 0,
$$
so that
\begin{equation} \label{eq:ungl6b}
\min\{0,T_2-T_1-1\}\ge 
\fl{\log_p (a+pj)}
-\fl{\log_p\big(Lj+\lfloor La/p\rfloor\big)}-1.
\end{equation}
Next, to get \eqref{eq:ungl5}, we used
\begin{equation} \label{eq:ungl6d}
\fl{\frac {N_i} {p^\ell }\cdot\frac {a+pj}
{p^{\fl{\log_p(a+pj)}}}}\ge
\fl{\frac {N_i} {p^\ell }}.
\end{equation}
To get \eqref{eq:ungl6}, we used
\begin{equation} \label{eq:ungl6c}
\fl{\log_p\big(Lj+\lfloor La/p\rfloor\big)}
\le\fl{\log_p L}+\fl{\log_p\left(j+\frac {1} {L}\lfloor La/p\rfloor\right)}+1 .
\end{equation}
Finally, we used $\frac {1} {L}\lfloor
La/p\rfloor<1$ in order to get
\eqref{eq:ungl7}. 

If we now repeat the arguments after \eqref{eq:unglA},
then we see that the above estimation implies
\begin{equation} \label{eq:ungl8} 
v_p\Big(B_{\mathbf N}(a+pj)\left(H_{Lj+\lfloor La/p\rfloor} - H_{Lj}\right)\Big)
\ge
v_p\big(M_{\mathbf{N}}/\Th_L\big).
\end{equation}
This almost proves \eqref{eq:congrconj1}, our lower
bound on the $p$-adic valuation of the number in \eqref{eq:congrconj1}
is just by $1$ too low. 

In order to establish that \eqref{eq:congrconj1} is indeed correct, we
assume by contradiction that all the inequalities in the estimations
leading to \eqref{eq:ungl7} and finally to \eqref{eq:ungl8} 
are in fact equalities. 
In particular, the estimations in \eqref{eq:ungl6a} hold with equality only
if $a=1$ and, if $L$ should be at least $p$, also $\fl{N_k/p}=\fl{L/p}$.
We shall henceforth assume both of these two conditions.

If we examine the arguments after \eqref{eq:unglA} that led us from
\eqref{eq:ungl7} to \eqref{eq:ungl8}, then we see that they prove in
fact 
\begin{equation} \label{eq:ungl11} 
v_p\Big(B_{\mathbf N}(a+pj)\left(H_{Lj+\lfloor La/p\rfloor} - H_{Lj}\right)\Big)
\ge
1+v_p\big(M_{\mathbf{N}}/\Th_L\big)
\end{equation}
except if:

\bigskip
{\sc Case 1:} $p=2$ and $\fl{L/2}=1$;

{\sc Case 2:} $p\ge3$ and $p\le L<2p$;

{\sc Case 3:} $p=3$ and $\fl{L/3}=2$;

{\sc Case 4:} $L<p$.
\bigskip

In all other cases, there holds either strict inequality in
\eqref{eq:log} with $k=\fl{L/p}$, or one has $v_p(\Th_L)\ge1$
and is able to show 
$$
v_p\Big(B_{\mathbf N}(a+pj)\left(H_{Lj+\lfloor La/p\rfloor} - H_{Lj}\right)\Big)
\ge
v_p\big(M_{\mathbf{N}}\big),
$$ 
so that \eqref{eq:ungl11} is satisfied, as desired.
We now show that \eqref{eq:ungl11} holds as well in Cases~1--4, thus
completing the proof of \eqref{eq:congrconj1}.

\medskip
{\sc Case 1}. Let first $p=2$ and $L=2$. We then have
\begin{align*}
\min\{0,T_2-T_1-1\}&=\min\{0,\fl{\log_2(2j+1)}-v_2(2j+1)-1\}\\
&=
\min\{0,\fl{\log_2(2j+1)}-1\}=0>-1,
\end{align*}
in contradiction to having equality in \eqref{eq:ungl6b}.

On the other hand, 
if $p=2$ and $L=3$ then, because of equality in the second estimation
in \eqref{eq:ungl6a}, we must have $N_k=3$. We have
$$H_{Lj+\lfloor La/p\rfloor} - H_{Lj}
=H_{3j+1} - H_{3j}=\frac {1} {3j+1}.$$
If there holds equality in \eqref{eq:ungl6b}, then $Lj+\fl{La/p}=3j+1$ 
must be a power of $2$, say $3j+1=2^e$ or, equivalently, 
$j=(2^e-1)/3$. It follows that
$$
\fl{\frac {N_k} {p }\cdot\frac {a+pj}
{p^{\fl{\log_p(a+pj)}}}}
=\fl{\frac {3} {2 }\cdot\frac {1+2j}
{2^{\fl{\log_2(1+2j)}}}}
=\fl{\frac {3} {2 }\cdot\frac {2^{e+1}+1}
{3\cdot 2^{e-1}}}
=2>1=\fl{\frac {3} {2}}=\fl{\frac {N_k} {p}},
$$
in contradiction to having equality in \eqref{eq:ungl6d} with $\ell=1$.

\medskip
{\sc Case 2}. Our assumptions $p\ge3$ and $p\le L<2p$ imply 
$$
H_{Lj+\lfloor La/p\rfloor} - H_{Lj}
=H_{Lj+1} - H_{Lj}=\frac {1} {Lj+1}.
$$
Arguing as in the previous case, in order to have equality in 
\eqref{eq:ungl6b}, we must have $Lj+1=f\cdot p^e$ for some positive
integers $e$ and $f$ with $0<f<p$. Thus, $j=(f\cdot p^e-1)/L$ and
$p<L$. (If $p=L$ then $j$ would be non-integral.) It follows that
\begin{equation} \label{eq:ungl12}
\fl{\frac {N_k} {p }\cdot\frac {a+pj}
{p^{\fl{\log_p(a+pj)}}}}=
\fl{\frac {N_k} {p }\cdot\frac {f\cdot p^{e+1}+L-p}
{L\cdot p^{\fl{\log_p((f\cdot p^{e+1}+L-p)/L)}}}}. 
\end{equation}
If $f=1$, then we obtain from \eqref{eq:ungl12} that
$$
\fl{\frac {N_k} {p }\cdot\frac {a+pj}
{p^{\fl{\log_p(a+pj)}}}}=
\fl{\frac {N_k} {p }\cdot\frac {p^{e+1}+L-p}
{L\cdot p^{e-1}}}
\ge\fl{\frac {p^{e+1}+L-p} {p^{e}}}
>1=\fl{\frac {L} {p}}=\fl{\frac {N_k} {p}}, 
$$
in contradiction with having equality in \eqref{eq:ungl6d} with $\ell=1$.
 
On the other hand, if $f\ge2$, then we obtain from \eqref{eq:ungl12} that
$$
\fl{\frac {N_k} {p }\cdot\frac {a+pj}
{p^{\fl{\log_p(a+pj)}}}}\ge
\fl{\frac {f\cdot p^{e+1}+L-p}
{p^{e+1}}}\ge f>1=\fl{\frac {L} {p}}=\fl{\frac {N_k} {p}}, 
$$
again in contradiction with having equality in \eqref{eq:ungl6d} with $\ell=1$.

\medskip
{\sc Case 3}. Our assumptions $p=3$ and $\fl{L/3}=2$ imply 
$$
H_{Lj+\lfloor La/p\rfloor} - H_{Lj}
=H_{Lj+2} - H_{Lj}=\frac {1} {Lj+1}+\frac {1} {Lj+2}.
$$
Similar to the previous cases, in order to have equality in 
\eqref{eq:ungl6b}, we must have $Lj+\ep=f\cdot 3^e$ for some positive
integers $\ep,e,f$ with $0<\ep,f<3$. The arguments from Case~2 can now
be repeated almost verbatim. We leave the details to the reader.

\medskip
{\sc Case 4}. If $L<p$, then $p/L>1=a$, a contradiction to the
assumption that we made at the very beginning of this section.

\medskip
This completes the proof of the lemma.

\section{Proof of Lemma~\ref{lem:10}} \label{sec:5}

We want to use Proposition~\ref{prop:dworkcongruence} with
$A_r(m)=g_r(m)=B_{\mathbf N}(m)$ for all $r$. 
Clearly, the proposition would imply
that $S(a,K,s,p,m)\in p^{s+1}
B_{\mathbf N}(m)\mathbb Z_p$, and, thus, the claim.
So, we need to
verify the conditions $(i)$--$(iii)$ in the statement of the proposition.

Condition~$(i)$ is true since $B_{\mathbf N}(0)=1$. Condition~$(ii)$ holds by the definitions of $A_r(m)$ and $g_r(m)$.
To check that Condition~$(iii)$ holds is more complicated.
The proof will be decomposed in three steps. 

\subsection{First step} Straightforward computations imply that, 
for any fixed $j\in\{1, 2, \ldots, k\}$, we have 
\begin{align}
&\frac{B_{N_j}(v+up+np^{s+1})}{B_{N_j}(up+np^{s+1})} 
= \frac{(N_jv+N_jup+N_jnp^{s+1})!}{(v+up+np^{s+1})!^{N_j}}\cdot 
\frac{(up+np^{s+1})!^{N_j}}{(N_jup+N_jnp^{s+1})!}\nonumber 
\\
&\kern1cm
= \frac{(N_jv+N_jup+N_jnp^{s+1})(N_jv-1+N_jup+N_jnp^{s+1})\cdots
(N_jup+1+N_jnp^{s+1})} 
{(v+up+np^{s+1})^{N_j}(v-1+up+np^{s+1})^{N_j}\cdots
(1+up+np^{s+1})^{N_j}} \nonumber 
\\
&\kern1cm
= \frac{(N_jv+N_jup)(N_jv-1+N_jup)\cdots (1+N_jup)+\mathcal{O}(p^{s+1})}
{(v+up)^{N_j}(v-1+up)^{N_j}\cdots (1+up)^{N_j}+\mathcal{O}(p^{s+1})} ,
\nonumber
\end{align}
where $\mathcal{O}(R)$ denotes an element of $R\mathbb Z_p$.
Now we claim that this implies
\begin{equation} \label{eq:1} 
\frac{B_{N_j}(v+up+np^{s+1})}{B_{N_j}(up+np^{s+1})} 
= \frac{B_{N_j}(v+up)}{B_{N_j}(up)}+\mathcal{O}(p^{s+1}). 
\end{equation}
Indeed, if $v=0$, then this holds trivially.  
If $v>0$, then, together with the hypothesis $ v<p$, we infer that
$(v+up)(v-1+up)\cdots (1+up)$  
is not divisible by $p$, and thus we have 
\begin{multline*}
 \frac{1}
{\big((v+up)(v-1+up)\cdots (1+up)\big)^{N_j}+\mathcal{O}(p^{s+1})} \\
=
\frac{1}{\big((v+up)(v-1+up)^{N_j}\cdots (1+up)\big)^{N_j}} \, \big(1+
\mathcal{O}(p^{s+1}) \big). 
\end{multline*}
Hence, 
\begin{multline*}
 \frac{(N_jv+N_jup)(N_jv-1+N_jup)\cdots (1+N_jup)+\mathcal{O}(p^{s+1})}
{(v+up)^{N_j}(v-1+up)^{N_j}\cdots (1+up)^{N_j}+\mathcal{O}(p^{s+1})} 
\\
=  \frac{(N_jv+N_jup)(N_jv-1+N_jup)\cdots (1+N_jup)}
{(v+up)^{N_j}(v-1+up)^{N_j}\cdots (1+up)^{N_j}} \\
+ \frac{\mathcal{O}(p^{s+1})}{(v+up)^{N_j}(v-1+up)^{N_j}\cdots (1+up)^{N_j}}, 
\end{multline*}
which proves~\eqref{eq:1} because 
\begin{equation} \label{eq:vup}
\frac{1}{(v+up)(v-1+up)\cdots (1+up)}\in \mathbb{Z}_p
\end{equation}
and  
\begin{equation}\label{eq:fix2}
 \frac{(N_jv+N_jup)(N_jv-1+N_jup)\cdots (1+N_jup)}
{(v+up)^{N_j}(v-1+up)^{N_j}\cdots (1+up)^{N_j}} = \frac{B_{N_j}(v+up)}{B_{N_j}(up)}.
\end{equation}
Note that 
\eqref{eq:vup} and \eqref{eq:fix2} also imply 
that $B_{N_j}(v+up)/B_{N_j}(up)\in
\mathbb{Z}_p$, a property that will be used below. 

We deduce from~\eqref{eq:1} that 
\begin{equation}\label{ex:fix1}
\prod_{j=1}^k \frac{B_{N_j}(v+up+np^{s+1})}{B_{N_j}(up+np^{s+1})} =
\prod_{j=1}^k\left(\frac{B_{N_j}(v+up)}{B_{N_j}(up)}+\mathcal{O}(p^{s+1})\right).
\end{equation}
By expanding the product on the right-hand side of~\eqref{ex:fix1} and
using that 
$$B_{N_j}(v+up)/B_{N_j}(up)\in \mathbb{Z}_p,$$
we obtain our first main equality: 
\begin{equation}\label{eq:fix3}
\frac{B_{\mathbf N}(v+up+np^{s+1})}{B_{\mathbf N}(up+np^{s+1})} = 
\frac{B_{\mathbf N}(v+up)}{B_{\mathbf N}(up)}+\mathcal{O}(p^{s+1}).
\end{equation}

\medskip

\subsection{Second step}
The properties of $\Gamma_p$ imply that 
\begin{align}
\frac{B_{N_j}(up+np^{s+1})}{B_{N_j}(u+np^s)} &= (-1)^{N_j-1}
\frac{\Gamma_p(1+N_jup+N_jnp^{s+1})}{\Gamma_p(1+up+np^{s+1})^{N_j}}
\label{eq:firstequality} \\
&= (-1)^{N_j-1} \frac{\Gamma_p(1+N_jup)+\mathcal{O}(p^{s+1})}
{\Gamma_p(1+up)^{N_j}+\mathcal{O}(p^{s+1})}
\label{eq:rajoutcorrectionbis} \\
&= (-1)^{N_j-1} \frac{\Gamma_p(1+N_jup)}
{\Gamma_p(1+up)^{N_j}}\,\big(1+\mathcal{O}(p^{s+1})\big)
 \label{eq:rajoutcorrection}\\
&= \frac{B_{N_j}(up)}{B_{N_j}(u)}\,\big(1+\mathcal{O}(p^{s+1})\big) ,
\label{eq:2}
\end{align}
where $(i)$ of Lemma~\ref{lem:gammap} is used 
to see~\eqref{eq:firstequality} and~\eqref{eq:2}, and 
$(ii)$ is used for~\eqref{eq:rajoutcorrectionbis}.
Equation~\eqref{eq:rajoutcorrection} holds because $\Gamma_p(1+up)$
and $\Gamma_p(1+N_jup)$  
are both not divisible by $p$. Taking the product over $j=1, 2, \dots,
k$, we obtain from \eqref{eq:2} our second main  
equality: 
\begin{equation}\label{eq:fix4}
\frac{B_{\mathbf N}(up+np^{s+1})}{B_{\mathbf N}(u+np^s)}
=\frac{B_{\mathbf N}(up)}{B_{\mathbf N}(u)}\,\big(1+\mathcal{O}(p^{s+1})\big).
\end{equation}

\subsection{Third step}
We now multiply the right-hand and left-hand sides
of the main equalities~\eqref{eq:fix3} and 
\eqref{eq:fix4}: we get after simplification
\begin{equation*}
\frac{B_{\mathbf N}(v+up+np^{s+1})}{B_{\mathbf N}(u+np^s)} = 
\frac{B_{\mathbf N}(v+up)}{B_{\mathbf N}(u)}\,\big(1+\mathcal{O}(p^{s+1})\big) 
+ \frac{B_{\mathbf N}(up)}{B_{\mathbf N}(u)} \,\mathcal{O}(p^{s+1}).
\end{equation*}
We can rewrite this as 
\begin{align}
\frac{B_{\mathbf N}(v+up+np^{s+1})}{B_{\mathbf N}(v+up)}
&=\frac{B_{\mathbf N}(u+np^s)}{B_{\mathbf N}(u)}\,\big(1+\mathcal{O}(p^{s+1})\big)
+\frac{B_{\mathbf N}(up)}{B_{\mathbf N}(u)}\cdot\frac{B_{\mathbf N}(u+np^s)}{B_{\mathbf N}(v+up)} \,\mathcal{O}(p^{s+1})
\notag
\\
&=\frac{B_{\mathbf N}(u+np^s)}{B_{\mathbf N}(u)}+\frac{B_{\mathbf N}(u+np^s)}{B_{\mathbf N}(u)}\,\mathcal{O}(p^{s+1})
+\frac{B_{\mathbf N}(u+np^s)}{B_{\mathbf N}(v+up)} \,\mathcal{O}(p^{s+1}),
\label{eq:quasifinale}
\end{align}
where the last line holds because $v_p\big(B_{\mathbf N}(up)/B_{\mathbf N}(u)\big)=0$.

If we compare \eqref{eq:DworkA} (with $A_r(m)=g_r(m)=B_{\mathbf N}(m)$) 
and \eqref{eq:quasifinale}, we see that it only
remains to prove that we have  
\begin{equation} 
\label{eq:Bg} 
\frac{B_{\mathbf N}(u+np^s)}{B_{\mathbf N}(u)} \in  \frac{B_{\mathbf N}(n)}{B_{\mathbf N}(v+up)}\,\mathbb{Z}_p 
\quad 
\textup{and}
\quad
\frac{B_{\mathbf N}(u+np^s)}{B_{\mathbf N}(v+up)} \in  \frac{B_{\mathbf N}(n)}{B_{\mathbf N}(v+up)}\,\mathbb{Z}_p.
\end{equation}
The first membership relation 
follows from Lemma~\ref{lem:6} because the latter implies the stronger property
\begin{equation}\label{eq:fix5}
\frac{B_{\mathbf N}(u+np^s)}{B_{\mathbf N}(u)} 
\in \bigg(\prod_{j=1}^k B_{N_j}(n)\bigg) \,\mathbb{Z}_p = B_{\mathbf N}(n)\,\mathbb{Z}_p.
\end{equation}
For the second membership relation, we note that 
\begin{align*}
\frac{B_{\mathbf N}(u+np^s)}{B_{\mathbf N}(v+up)} &= B_{\mathbf N}(u) \cdot \frac{B_{\mathbf N}(u+np^s)}{B_{\mathbf N}(u)} \cdot
\frac{1}{B_{\mathbf N}(v+up)}  \\
&\kern1cm
\in B_{\mathbf N}(u) \,\frac{B_{\mathbf N}(n)}{B_{\mathbf N}(v+up)} \,\mathbb{Z}_p \subset
\frac{B_{\mathbf N}(n)}{B_{\mathbf N}(v+up)} \, \mathbb{Z}_p,
\end{align*}
where we have used~\eqref{eq:fix5}.

It now follows from~\eqref{eq:quasifinale} and~\eqref{eq:Bg} that we have 
$$\frac{B_{\mathbf N}(v+up+np^{s+1})}{B_{\mathbf N}(v+up)} - 
\frac{B_{\mathbf N}(u+np^s)}{B_{\mathbf N}(u)} \in  \frac{B_{\mathbf N}(n)}{B_{\mathbf N}(v+up)}\,\mathbb{Z}_p,
$$
i.e., we have checked that hypothesis $(iii)$ in
Proposition~\ref{prop:dworkcongruence} holds in our situation.  
We can therefore apply this proposition and obtain exactly the
statement of the lemma.

\section{Proof of Lemma~\ref{lem:11}} \label{sec:6}

The claim is trivially true if $p$ divides $m$.
We may therefore assume that $p$ does not divide $m$ for the rest of the proof.
Let us write $m=a+pj$, with $0< a<p$. Then
comparison with \eqref{eq:congrconj1} shows that we are in a
very similar situation here. Indeed, we may derive \eqref{eq:110}
from Lemma~\ref{lem:12}. In order to see this, we observe that
\begin{align*}
H_{Lmp^s}-H_{L\fl{m/p}p^{s+1}}&=
\sum _{\ep=1} ^{Lap^s}\frac {1} {Ljp^{s+1}+\ep}\\
&=
\sum _{\ep=1} ^{\fl{La/p}}\frac {1} {Ljp^{s+1}+\ep p^{s+1}}
+
\underset{p^{s+1}\nmid \ep}{\sum _{\ep=1} ^{Lap^s}}\frac {1}
{Ljp^{s+1}+\ep}\\
&=\frac {1} {p^{s+1}}(H_{Lj+\fl{La/p}}-H_{Lj})+
\underset{p^{s+1}\nmid \ep}{\sum _{\ep=1} ^{Lap^s}}\frac {1}
{Ljp^{s+1}+\ep}.
\end{align*}
Because of $v_p(x+y)\ge\min\{v_p(x),v_p(y)\}$, this implies
$$
v_p(H_{Lmp^s}-H_{L\fl{m/p}p^{s+1}})\ge
\min\{ -1-s+v_p(H_{Lj+\fl{La/p}}-H_{Lj}),-s\}.
$$
It follows that
\begin{multline*}
v_p\Big(B_{\mathbf N}(m)\big(H_{Lmp^s}-H_{L\fl{m/p}p^{s+1}}\big)\Big)\\
\ge
-1-s+\min\left\{v_p\Big(B_{\mathbf N}(a+pj)(H_{Lj+\fl{La/p}}-H_{Lj})\big),
1+v_p\big(B_{\mathbf N}(a+pj)\Big)\right\}.
\end{multline*}
Use of Lemmas~\ref{lem:multinomial/N!} and \ref{lem:12} 
then completes the proof.

\section{Outline of the proof of Theorem~\ref{thm:3}}
\label{sec:Xi}

In this section, we provide a brief outline of the proof of
Theorem~\ref{thm:3}, reducing it to several lemmas and their corollaries.
These are subsequently proved in the next section.

We must slightly ``upgrade" the proof of Theorem~\ref{thm:2}
in the special case that $\mathbf N=(N,N,\dots,N)$ (with $k$
occurrences of $N$) and $L=N$. To be
precise, for all primes $p$, and for all 
non-negative integers $K$, $a$, and $j$ with $0\le a<p$, we must prove 
\begin{equation} \label{eq:CcongU}
C(a+Kp)=\sum_{j=0}^K B_{\mathbf N}(a+jp)B_{\mathbf N}(K-j)
(H_{N(K-j)}-pH_{Na+Njp})  \in p \Xi_N N!^k\mathbb{Z}_p,
\end{equation}
where $B_{\mathbf N}(m)=\frac{(Nm)!^k}{m!^{kN}}$, 
instead of ``just" \eqref{eq:Ccong} for the present choice of
parameters. We recall that, in view of a comparison of
\eqref{eq:ThL} with the definition of $\Xi_N$, there is only something
to prove if $v_p(H_N)>0$. 
Hence, Lemma~\ref{lem:H_L} tells that there is nothing to prove if
$p=2$, and Lemma~\ref{lem:3} together with Remarks~\ref{rem:Xi7}(a) 
in the Introduction tells that, if $p=3$, only the case
$N=22$ need to be considered. 
Therefore, in the remainder of this section, we always assume that $p$
is a prime with 
$3\le p\le N$ and $v_p(H_N)>0$, and if $p=3$ then $N=22$.

There are two cases which can be treated directly:
if $K=a=0$, then $C(0)=0$, and thus \eqref{eq:CcongU} holds trivially,
whereas if $K=0$ and $a=1$, then $C(1)=-pB_{\mathbf
N}(1)H_N=-pN!^kH_N$ and thus \eqref{eq:CcongU} holds by definition of
$\Xi_N$. We therefore assume 
in addition $a+Kp\ge2$ for the remainder of this section.

Going through the outline of the proof of Theorem~\ref{thm:2} in
Section~\ref{sec:2}, we see that, in order to establish \eqref{eq:CcongU},
we need to prove corresponding stronger versions of
\eqref{eq:firstreduction} and Lemmas~\ref{lem:12} and \ref{lem:11}. 
To be precise, given
non-negative integers $m$, $K$, $a$, and $j$ with $0\le a<p$ and
$a+Kp\ge2$, for $p\ge5$ we should prove
\begin{equation}\label{eq:firstreductionU}
C(a+Kp) \equiv \sum_{j=0}^K B_{\mathbf N}(a+jp)B_{\mathbf N}(K-j) (H_{N(K-j)}-H_{\lfloor
Na/p\rfloor+Nj}) \mod p^3 N!^k\mathbb{Z}_p,
\end{equation}
respectively, for $p=3$ and $N=22$,
\begin{equation}\label{eq:firstreductionU3}
C(a+3K) \equiv \sum_{j=0}^K B_{\mathbf N}(a+3j)B_{\mathbf N}(K-j) (H_{N(K-j)}-H_{\lfloor
Na/3\rfloor+Nj}) \mod 3^2 N!^k\mathbb{Z}_3,
\end{equation}
respectively, if $v_p(H_N)\ge3$ and $p$ a Wolstenholme prime or $p\mid N$, 
\begin{equation}\label{eq:firstreductionUW}
C(a+Kp) \equiv \sum_{j=0}^K B_{\mathbf N}(a+jp)B_{\mathbf N}(K-j) (H_{N(K-j)}-H_{\lfloor
Na/p\rfloor+Nj}) \mod p^4 N!^k\mathbb{Z}_p,
\end{equation}
and we need to prove
\begin{equation} \label{eq:congrconj1U}
B_{\mathbf N}(a+pj)\left(H_{Nj+\lfloor Na/p\rfloor} - H_{Nj}\right)
\in p \Xi_N N!^k \mathbb{Z}_p
\end{equation}
and
\begin{equation} \label{eq:110U}
B_{\mathbf N}(m)\big(H_{Nmp^s}-H_{N\fl{m/p}p^{s+1}}\big)\in 
p^{-s}\Xi_N N!^k \mathbb Z_p\ .
\end{equation}
If these five relations are used in Section~\ref{sec:2} instead of
their weaker counterparts 
\eqref{eq:firstreduction}, \eqref{eq:congrconj1}, and
\eqref{eq:110}, respectively, then the proof in Section~\ref{sec:2} 
yields \eqref{eq:CcongU}, as required.

\subsection{Proof of
\eqref{eq:firstreductionU}--\eqref{eq:firstreductionUW}}
We recall that the congruence \eqref{eq:firstreduction} followed from
the congruence \eqref{eq:J}. Clearly, one cannot hope for
improving \eqref{eq:J} to $pH_{J} \equiv H_{\fl{J/p}} \text{ mod }
p^3\mathbb{Z}_p,$ there are counterexamples. However, if we combine
Corollary~\ref{cor:C1} and Lemma~\ref{lem:p<L}(1) with $L=N$, 
then we see that
$$
B_{\mathbf N}(a+jp)B_{\mathbf N}(K-j)pH_{Na+Njp} 
\equiv B_{\mathbf N}(a+jp)B_{\mathbf N}(K-j)H_{\fl{Na/p}+Nj} \mod
p^3N!^k\mathbb{Z}_p
$$
as long as $a+Kp\ge2$ and $a\ne0$. Moreover, 
due to
Lemma~\ref{lem:multinomial/N!} and Corollary~\ref{cor:W2}, the above
congruence is even true if $a=0$
and $p\ge5$.
This implies \eqref{eq:firstreductionU}.
The congruence \eqref{eq:firstreductionU3} follows in the same way by
the slightly weaker assertion for $p=3$ in Corollary~\ref{cor:W2}.

For the congruence \eqref{eq:firstreductionUW}, one needs to combine
Corollary~\ref{cor:C1} and Lemma~\ref{lem:p<L}(4) with $L=N$, to see
that
\begin{equation} \label{eq:congW}
B_{\mathbf N}(a+jp)B_{\mathbf N}(K-j)pH_{Na+Njp} 
\equiv B_{\mathbf N}(a+jp)B_{\mathbf N}(K-j)H_{\fl{Na/p}+Nj} \mod
p^4N!^k\mathbb{Z}_p
\end{equation}
as long as $a+Kp\ge2$ and $a\ne0$. 
If $a=0$, then the congruence \eqref{eq:congW}
still holds as long as $j<K$ because of 
Lemma~\ref{lem:multinomial/N!}, Corollary~\ref{cor:W2}, and the fact
that the term $B_{\mathbf N}(K-j)$ contributes at least one factor $p$. 
The only remaining case to be discussed
is $a=0$ and $j=K$.
If we apply the simple observation that
$v_p(B_{\mathbf N}(p^eh))=v_p(B_{\mathbf N}(h))$ for any positive
integers $e$ and $h$
to $B_{\mathbf N}(a+jp)=B_{\mathbf N}(jp)$, then, making again appeal to
Corollary~\ref{cor:C1} and Lemma~\ref{lem:p<L}(4) with $L=N$,
we see that \eqref{eq:congW} holds as well
as long as $j=K$ is no prime power. Finally, let $a=0$ and $j=K$ be a prime
power, $j=K=p^e$ say. If $e\ge1$, then we may use Lemma~\ref{lem:W3}
with $J=a+jp=a+Kp=p^{e+1}$ and $v_p(B_{\mathbf N}(p^{e+1}))=v_p(B_{\mathbf
N}(1))=v_p(N!^k)$ to conclude that \eqref{eq:congW} also holds in this case.
On the other hand, if $j=K=1$, then Lemma~\ref{lem:congH}, the
assumption that $p$ is a Wolstenholme prime or that $p$ divides $N$,
the fact that $v_p(B_{\mathbf N}(p))=v_p(B_{\mathbf
N}(1))=v_p(N!^k)$,
altogether yield the congruence \eqref{eq:congW} in this case as well.
 
\subsection{Proof of \eqref{eq:congrconj1U}}
The proof of Lemma~\ref{lem:12} in Section~\ref{sec:4} can be used
verbatim. The place where the necessary improvement is possible is
\eqref{eq:log} (which was used there with $k=\fl{\log_pL}$). 
Clearly, we have $\fl{22/3}\ge\fl{\log_3\fl{22/3}}+5$. Moreover,
we claim that in our more special context we have
\begin{equation} \label{eq:logU5}
\fl{N/5}\ge\fl{\log_5\fl{N/5}}+4 
\end{equation}
if $p=5$, we have
\begin{equation} \label{eq:logUp}
\fl{N/p}\ge\fl{\log_p\fl{N/p}}+5
\end{equation}
if $p>5$, and we have
\begin{equation} \label{eq:logUW}
\fl{N/p}\ge\fl{\log_p\fl{N/p}}+6
\end{equation}
if $v_p(H_N)>2$.
Indeed, if $p=5$ then Lemma~\ref{lem:5} says that $N$ must be at
least $20$ because otherwise $v_5(H_N)\le 0$. The inequality
\eqref{eq:logU5} follows immediately. On the other hand, if $p=11$
then it can be checked 
(by our table mentioned in Remarks~\ref{rem:Xi7}(c) in the
Introduction, for example) that $N$ must be at least $77$ because otherwise
$v_{11}(H_N)\le 0$. If $p$ is different from $5$ and $11$ (in
addition to being at least $5$) then, according to
Lemma~\ref{lem:p<L}(3), we must have $N/p\ge5$. In both of the latter
cases, the inequality \eqref{eq:logUp} follows immediately.
Finally, if $v_p(H_N)>2$ then, according to
Lemmas~\ref{lem:5} and \ref{lem:p<L}(4), we must have 
$p>5$ and $N/p\ge6$, thus implying \eqref{eq:logUW}.

The effect of the above improvement over \eqref{eq:log} is that the
proof in Section~\ref{sec:4} shows that in our more special context we
have
$$
B_{\mathbf N}(a+pj)\left(H_{Nj+\lfloor Na/p\rfloor} - H_{Nj}\right)
\in \begin{cases} p^2N!^k\mathbb{Z}_p&\text{if }p=5,\\
p^4N!^k\mathbb{Z}_p&\text{if }v_p(H_N)>2,\\
p^3N!^k\mathbb{Z}_p&\text{otherwise.}\end{cases}
$$
This implies \eqref{eq:congrconj1U} since $v_5(H_N)\le 1$ according to
Lemma~\ref{lem:5}. 

\subsection{Proof of \eqref{eq:110U}} \label{sec:m=1}
Again, the proof of Lemma~\ref{lem:11} in Section~\ref{sec:6} can be
used verbatim. The only differences  
arise at the very end. Namely, we have
to use Corollary~\ref{cor:C1} combined with Lemma~\ref{lem:p<L}(1) 
instead of Lemma~\ref{lem:multinomial/N!},
and \eqref{eq:firstreductionU} instead of Lemma~\ref{lem:12}, to obtain
our assertion for $m\ge2$. In the remaining case $m=1$, we compute directly:
\begin{align*}
B_{\mathbf N}(1)H_{Np^s}
&=N!^k\Bigg(\frac {1} {p^s}H_N+
\frac {1} {p^{s-1}}
\underset{p\nmid\ep}{\sum _{\ep=1} ^{Np}}\frac {1} {\ep}+
\frac {1} {p^{s-2}}
\underset{p\nmid\ep}{\sum _{\ep=1} ^{Np^2}}\frac {1} {\ep}+
\underset{p^{s-2}\nmid \ep}{\sum _{\ep=1} ^{Np^s}}\frac {1} {\ep}\Bigg).
\end{align*}
The last sum over $\ep$ is clearly an element of
$p^{-s+3}\mathbb Z_p$. Moreover, 
if $p\ge5$, Lemma~\ref{lem:W1} implies that the
next-to-last expression between parentheses 
is an element of $p^{-s+4}\mathbb Z_p\subset p^{-s+3}\mathbb Z_p$,
and that the second expression between parentheses is an element of
$p^{-s+3}\mathbb Z_p$.
If $p=3$, it can be checked directly that 
the expression between parentheses is an element of $3^{-s+2}\mathbb Z_3$.
Putting everything together, we conclude that
$$
B_{\mathbf N}(1)H_{Np^s}\in p^{-s} N!^k\Xi_N\mathbb Z_p,
$$
which finishes the proof of \eqref{eq:110U}.

\section{More auxiliary lemmas}
\label{sec:aux}

In this section, we prove the auxiliary results necessary 
for the proof of Theorem~\ref{thm:3}, of which the outline was given in the
previous section. These are, on the one hand,
improvements of Lemma~\ref{lem:multinomial/N!}
(see Lemmas~\ref{lem:B1}, \ref{lem:B2} and Corollary~\ref{cor:C1}), 
and, on the other hand, assertions addressing specific $p$-adic properties
of harmonic numbers 
(see Lemmas~\ref{lem:H_L}--\ref{lem:congH} and Corollary~\ref{cor:W2}).
Some of the results of this section are also referred to in the next section.

We begin by two lemmas improving on Lemma~\ref{lem:multinomial/N!}.

\begin{lem} \label{lem:B1}
For all positive integers $N$, $a$, and primes $p$ with $2\le a<p$, we
have
$$v_p(B_{\mathbf N}(a))\ge \fl{\frac {N} {p}}+v_p(N!^k),$$
where $\mathbf N=(N,N,\dots,N)$ {\em(}with $k$
occurrences of $N${\em)}.
\end{lem}

\begin{proof}
This was implicitly proved by the estimations leading to
\eqref{eq:unglA0} in the case that $N_i=N$ for all $i$. 
\end{proof}

\begin{lem} \label{lem:B2}
For all positive integers $N$, $a$, $j$, and primes $p$ with $1\le a<p$, we
have
$$v_p(B_{\mathbf N}(a+jp))\ge \fl{\frac {N} {p}}+\min\{1+T_1,T_2\}-1+
v_p(N!^k),$$
where $\mathbf N=(N,N,\dots,N)$ {\em(}with $k$
occurrences of $N${\em)}, and where
$$T_1=\max_{1\le\ep\le \lfloor
Na/p\rfloor}v_p(Nj+\ep)\quad \text{and}\quad  T_2=\fl{\log_p (a+pj)}.$$
\end{lem}

\begin{proof}
This is seen by going through the estimations leading to
\eqref{eq:ungl7} in the case that $N_i=N$ for all $i$,
without employing \eqref{eq:ungl6a}, \eqref{eq:ungl100},
and \eqref{eq:ungl6b}.
\end{proof}

As a corollary to Lemmas~\ref{lem:B1} and \ref{lem:B2}, we obtain the
following succinct $p$-adic estimation for $B_{\mathbf N}(m)$,
which is needed in the proofs of \eqref{eq:firstreductionU},
\eqref{eq:firstreductionUW}, and \eqref{eq:110U}.

\begin{coro} \label{cor:C1}
Let $N$ and $m$ be positive integers and $p$ be a prime such that
$m$ is at least $2$ and not divisible by $p$. Then we
have
$$v_p(B_{\mathbf N}(m))\ge \fl{\frac {N} {p}}+v_p(N!^k),$$
where $\mathbf N=(N,N,\dots,N)$ {\em(}with $k$
occurrences of $N${\em)}.
\end{coro}

The next three lemmas of this section 
provide elementary information on the $p$-adic
valuation of harmonic numbers for $p=2,3,5$ which is needed
in the proof of Lemma~\ref{lem:p<L} and is also referred to
frequently at other places. 
(For example, Lemma~\ref{lem:H_L} was used in the proof of~\eqref{eq:unglB}.)
The proofs are not
difficult (cf.\ \cite{boyd}) and are therefore omitted. 

\begin{lem} \label{lem:H_L}
For all positive integers $L$, we have $v_2(H_L)=-\fl{\log_2L}$.
\end{lem}

\begin{lem} \label{lem:3}
We have $v_3(H_2)=v_3(H_7)=v_3(H_{22})=1$.
For positive integers $L\notin\{2,7,22\}$, we have $v_3(H_L)\le0$.
\end{lem}

\begin{lem} \label{lem:5}
We have $v_5(H_4)=2$ and $v_5(H_{20})=v_5(H_{24})=1$.
For positive integers $L\notin\{4,20,24\}$, we have $v_5(H_L)\le0$.
\end{lem}

Next, we record some properties of integers $L$ and primes $p$
for which $v_p(H_L)>0$. These are needed throughout
Section~\ref{sec:Xi}. 

\begin{lem} \label{lem:p<L}
Let $p$ be a prime, and let $L$ be an integer with $p\le L$.
Then the following assertions hold true:

\begin{enumerate} 
\item If $v_p(H_L)>0$ then $L\ge 2p$.
\item If $v_p(H_L)>0$ and $p\ne 3$ then $L\ge 3p$.
\item If $v_p(H_L)>0$ and $p\notin\{3,5,11\}$ then $L\ge 5p$.
\item If $v_p(H_L)>2$ then $L\ge 6p$.
\end{enumerate}
\end{lem}

\begin{proof}
We have
$$H_L=\frac {1} {p}H_{\fl{L/p}}+
\underset{p\nmid\ep}{\sum _{\ep=1} ^{L}}\frac {1} {\ep}.$$
Since the sum over $\ep$ is in $\mathbb Z_p$, in order to have
$v_p(H_L)>0$ we must have $v_p(H_{\fl{L/p}})>0$. 
Clearly, $v_p(H_1)=0$ so that (1) follows. If $\fl{L/p}<3$ and $p\ne3$, 
then $v_p(H_{\fl{L/p}})$ cannot be positive since
$H_1=1$ and $H_2=\frac {3} {2}$. This implies (2). 
If $\fl{L/p}<5$ and $p\notin \{3,5,11\}$, 
then, again, $v_p(H_{\fl{L/p}})$ cannot be positive since, as we
already noted, $H_1=1$ and $H_2=\frac {3} {2}$, and since
$H_3=\frac {11} {6}$ and $H_4=\frac {25} {12}.$
This yields (3).

To see (4), we observe that, owing to
Lemmas~\ref{lem:H_L}--\ref{lem:5}, we may assume that
$p\notin\{2,3,5\}$. Furthermore, if $p=11$ then, according to our
table referred to in Remarks~\ref{rem:Xi7}(c) in the Introduction
(see also \cite{boyd}), 
we have $L\ge848$. Similarly, if $p=137$ then, according to
our table, we have $L>500000$. The claim can now be established in
the style of the proofs of (1)--(3) upon observing
that $H_5=\frac {137} {60}$.
\end{proof}

We turn to a slight generalisation of Wolstenholme's theorem
on harmonic numbers. 
(We refer the reader to \cite[Chapter VII]{hw} for information on 
this theorem, which corresponds to the case $r=1$ in the lemma
below.) 

\begin{lem} \label{lem:W1}
For all primes $p\ge5$ and positive integers $r$, we have
$$v_p(H_{rp-1}-H_{rp-p})\ge2.$$
\end{lem}
\begin{proof}
By simple rearrangement, we have
\begin{align*}
H_{rp-1}-H_{rp-p}&=
\sum _{\ep=1} ^{p-1}\frac {1} {rp-p+\ep}=
\sum _{\ep=1} ^{(p-1)/2}\left(\frac {1} {rp-p+\ep}+\frac {1}
{rp-\ep}\right)\\
&=
p(2r-1)\sum _{\ep=1} ^{(p-1)/2}\frac {1} {(rp-p+\ep)(rp-\ep)}.
\end{align*}
It therefore suffices to consider the last sum over $\ep$ in 
$\mathbb Z/p\mathbb Z$
and show that it is $\equiv 0$~mod~$p$. If we reduce
this sum mod~$p$, then we are left with
$$-\sum _{\ep=1} ^{(p-1)/2}\frac {1} {\ep^2},$$
which is, up to the sign, the sum of all quadratic residues in
$\mathbb Z/p\mathbb Z$, 
that is, equivalently, 
$$-\sum _{\ep=1} ^{(p-1)/2}\ep^2=\frac {p(p-1)(p+1)} {24}.$$
Clearly, this is divisible by $p$ for all primes $p\ge5$.
\end{proof}

As a corollary, we obtain strengthenings of \eqref{eq:J} that we
need in the proof of \eqref{eq:firstreductionU} and
\eqref{eq:firstreductionU3}.

\begin{coro} \label{cor:W2}
For all primes $p\ge5$ and positive integers $J$ divisible by $p$, we have
$$
pH_{J} \equiv H_{{J/p}} \mod p^3\mathbb{Z}_p.
$$
Moreover, for all positive integers $J$ divisible by $3$, we have
$$
3H_{J} \equiv H_{{J/3}} \mod 3^2\mathbb{Z}_3.
$$
\end{coro}
\begin{proof}
By simple rearrangement, we have
$$pH_{J} - H_{{J/p}}= 
p\sum _{r=1} ^{J/p}(H_{rp-1}-H_{rp-p}).$$
Due to Lemma~\ref{lem:W1}, the $p$-adic valuation of this
expression is at least $3$
if $p\ge5$. If $p=3$, it is easily seen
directly that the $3$-adic valuation of this expression is at least $2$.
\end{proof}

Further strengthenings of \eqref{eq:J}, needed 
in the proof of \eqref{eq:firstreductionUW}, are given in the final
two lemmas of this section.

\begin{lem} \label{lem:W3}
For all primes $p\ge5$ and positive integers $J$ divisible by $p^2$, we have
$$
pH_{J} \equiv H_{{J/p}} \mod p^5\mathbb{Z}_p.
$$
\end{lem}
\begin{proof}
Again, by simple rearrangement, we have
\begin{align*}
pH_{J}-H_{J/p}&=p
\underset{p\nmid\ep}{\sum _{\ep=1} ^{J-1}}\frac {1} {\ep}
=p\sum _{r=1} ^{J/p^2}
\underset{p\nmid\ep}{\sum _{\ep=1} ^{p^2-1}}\frac {1} {rp^2-p^2+\ep}\\
&=p^3\sum _{r=1} ^{J/p^2}(2r-1)
\underset{p\nmid\ep}{\sum _{\ep=1} ^{(p^2-1)/2}}
\frac {1} {(rp^2-p^2+\ep)(rp^2-\ep)}.
\end{align*}
It therefore suffices to consider the last sum over $\ep$ in 
$\mathbb Z/p^2\mathbb Z$
and show that it is $\equiv 0$~mod~$p^2$. If we reduce
this sum mod~$p^2$, then we are left with
$$-\underset{p\nmid\ep}{\sum _{\ep=1} ^{(p^2-1)/2}}\frac {1} {\ep^2},$$
which is, up to the sign, the sum of all quadratic residues in
$\mathbb Z/p^2\mathbb Z$, 
that is, equivalently, 
$$-\underset{p\nmid\ep}{\sum _{\ep=1} ^{(p^2-1)/2}}\ep^2
=-\sum _{\ep=1} ^{(p^2-1)/2}\ep^2
+\sum _{\ep=1} ^{(p-1)/2}(p\ep)^2
=-\frac {p^2(p^2-1)(p^2+1)} {24}+p^2
\frac {p(p-1)(p+1)} {24}.$$
Clearly, this is divisible by $p^2$ for all primes $p\ge5$.
\end{proof}

\begin{lem} \label{lem:congH}
For all primes $p\ge5$ and positive integers $N$, we have
\begin{equation} \label{eq:congH1}
pH_{pN}\equiv H_N\mod p^4\mathbb Z_p
\end{equation}
if and only if $p$ is a Wolstenholme prime or $p$ divides $N$.
\end{lem}
\begin{proof}
Using a rearrangement in the spirit of Lemma~\ref{lem:W1}, we obtain
$$pH_{pN}-H_N=
p\sum _{r=1} ^{N}\sum _{\ep=1} ^{p-1}\frac {1} {rp-p+\ep}.
$$
We consider the sum over $r$ in $\mathbb Z/p^3\mathbb Z$. This leads
to
\begin{align} \notag
\sum _{r=1} ^{N}\sum _{\ep=1} ^{p-1}(rp-p+\ep)^{-1}
&\equiv
\sum _{r=1} ^{N}\sum _{\ep=1} ^{p-1}\ep^{-1}\big(1+p(r-1)\ep^{-1}\big)^{-1}\\
\notag
&\equiv
\sum _{r=1} ^{N}\sum _{\ep=1} ^{p-1}
\left(\ep^{-1}-p(r-1)\ep^{-2}+p^2(r-1)^2\ep^{-3}\right)\\
&\equiv NH_{p-1}-p\binom N2H_{p-1}^{(2)}+p^2\frac {N(N-1)(2N-1)}
{6}H_{p-1}^{(3)}\kern-5pt\mod \mathbb Z/p^3\mathbb Z,
\label{eq:p3a}
\end{align}
where $H_{m}^{(\al)}$ denotes the higher harmonic number defined by
$H_{m}^{(\al)}=\sum_{n=1}^{m} \frac{1}{n^\al}$. By a rearrangement
analogous to the one in the proof of Lemma~\ref{lem:W1}, one sees
that $v_p(H_{p-1}^{(3)})\ge1$, whence we may disregard the last term
in the last line of \eqref{eq:p3a}. As it turns out, $H_{p-1}$ and
$H_{p-1}^{(2)}$ are directly related modulo $\mathbb Z/p^3\mathbb Z$.
Namely, we have
\begin{align} 
\notag
H_{p-1}&\equiv\sum _{\ep=1} ^{p-1}\ep^{-1}
\equiv 2^{-1}\sum _{\ep=1} ^{p-1}\left(\ep^{-1} + (p-\ep)^{-1}\right)
\equiv p2^{-1}\sum _{\ep=1} ^{p-1}\big(\ep(p-\ep)\big)^{-1}\\
\notag
&\equiv -p2^{-1}\sum _{\ep=1} ^{p-1}\ep^{-2}\big(1+p\ep^{-1}\big)
\equiv -p2^{-1}\sum _{\ep=1} ^{p-1}\big(\ep^{-2}+p\ep^{-3}\big)\\
&\equiv -p2^{-1}H_{p-1}^{(2)}-p^22^{-1}H_{p-1}^{(3)}\mod \mathbb Z/p^3\mathbb Z.
\label{eq:HnHn2}
\end{align}
Again, we may disregard the last term. If we substitute this
congruence in \eqref{eq:p3a}, then we obtain
$$
\sum _{r=1} ^{N}\sum _{\ep=1} ^{p-1}(rp-p+\ep)^{-1}
\equiv 
N^2H_{p-1}
\mod \mathbb Z/p^3\mathbb Z.
$$
Hence, the sum is congruent to zero mod~$\mathbb Z/p^3\mathbb Z$ if and only
if $N$ is divisible by $p$
(recall Lemma~\ref{lem:W1} with $r=1$) or $v_p(H_{p-1})\ge3$,
that is, if $p$ is a Wolstenholme prime.
This establishes the assertion of the lemma.
\end{proof}

\section{The sequence $(t_N)_{N\ge1}$}
\label{sec:DworkKont}

The purpose of this section is to report on the evidence for our
conjecture that the largest integer $t_N$ such that
$q_{N,N}(z)^{1/t_N}\in\mathbb Z[[z]]$ is given by $t_N=\Xi_NN!$, 
that is, that Theorem~\ref{thm:3} with $k=1$ is optimal.
We assume $k=1$ throughout this section.

We first prove that there cannot be any prime number $p$ larger than
$N$ which divides $t_N$.

\begin{prop} \label{prop:p>N}
Let $p$ be a prime number and $N$ a positive integer with $p>N$.
Then there exists a positive integer $a<p$ such that
\begin{equation} \label{eq:p>N} 
v_p\big(B_{N}(a)H_{Na}\big)=0,
\end{equation}
where $B_N(m)=\frac {(Nm)!} {m!^N}$.
In particular, $p$ does not divide $t_N$. 
\end{prop}

\begin{proof}
If $N=1$, we choose $a=1$ to obtain $B_{1}(1)H_{1}=1$. On the other
hand, if $N>1$, we choose $a$ to be the least integer such that $aN\ge p$.
Since then $a<p$, we have
$$v_p\big(B_{N}(a)H_{Na}\big)=
\sum _{\ell = 1}^{\infty}\left(\fl{\frac {Na} {p^\ell}}-
N\fl{\frac {a} {p^\ell}}\right)
+v_p(H_{Na})=1-1=0.$$

To see that $p$ cannot divide $t_N$, we observe that
we have $C(a)=-pB_{N}(a)H_{Na}$ (for the sum
$C(\cdot)$ defined in \eqref{eq:C} in the case $k=1$ and $L=N$) because $a<p$. 
Hence, the assertion \eqref{eq:p>N} can be
reformulated as $v_p\big(C(a)\big)=1$.
Since, because of $p>N$,
we have $v_p(M_N)=v_p(N!)=0$ and $v_p(\Xi_N)=0$, it follows that 
$$C(a)\notin p^2\Xi_N N!\, \mathbb Z_p.$$
This means that 
one cannot increase the exponent of $p$ in \eqref{eq:CcongU} 
(with $k=1$) in our special case, and thus $p$ cannot divide $t_N$.
\end{proof}

So, if we hope to improve Theorem~\ref{thm:3} with $k=1$, then it must be by
increasing exponents of prime numbers $p\le N$ in \eqref{eq:Xi}. 
It can be checked directly that the exponent of $3$ cannot be
increased if $N=7$. (The reader should recall Remarks~\ref{rem:Xi7}(a)
in the Introduction.)
According to Remarks~\ref{rem:Xi7}(b), an
improvement is therefore only possible if $v_p(H_N)>2$ for some $p\le N$. 
Lemmas~\ref{lem:H_L}--\ref{lem:5} in Section~\ref{sec:aux} tell
that this does not happen with $p=2,3,5$, so that the exponents of
$2,3,5$ cannot be improved. 
(The same conclusion can also be drawn from \cite{boyd} for many
other prime numbers, but so far not for $83$, for example.)

We already discussed in Remarks~\ref{rem:Xi7}(c) 
whether at all there are primes $p$ and integers $N$
with $p\le N$ and $v_p(H_N)\ge3$. We recall that, so far, only five
examples are known, four of them involving $p=11$. 

The final result of this section shows that, if $v_p(H_N)=3$, the exponent
of $p$ in the definition of $\Xi_N$ in \eqref{eq:Xi} cannot be increased
so that Theorem~\ref{thm:3} would still hold. 
(The reader should recall Remarks~\ref{rem:Xi7}(b).)
\begin{prop} \label{prop:vp=3}
Let $p$ be a prime number and $N$ a positive integer with $p\le N$
and $v_p(H_N)=3$. If $p$ is not a Wolstenholme prime and $p$ does not
divide $N$, then\break
$q_{N,N}(z)^{1/p\Xi_N N!}\notin \mathbb Z[[z]]$.
\end{prop}

\begin{proof}
We assume that $p$ is not a Wolstenholme prime and that $p$ does not
divide $N$.
By Lemmas~\ref{lem:H_L}--\ref{lem:5}, we can furthermore
assume that $p\ge7$.

Going back to the outline of the proof of Theorem~\ref{thm:3} in
Section~\ref{sec:Xi}, we claim that
\begin{equation} \label{eq:Cp}
C(p)=\big(B_N(1)H_N-B_N(p)pH_{Np}\big)\notin p^4 N!\,\mathbb{Z}_p,
\end{equation}
where $B_N(m)=\frac {(Nm)!} {m!^N}$.
(The claim here is the non-membership relation; the equality holds by
the definition of $C(\cdot)$ in \eqref{eq:CcongU}.)
This would imply that $C(p)\notin p^2 \Xi_NN!\,\mathbb Z_{p}$, 
and thus, by Lemma~\ref{lem:4}
(recall the outlines of the proofs of Theorems~\ref{thm:2} and
\ref{thm:3} in Sections~\ref{sec:2} and \ref{sec:Xi}, respectively), that
$q_{N,N}(z)^{1/p\Xi_N N!}\notin \mathbb Z[[z]]$, as desired.

To establish \eqref{eq:Cp}, we consider 
\begin{multline*}
H_N(B_N(1)-B_N(p))\\
=H_NN!\bigg(1-\frac 1{(p-1)!^N}\big(1\cdot 2\cdots (p-1)\big)
\big((p+1)\cdot (p+2)\cdots (2p-1)\big)\cdots\\
\times\cdots
\big((pN-p+1)\cdot (pN-p+2)\cdots (pN-1)\big) \bigg).
\end{multline*}
Using $v_{p}(H_N)=3$ and Wilson's theorem, we deduce 
\begin{equation} \label{eq:Zp} 
H_N(B_N(1)-B_N(p))\in p^4N!\,\mathbb Z_{p}.
\end{equation}
However, by Lemma~\ref{lem:congH} and the fact that 
$v_{p}\big(B_N(p)\big)=v_{p}\big(B_N(1)\big)=v_{p}(N!)$, we obtain
$$B_N (1)H_N-B_N(p)pH_{pN}\hbox{${}\not\equiv{}$}
H_N(B_N(1)-B_N(p))\mod p^4N!\,\mathbb Z_{p}.$$
Together with \eqref{eq:Zp}, this yields \eqref{eq:Cp}.
\end{proof}

To summarise the discussion of this section: if the conjecture in
Remarks~\ref{rem:Xi7}(c) that no prime $p$ and integer $N$ exists with
$v_p(H_N)\ge4$ should be true, then Theorem~\ref{thm:3} with $k=1$ is
sharp, that is, the sequence $(t_N)_{N\ge1}$ is given by
$t_N=\Xi_NN!$.

\section{Sketch of the proof of Theorem~\ref{thm:3a}}
\label{sec:Om}

In this section we discuss the proof of Theorem~\ref{thm:3a}.
Since it is completely analogous to the proof of Theorem~\ref{thm:3}
(see Section~\ref{sec:Xi}), we content ourselves to point out the
differences. At the end of the section, we present analogues of
Propositions~\ref{prop:p>N} and \ref{prop:vp=3}, addressing the
question of optimality of Theorem~\ref{thm:3a} with $k=1$.

First of all, by \eqref{eq:truemap}, we have
$$
\big(z^{-1}q_N(z)\big)^{1/kN}=\exp(\tilde G_N(z)/F_N(z)),
$$
where $F_N(z)$ is the series from the Introduction and
\begin{equation*}
\tilde G_N(z):=\sum_{m=1}^{\infty} \frac{(Nm)!^k}{m!^{kN}}
\bigg(H_{Nm}-H_m\bigg)\,z^m.
\end{equation*}

As in Section~\ref{sec:Xi},
we must ``upgrade" the proof of Theorem~\ref{thm:2}
in the special case that $\mathbf N=(N,N,\dots,N)$. Writing as before
$B_{\mathbf N}(m)=\frac{(Nm)!^k}{m!^{kN}}$, we must consider the sum
\begin{multline} \label{eq:Summe}
\tilde C(a+Kp):=\sum_{j=0}^K B_{\mathbf N}(a+jp)B_{\mathbf N}(K-j)
\big((H_{N(K-j)}-H_{K-j})-p(H_{Na+Njp}-H_{a+jp})\big)  \\
=\sum_{j=0}^K B_{\mathbf N}(a+jp)B_{\mathbf N}(K-j)
(H_{N(K-j)}-pH_{Na+Njp}) \kern4cm\\-
\sum_{j=0}^K B_{\mathbf N}(a+jp)B_{\mathbf N}(K-j)
(H_{K-j}-pH_{a+jp}) 
\end{multline}
and show that it is in $p\Om_NN!^k\mathbb Z_p$ for all primes $p$, and for all 
non-negative integers $K$, $a$, and $j$ with $0\le a<p$. 
The special cases $K=a=0$, respectively $K=0$
and $a=1$, are equally simple to be handled directly here. 
We leave their verification to the reader and assume $a+Kp\ge2$ from now on.

Following the outline of the proof of Theorem~\ref{thm:3} in
Section~\ref{sec:Xi}, given
non-negative integers $m$, $K$, $a$, and $j$ with $0\le a<p$ and
$a+Kp\ge2$, we should prove
\begin{multline}\label{eq:firstreductionU2}
\tilde C(a+Kp) \equiv \sum_{j=0}^K B_{\mathbf N}(a+jp)B_{\mathbf N}(K-j) 
\big((H_{N(K-j)}-H_{K-j})-(H_{\lfloor
Na/p\rfloor+Nj}-H_j)\big) \\
\mod p^3 N!^k\mathbb{Z}_p,
\end{multline}
respectively, if $v_p(H_N-1)\ge3$ and $p$ a Wolstenholme prime or
$N\equiv\pm1$~mod~$p$, 
\begin{multline}\label{eq:firstreductionUW2}
\tilde C(a+Kp) \equiv \sum_{j=0}^K B_{\mathbf N}(a+jp)B_{\mathbf N}(K-j) 
\big((H_{N(K-j)}-H_{K-j})-(H_{\lfloor
Na/p\rfloor+Nj}-H_j)\big)\\
\mod p^4 N!^k\mathbb{Z}_p,
\end{multline}
and we need to prove
\begin{equation} \label{eq:congrconj1U2}
B_{\mathbf N}(a+pj)\left(H_{Nj+\lfloor Na/p\rfloor} - H_{Nj}\right)
\in p \Om_N N!^k \mathbb{Z}_p
\end{equation}
and
\begin{equation} \label{eq:110U2}
B_{\mathbf N}(m)\big((H_{Nmp^s}-H_{mp^s})-
(H_{N\fl{m/p}p^{s+1}}-H_{\fl{m/p}p^{s+1}})\big)\in 
p^{-s}\Om_N N!^k \mathbb Z_p\ .
\end{equation}

As the second line of \eqref{eq:Summe} shows, $\tilde C(a+Kp)$ 
is a difference of 
two sums $C(a+Kp)$ as in \eqref{eq:C}, one with $L=N$, the other with
$L=1$. Hence, all the arguments in Section~\ref{sec:Xi} which are not
based on the assumption that $v_p(H_N)>0$ can be used with only
little modification.

On the other hand, 
the assumption $v_p(H_N)>0$ is not used at many places.
First, we need substitutes for Lemmas~\ref{lem:H_L}--\ref{lem:p<L}.

\begin{lem} \label{lem:H_L2}
For all positive integers $L\ge2$, we have $v_2(H_L-1)=-\fl{\log_2L}$.
\end{lem}

\begin{lem} \label{lem:32}
We have $v_3(H_{66}-1)=v_3(H_{68}-1)=1$.
For positive integers $L\notin\{1,66,68\}$, we have $v_3(H_L-1)\le0$.
\end{lem}

\begin{lem} \label{lem:52}
We have $v_5(H_3-1)=v_5(H_{21}-1)=v_5(H_{23}-1)=1$.
For positive integers $L\notin\{1,3,21,23\}$, we have $v_5(H_L-1)\le0$.
\end{lem}

\begin{lem} \label{lem:p<L2}
Let $p$ be a prime, and let $L$ be an integer with $L\ge2$ and $p\le L$.
Then the following assertions hold true:

\begin{enumerate} 
\item If $v_p(H_L-1)>0$ then $L\ge 4p$.
\item If $v_p(H_L-1)>0$ and $p\ne5$ then $L\ge 6p$.
\end{enumerate}
\end{lem}

These results can be proved in exactly the same way as 
Lemmas~\ref{lem:H_L}--\ref{lem:p<L}, respectively.
In comparison to Lemma~\ref{lem:p<L}, the statement of
Lemma~\ref{lem:p<L2} is in fact simpler, so that complications
that arose in Section~\ref{sec:Xi}
(such as \eqref{eq:firstreductionU3}, for example) do not arise here.

Second, we need a substitute for Lemma~\ref{lem:congH}.

\begin{lem} \label{lem:congH2}
For all primes $p\ge5$ and positive integers $N$, we have
\begin{equation} \label{eq:congH12}
p(H_{pN}-H_p)\equiv H_N-1\mod p^4\mathbb Z_p
\end{equation}
if and only if $p$ is a Wolstenholme prime or $N\equiv\pm1$~{\em mod}~$p$.
\end{lem}

\begin{proof}
From the proof of Lemma~\ref{lem:congH}, we know that 
$$
pH_{pN}-H_N\equiv pN^2H_{p-1}\mod \mathbb Z/p^4\mathbb Z.$$
As a consequence, we obtain
$$
p(H_{pN}-H_p)-(H_N-1)\equiv p(N^2-1)H_{p-1}\mod \mathbb Z/p^4\mathbb Z.$$
The assertion of the lemma follows now immediately.
\end{proof}

Finally, the computation for $m=1$ in Subsection~\ref{sec:m=1}
must be replaced by
\begin{align*}
B_{\mathbf N}(1)(H_{Np^s}-H_{p^s})
&=N!^k\Bigg(\frac {1} {p^s}(H_N-1)+
\frac {1} {p^{s-1}}
\underset{p\nmid\ep}{\sum _{\ep=p+1} ^{Np}}\frac {1} {\ep}+
\frac {1} {p^{s-2}}
\underset{p\nmid\ep}{\sum _{\ep=p^2+1} ^{Np^2}}\frac {1} {\ep}+
\underset{p^{s-2}\nmid \ep}{\sum _{\ep=p^s+1} ^{Np^s}}\frac {1}
{\ep}\Bigg),
\end{align*}
with the conclusion that
$$
B_{\mathbf N}(1)(H_{Np^s}-H_{p^s})\in p^{-s} N!^k\Om_N\mathbb Z_p
$$
being found in a completely analogous manner.

Altogether, this leads to a proof of Theorem~\ref{thm:3a}.

\section{The Dwork--Kontsevich sequence}
\label{sec:DK}

In this section, we address the question
of optimality of Theorem~\ref{thm:3a} when $k=1$, that is, whether,
given that $k=1$, the largest integer $u_N$ such that
$\big(z^{-1}q_{N}(z)\big)^{\frac{1}{Nu_N}} \in\mathbb{Z}[[z]]$ 
is given by $\Om_N N!$. 
Let us write $\tilde q_N(z)$ for $\big(z^{-1}q_N(z)\big)^{1/N}$
with $k$ {\it being fixed to $1$.}
The first proposition shows
that there cannot be any prime number $p$ larger than
$N$ which divides $u_N$. We omit the proof since it is
completely analogous to the proof of Proposition~\ref{prop:p>N}
in Section~\ref{sec:DworkKont}.

\begin{prop} \label{prop:p>N2}
Let $p$ be a prime number and $N$ a positive integer with $p>N\ge2$.
Then there exists a positive integer $a<p$ such that
\begin{equation} \label{eq:p>N2} 
v_p\big(B_{N}(a)(H_{Na}-H_a)\big)=0.
\end{equation}
In particular, $p$ does not divide $u_N$. 
\end{prop}

So, if we hope to improve Theorem~\ref{thm:3a} with $k=1$, then it must be by
increasing exponents of prime numbers $p\le N$ in \eqref{eq:Om}. 
According to Remarks~\ref{rem:Om}(b) in the Introduction, an
improvement is therefore only possible if 
$v_p(H_N-1)>2$ for some $p\le N$. 

The final result of this section shows that, if $v_p(H_N-1)=3$
(for which, however, so far no examples are known; see
Remarks~\ref{rem:Om}(b)), the exponent
of $p$ in the definition of $\Om_N$ in \eqref{eq:Om} cannot be increased
so that Theorem~\ref{thm:3a} with $k=1$ would still hold. 

\begin{prop} \label{prop:vp=32}
Let $p$ be a prime number and $N$ a positive integer with $p\le N$
and $v_p(H_N)=3$. If $p$ is not a Wolstenholme prime and 
$N\hbox{${}\not\equiv{}$}\pm1$~{\em mod}~$p$, then
$\tilde q_N(z)^{\frac{1}{p\Om_{N}N!}}
\notin\mathbb{Z}[[z]]$.
\end{prop}

Again, the proof is completely analogous to the proof of 
Proposition~\ref{prop:vp=3} in Section~\ref{sec:DworkKont}, 
which we therefore omit.

So, if the conjecture in
Remarks~\ref{rem:Om}(c) that no prime $p$ and integer $N$ exists
with\break
$v_p(H_N-1)\ge4$ should be true, then Theorem~\ref{thm:3a} with $k=1$ is
optimal, that is, the Dwork--Kontsevich sequence $(u_N)_{N\ge1}$ is given by
$u_N=\Om_NN!$.

\section{Outline of the proof of Theorem~\ref{thm:4}}
\label{sec:7}

In this section, we provide a brief outline of the proof of
Theorem~\ref{thm:4}, reducing it to
Lemmas~\ref{lem:12a}--\ref{lem:strat4}. These lemmas are subsequently
proved in Sections~\ref{sec:9} and \ref{sec:10}, with two auxiliary lemmas
being the subject of the subsequent section.

We follow the strategy that we used in Section~\ref{sec:2} to prove
Theorem~\ref{thm:2}; that is, 
by the consequence of Dwork's Lemma given in
Lemma~\ref{lem:4}, we want to prove that  
$$\mathbf F_{\mathbf N}(z)\mathbf G_{L,\mathbf N}(z^p)
-p\mathbf F_{\mathbf N}(z^p)\mathbf G_{L,\mathbf
N}(z) \in p z \mathbb{Z}_p[[z]].$$ 

The $(a+Kp)$-th Taylor coefficient of
$\mathbf F_{\mathbf N}(z)\mathbf G_{L,\mathbf N}(z^p)
-p\mathbf F_{\mathbf N}(z^p)\mathbf G_{L,\mathbf
N}(z)$ is  
$$
\mathbf C(a+Kp):=\sum_{j=0}^K 
\mathbf B_{\mathbf N}(a+jp)\mathbf B_{\mathbf N}(K-j) 
(H_{L(K-j)}-pH_{La+Ljp}),
$$
where $\mathbf B_{\mathbf N}(m):=\prod_{j=1}^k \mathbf B_{N_j}(m),$
the quantities $\mathbf B_{N_j}(m)$ being defined in \eqref{eq:Bzudilin}.
In view of Lemma~\ref{lem:4}, proving 
Theorem~\ref{thm:4} is equivalent to proving that 
$$\mathbf C(a+Kp) \in p \mathbb{Z}_p$$
for all primes $p$ and non-negative integers $a$ and $K$ with $0\le a<p$.
Again using \eqref{eq:J} with $J=La+Ljp$, we have
$$\mathbf C(a+Kp)  \equiv \sum_{j=0}^K 
\mathbf B_{\mathbf N}(a+jp)\mathbf B_{\mathbf N}(K-j) (H_{L(K-j)}-H_{\lfloor
La/p\rfloor+Lj}) \mod p \mathbb{Z}_p.$$

The analogue of Lemma~\ref{lem:12} in the present context, which allows
us to get rid of the floor function $\fl{La/p}$ and rearrange the sum
over $j$, is the following lemma. The proof can be found in
Section~\ref{sec:9}.

\begin{lem}\label{lem:12a}
For any prime $p$, non-negative integers $a$ and $j$ with $0\le a<p$, 
positive integers $N_1,N_2,\dots,N_k$, and $L\in\{1,2, \ldots, 
\max(N_1, \ldots, N_k)\}$, we have 
$$
\mathbf{B}_{\mathbf{N}}(a+pj)\big(H_{Lj+ \lfloor La/p\rfloor}
-H_{Lj}\big) \in p\mathbb{Z}_p. 
$$
\end{lem}

We now do the same rearrangements as those after Lemma~\ref{lem:12}
to conclude that
$$\mathbf C(a+Kp)\equiv -\sum _{s=0} ^{r}
\sum _{m=0} ^{p^{r+1-s}-1}\mathbf Y_{m,s} \mod p \mathbb{Z}_p,$$
where $r$ is such that $K<p^r$, and 
\begin{equation*} 
\mathbf Y_{m,s}:=\big(H_{Lmp^s}-H_{L\fl{m/p}p^{s+1}}\big)\mathbf S(a,K,s,p,m),
\end{equation*}
the expression $\mathbf S(a,K,s,p,m)$ being defined by
$$
\mathbf S(a,K,s,p,m):=\sum _{j=mp^s} ^{(m+1)p^s-1}
\big(\mathbf B_{\mathbf N}(a+jp)\mathbf B_{\mathbf N}(K-j)
-\mathbf B_{\mathbf N}(j)\mathbf B_{\mathbf N}(a+(K-j)p)
\big).
$$
In this expression for $\mathbf{S}(a,K,s,p,m)$, it is assumed that
$\mathbf{B}_{\mathbf N}(n)=0$ for negative integers~$n.$ 
If we prove that 
\begin{equation} \label{eq:Yms} 
\mathbf Y_{m,s}\in p\mathbb Z_p
\end{equation}
for all $m$ and
$s$, then $\mathbf C(a+Kp)\in p\mathbb Z_p$, as desired.

Now, the last assertion follows from the following two lemmas, with
the proof of the first given in Section~\ref{sec:10}, while the proof
of the second is easily accomplished by 
a trivial adaptation of the proof of Lemma~\ref{lem:11} given  
in Section~\ref{sec:6}, where we change all occurrences 
of $B_{\mathbf N}$ by $\mathbf{B}_{\mathbf{N}}$ and apply Lemma~\ref{lem:12a}
instead of Lemma~\ref{lem:12}. 
The use of Lemma~\ref{lem:multinomial/N!} can be dropped without 
substitute.~(\footnote{\label{foot:1}To prove the refinement announced 
at the end of the Introduction that 
${\mathbf q}_{1,\mathbf{N}}(z)^{1/\mathbf{B}_{\mathbf{N}}(1)}\in
\mathbb{Z}[[z]]$, the use of Lemma~\ref{lem:multinomial/N!} must be 
everywhere replaced by the use of Lemma~\ref{lem:diviBB}, 
the latter being proved in Section~\ref{sec:8}.})
Lemma~\ref{lem:strat3} is the analogue
of Lemma~\ref{lem:10} in the present context, 
while Lemma~\ref{lem:strat4} is the analogue of Lemma~\ref{lem:11}. 

\begin{lem} \label{lem:strat3} 
For all primes $p$ and
non-negative integers $a,m,s,K$ with $0\le a<p$, we have 
$$
\mathbf S(a,K,s,p,m) \in p^{s+1}\mathbf{B}_{\mathbf{N}}(m) \mathbb{Z}_p.
$$
\end{lem}

\begin{lem} \label{lem:strat4} 
For all primes $p$, non-negative integers $m$, 
positive integers $N_1,N_2,\dots,N_k$, and $L\in\{1,2, \ldots, 
\max(N_1, \ldots, N_k)\}$, we have
$$
\mathbf{B}_{\mathbf{N}}(m) \left(H_{Lmp^s}-H_{L\lfloor m/p\rfloor
p^{s+1}}\right) \in \frac{1}{p^s} \,\mathbb{Z}_p. 
$$
\end{lem}
It is clear that Lemmas~\ref{lem:strat3} and~\ref{lem:strat4}
imply~\eqref{eq:Yms}. 
This completes the outline of the proof of Theorem~\ref{thm:4}.

\section{Further auxiliary lemmas}
\label{sec:8}

In this section, we establish four
auxiliary results. The first one,
Lemma~\ref{lem:rajout2}, provides  
the observation~\eqref{eq:rajout2} that reduces Zudilin's
Conjecture~\ref{conj:zudilin2} to Theorem~\ref{thm:4}. 
The second one, Lemma~\ref{lem:10a}, is required for the proof of
Lemma~\ref{lem:12a} in 
Section~\ref{sec:9}, while the third one, Lemma~\ref{lem:ultime},
is required for the proof of Lemma~\ref{lem:strat3} in
Section~\ref{sec:10}. The last result, Lemma~\ref{lem:diviBB}, 
justifies an assertion made in the Introduction. 
Moreover, the proofs of Lemmas~\ref{lem:ultime} and~\ref{lem:diviBB}
make themselves use of Lemma~\ref{lem:10a}.

\begin{lem} \label{lem:rajout2} 
Let $m$ be a non-negative integer, and let $N$ be a positive integer
with associated parameters $\al_i, \beta_i, 
\mu, \eta$ {\em(}that is, given by
\eqref{eq:aj} and \eqref{eq:bj}, respectively{\em)}. Then
\begin{equation*}
\mathbf{H}_N(m)
= \sum_{j=1}^{\mu} \alpha_j H_{\alpha_jm} -  \sum_{j=1}^{\eta} \beta_j
H_{\beta_jm}. 
\end{equation*}
\end{lem}
\begin{proof} 
For $N=1$, we have $\mathbf{H}_1(m)=0$, so that the assertion of
the lemma holds trivially. Therefore, from now on, we assume $N\ge2$.

We claim that, for any real number $m\ge0$, we have
\begin{equation}\label{eq:rajout3}
\frac{C_N^m}{
\Gamma(m+1)^{\varphi(N)}}\prod_{j=1}^{\varphi(N)}\frac{\Gamma(m+r_j/N)}{\Gamma(r_j/N)}
=  
\frac{\prod_{j=1}^\mu \Gamma(\alpha_j m+1)}{\prod_{j=1}^\eta
\Gamma(\beta_j m+1)}, 
\end{equation}
where $\Ga(x)$ denotes the gamma function.
This generalises Zudilin's identity~\eqref{eq:Bzudilin}
to real values of $m$. We essentially extend his proof to
real $m$, using the well-known formula~\cite[p.~23, Theorem~1.5.2]{AAR}
\begin{equation} \label{eq:Ga}
\Ga(a)\,
\Ga\left(a+\frac {1} {n}\right)\,
\Ga\left(a+\frac {2} {n}\right)\cdots
\Ga\left(a+\frac {n-1} {n}\right)=
n^{\frac {1} {2}-an}(2\pi)^{(n-1)/2}\Ga(an),
\end{equation}
valid for real numbers $a$ and positive integers $n$ such that $aN$
is not an integer $\le0$. Indeed, as in the Introduction,
let $p_1, p_2, \dots, p_\ell$ denote the distinct prime factors   
of $N$. 
(It should be noted that there is at least one such prime factor due
to our assumption $N\ge2$.)
Furthermore, for a subset $J$ of $\{1,2,\dots,\ell\}$, let
$p_J$ denote the product $\prod _{j\in J} ^{}p_j$ of corresponding prime
factors of $N$. 
(In the case that $J=\emptyset$, the empty product must be interpreted
as $1$.) Then, by the principle of
inclusion-exclusion, we can rewrite the left-hand side of
\eqref{eq:rajout3} in the form
$$
\frac{C_N^m}{
\Gamma(m+1)^{\varphi(N)}}\cdot
\frac {
\underset{\vert J\vert\text { even}}
{\prod _{J\subseteq \{1,2,\dots,\ell\}} ^{}}
\prod _{i=1} ^{N/p_J}\Ga\left(m+\frac {ip_J} {N}\right)} 
{
\underset{\vert J\vert\text { odd}}
{\prod _{J\subseteq \{1,2,\dots,\ell\}} ^{}}
\prod _{i=1} ^{N/p_J}\Ga\left(m+\frac {ip_J} {N}\right)}
\cdot
\frac {
\underset{\vert J\vert\text { odd}}
{\prod _{J\subseteq \{1,2,\dots,\ell\}} ^{}}
\prod _{i=1} ^{N/p_J}\Ga\left(\frac {ip_J} {N}\right)} 
{
\underset{\vert J\vert\text { even}}
{\prod _{J\subseteq \{1,2,\dots,\ell\}} ^{}}
\prod _{i=1} ^{N/p_J}\Ga\left(\frac {ip_J} {N}\right)}  .
$$
To each of the products over $i$, formula~\eqref{eq:Ga} can be
applied. As a result, we obtain the expression
\begin{multline} \label{eq:Gaexpr}
\frac{C_N^m}{
\Gamma(m+1)^{\varphi(N)}}\cdot
\frac {
\underset{\vert J\vert\text { even}}
{\prod _{J\subseteq \{1,2,\dots,\ell\}} ^{}}
\left(\frac {N} {p_J}\right)^{-\left(m+\frac {p_J} {N}\right)\frac {N}
{p_J}}
\Ga\left(m\frac {N} {p_J}+1\right)} 
{
\underset{\vert J\vert\text { odd}}
{\prod _{J\subseteq \{1,2,\dots,\ell\}} ^{}}
\left(\frac {N} {p_J}\right)^{-\left(m+\frac {p_J} {N}\right)\frac {N}
{p_J}}
\Ga\left(m\frac {N} {p_J}+1\right)} 
\cdot
\frac {
\underset{\vert J\vert\text { odd}}
{\prod _{J\subseteq \{1,2,\dots,\ell\}} ^{}}
\Ga\left(1\right)} 
{
\underset{\vert J\vert\text { even}}
{\prod _{J\subseteq \{1,2,\dots,\ell\}} ^{}}
\Ga\left(1\right)} \\
=
\frac{C_N^m}{
\Gamma(m+1)^{\varphi(N)}}\cdot
\frac {
\underset{\vert J\vert\text { even}}
{\prod _{J\subseteq \{1,2,\dots,\ell\}} ^{}}
\left(\frac {N} {p_J}\right)^{-m {N}/{p_J}}
\Ga\left(m\frac {N} {p_J}+1\right)} 
{
\underset{\vert J\vert\text { odd}}
{\prod _{J\subseteq \{1,2,\dots,\ell\}} ^{}}
\left(\frac {N} {p_J}\right)^{-m {N}/{p_J}}
\Ga\left(m\frac {N} {p_J}+1\right)} ,
\end{multline}
where the simplification in the exponent of $N/p_J$ is due to the
fact that there are as many subsets of even cardinality of a given 
non-empty set as there are subsets of odd cardinality. 
Since, again by inclusion-exclusion, 
\begin{equation} \label{eq:phi}
\underset{\vert J\vert\text { even}}
{\sum _{J\subseteq \{1,2,\dots,\ell\}} ^{}}\frac {N} {p_J}
-\underset{\vert J\vert\text { odd}}
{\sum _{J\subseteq \{1,2,\dots,\ell\}} ^{}}\frac {N} {p_J}
=N\prod _{p\mid N} ^{}\left(1-\frac {1} {p}\right)
=\ph(N),
\end{equation}
we have
$$
\frac{1}{
\Gamma(m+1)^{\varphi(N)}}\cdot
\frac {
\underset{\vert J\vert\text { even}}
{\prod _{J\subseteq \{1,2,\dots,\ell\}} ^{}}
\Ga\left(m\frac {N} {p_J}+1\right)} 
{
\underset{\vert J\vert\text { odd}}
{\prod _{J\subseteq \{1,2,\dots,\ell\}} ^{}}
\Ga\left(m\frac {N} {p_J}+1\right)} =
\frac{\prod_{j=1}^\mu \Gamma(\alpha_j m+1)}{\prod_{j=1}^\eta
\Gamma(\beta_j m+1)}
$$
and
$$
\frac {
\underset{\vert J\vert\text { even}}
{\prod _{J\subseteq \{1,2,\dots,\ell\}} ^{}}
{N}^{-m {N}/{p_J}}}
{
\underset{\vert J\vert\text { odd}}
{\prod _{J\subseteq \{1,2,\dots,\ell\}} ^{}}
{N}^{-m {N}/{p_J}}}=
N^{-m\ph(N)}.
$$
Finally, consider a fixed prime number dividing $N$, 
$p_j$ say. Then, using again \eqref{eq:phi}, we
see that the exponent of $p_j$ in the expression
\eqref{eq:Gaexpr} is
$$
-\frac {m} {p_j}\underset{\vert J\vert\text { odd},\,j\notin J}
{\sum _{J\subseteq \{1,2,\dots,\ell\}} ^{}}\frac {N} {p_J}
+\frac {m} {p_j}\underset{\vert J\vert\text { even},\,j\notin J}
{\sum _{J\subseteq \{1,2,\dots,\ell\}} ^{}}\frac {N} {p_J}
=\frac {mN} {p_j}\underset{p\ne p_j}
{\prod _{p\mid N} ^{}}\left(1-\frac {1} {p}\right)
=\frac {m} {p_j}\frac {\ph(N)} {1-\frac {1} {p_j}}=\frac {m\varphi(N)}
{p_j-1}.
$$
If all these observations are used in \eqref{eq:Gaexpr}, we arrive at
the right-hand side of \eqref{eq:rajout3}.

Now, let us call $b(m)$ the function defined by both sides
of~\eqref{eq:rajout3}, and
let $\psi(x)=\Gamma'(x)/\Gamma(x)$ be the
digamma function.  
We will use the well-known property 
(see~\cite[p.~13, Theorem~1.2.7]{AAR}) that
$
\psi(x+n) - \psi(x) = H(x,n)
$
for real numbers $x>0$ and integers $n\ge 0$.

By taking the logarithmic derivative of the right
hand side of~\eqref{eq:rajout3}, we have  
\begin{align}
\frac{b'(m)}{b(m)} &= \sum_{j=1}^\mu \al_j \psi(\al_jm+1) -
\sum_{j=1}^\eta \be_j \psi(\be_jm+1) \notag 
\\
         &= \sum_{j=1}^\mu \al_j \big(\psi(1)+H_{\al_jm}\big) -
\sum_{j=1}^\eta \be_j\big(\psi(1)+ H_{\be_jm}\big) \notag 
\\
 &= \sum_{j=1}^\mu \al_j H_{\al_jm} - \sum_{j=1}^\eta \be_j
H_{\be_jm}, \label{eq:rajout:4} 
\end{align}
because $\sum_{j=1}^\mu \al_j = \sum_{j=1}^\eta \be_j$. It also
follows that $b'(0)/b(0)=0$. 

On the other hand,  by taking the logarithmic derivative of the  left
hand side of~\eqref{eq:rajout3}, we also have  
\begin{equation*}
\frac{b'(m)}{b(m)} = \log(C_N) + \sum_{j=1}^{\varphi(N)} \psi(m+r_j/N)
- \varphi(N) \psi(m+1). 
\end{equation*}
Since $b'(0)/b(0)=0$, we have 
$
 \log(C_N) = -\sum_{j=1}^{\varphi(N)} \psi(r_j/N) + \varphi(N) \psi(1)
$
and therefore, 
\begin{align}
\frac{b'(m)}{b(m)} &= \sum_{j=1}^{\varphi(N)} \big(\psi(m+r_j/N)-
\psi(r_j/N)\big) - \varphi(N) \big(\psi(m+1)-\psi(1)\big)\notag 
\\
      &=  \sum_{j=1}^{\varphi(N)} H(r_j/N,m) - \varphi(N) H(1,m)\notag
\\ 
      &= \mathbf{H}_N(m).    \label{eq:rajout:5}
\end{align}
The lemma follows by equating the expressions~\eqref{eq:rajout:4}
and~\eqref{eq:rajout:5} obtained for $b'(m)/b(m)$. 
\end{proof}

\begin{lem} \label{lem:10a}
For any integer $N\ge 1$ with associated parameters $\al_i, \beta_i,
\mu, \eta$, the function 
$$\Delta(x):=
 \sum_{i=1}^{\mu}
\left\lfloor \al_{i}x \right\rfloor - \sum_{i=1}^{\eta} 
\left\lfloor \be_{i}x\right\rfloor
$$
has the following properties:
\begin{enumerate}
\item [$(i)$] $\De$ is $1$-periodic. In particular, $\De(n)=0$ for
all integers $n$.
\item [$(ii)$] For all integers $n$, $\De$ is weakly increasing on
intervals of the form $[n,n+1)$.
\item [$(iii)$] For all real numbers $x$, we have $\De(x)\ge0$.
\item [$(iv)$] For all rational numbers $r\neq 0$ whose denominator 
is an element of $\{1,2,\dots,N\}$, we have $\Delta(r)\ge1$.
\end{enumerate}
\end{lem}

\begin{Remark} 
Clearly, the function $\Delta$ is a step function.
The proof below shows that, in fact, all the jumps of $\Delta$ 
at non-integral places have the 
value $+1$ and occur exactly at rational numbers of the form 
$r/N$, with $r$ coprime to $N$. 
\end{Remark} 

\begin{proof}
Property~$(i)$ follows from 
the equality $\sum_{i=1}^\mu \al_i = \sum_{i=1}^\eta \beta_i$
and the trivial fact that $\De(0)=0$.

We turn our attention to property~$(ii)$.
For convenience of notation, let
$$N=p_1^{e_1}p_2^{e_2}\cdots p_\ell ^{e_\ell }$$
be the prime factorisation of $N$, where,
as before, $p_1,p_2,\dots,p_\ell $ are the
distinct prime factors of $N$, and where $e_1,e_2,\dots,e_\ell $ are
positive integers. 

As we already observed in the remark above, 
the function $\Delta$ is a step function. Moreover,
jumps of $\Delta$ can only occur at values of $x$
where some of the $\al_ix$, $1\le i\le\mu$, or some of the $\be_jx$,
$1\le j\le \eta$, (or both) are integers.
At these values of $x$, the value of a (possible) jump is the difference
between the number of $i$'s for which $\al_ix$ is integral and the number
of $j$'s for which $\be_jx$ is integral. In symbols, the value of the
jump is
\begin{equation} \label{eq:diff} 
\#\{i:1\le i\le\mu\text{ and }\al_ix\in\mathbb Z\}
-\#\{j:1\le j\le\eta\text{ and }\be_jx\in\mathbb Z\}.
\end{equation}

Let $X$ be
the place of a jump, $X$ not being an integer. Then we can write $X$
as
$$X=\frac {Z} {p_1^{f_1}p_2^{f_2}\cdots p_\ell ^{f_\ell }},$$
where $f_1,f_2,\dots,f_\ell$ are non-negative integers, not all zero,
and where $Z$ is a non-zero integer relatively prime to
${p_1^{f_1}p_2^{f_2}\cdots p_\ell ^{f_\ell }}$.
Given
$$\al_i=p_1^{a_1}p_2^{a_2}\cdots p_\ell ^{a_\ell }$$
with $e_1+e_2+\dots+e_\ell -(a_1+a_2+\dots+a_\ell )$ even
and $0\le e_k-a_k\le 1$ for each $k=1, 2, \ldots, \ell $,
the number $\al_iX$ will be integral if and only if $a_k\ge f_k$ for all
$k\in\{1,2,\dots,\ell \}$. Similarly, given
$$\be_j=p_1^{b_1}p_2^{b_2}\cdots p_\ell ^{b_\ell }$$
with $e_1+e_2+\dots+e_\ell -(b_1+b_2+\dots+b_\ell )$ odd
and $0\le e_k-b_k\le 1$ for each $k=1, 2, \ldots, \ell $,
the number 
$\be_jX$ will be integral if and only if $b_k\ge f_k$ for all
$k\in\{1,2,\dots,\ell \}$. 
We do not have to take into account the $\be_j$'s which are $1$,
because $1\cdot X=X$ is not an integer by assumption.
For the generating function of vectors
$(c_1,c_2,\dots,c_\ell )$ with $e_k\ge c_k\ge f_k$ 
and $e_{k} - c_{k}\le 1$, 
we have
$$
\sum _{c_1=\max\{e_1-1, f_1\}} ^{e_1}\cdots
\sum _{c_\ell =\max\{e_\ell -1, f_\ell \}} ^{e_\ell }
z^{e_1+\dots+e_\ell -(c_1+\dots+c_\ell )}=
\prod _{k=1} ^{\ell }\left(1+z\cdot\min\{1,e_k-f_k\}\right).
$$
We obtain the difference in \eqref{eq:diff} (with $X$ in place of $x$) 
by putting $z=-1$ on the
left-hand side of this relation. The product on the right-hand side
tells us that this difference is $0$ in case that 
$e_k\ne f_k$ for some $k$, while it is $1$ otherwise. Thus, all the jumps 
of the function $\Delta$ at non-integral places have the value $+1$.

Property~$(iii)$ follows now easily from $(i)$ and $(ii)$.

In order to prove 
$(iv)$, we observe that the first jump of $\Delta$ in the interval
$[0,1)$ occurs at $x=1/N$. Thus, $\Delta(x)\ge1$ for all $x$ in $[1/N,1)$.
This implies in particular that $\Delta(r)\ge1$ for all the above rational
numbers $r$ in the interval $[1/N,1)$. That the same assertion holds
in fact for {\it all\/} the above rational numbers $r$ (not
necessarily restricted to $[1/N,1)$) follows now from the
1-periodicity of the function $\Delta$.
\end{proof}

\begin{lem} \label{lem:ultime} 
For any integers $m,r,w\ge 0$ such that
$0\le w<p^r$, we have  
\begin{equation} \label{eq:ultime}
\frac{\mathbf{B}_{\mathbf{N}}(w+mp^{r})}{\mathbf{B}_{\mathbf{N}}(m)}
\in \mathbb{Z}_p. 
\end{equation}
\end{lem}

\begin{proof}
We first show that we can assume that $m$ is coprime to
$p$. Indeed, let us write $m=hp^t$ with $\gcd(h,p)=1$. We have to prove
that  
$$
\frac{\mathbf{B}_{\mathbf{N}}(w+hp^{r+t})}{\mathbf{B}_{\mathbf{N}}(hp^t)}
\in \mathbb{Z}_p. 
$$
Since
$v_p\big(\mathbf{B}_{\mathbf{N}}(hp^t)/\mathbf{B}_{\mathbf{N}}(h)
\big)=0$
(as can be easily seen from \eqref{eq:Bzudilin} and Legendre's
formula~\eqref{eq:Leg}), this amounts to prove 
that 
$$
\frac{\mathbf{B}_{\mathbf{N}}(w+hp^{r+t})}{\mathbf{B}_{\mathbf{N}}(h)}
\in \mathbb{Z}_p, 
$$
which is the content of the lemma with $r+t$ instead of $r$ and $h$
instead of $m$, with  
$w<p^r<p^{r+t}$.

Therefore, from now on, we assume that $\gcd(m,p)=1$ (however, 
this assumption will only be used after~\eqref{eq:ultime3}). Since
$v_p\big(\mathbf{B}_{\mathbf{N}}(mp^r)/\mathbf{B}_{\mathbf{N}}(m)
\big)=0$, we have 
to prove that 
$$
v_p\left(
\frac{\mathbf{B}_{\mathbf{N}}(w+mp^{r})}{\mathbf{B}_{\mathbf{N}}(mp^r)}\right)
\ge 0 
$$
or, in an equivalent form, that
\begin{multline}\label{eq:ultime2}
\sum_{j=1}^k\sum_{\ell=1}^{\infty} \Bigg( \bigg(\sum_{i=1}^{\mu_j} 
 \left\lfloor \frac{\al_{i,j}(w+mp^r)}{p^\ell}\right\rfloor -
\sum_{i=1}^{\eta_j}  \left\lfloor
\frac{\be_{i,j}(w+mp^r)}{p^\ell}\right\rfloor 
\bigg)
\\
- 
\bigg( \sum_{i=1}^{\mu_j} \left\lfloor \frac{\al_{i,j}
mp^r}{p^\ell}\right\rfloor - \sum_{i=1}^{\eta_j}  \left\lfloor  
\frac{\be_{i,j} mp^r}{p^\ell}\right\rfloor\bigg)
\Bigg) \ge 0,
\end{multline}
where $\al_{i,j}, \beta_{i,j}, \mu_j, \eta_j$ are the parameters
associated to $N_j$.

If $\ell\le r$, then for any $j\in\{1, 2,\ldots, k \}$,  
$$
\sum_{i=1}^{\mu_j} \left\lfloor \frac{\al_{i,j} mp^r}{p^\ell}\right\rfloor 
- \sum_{i=1}^{\eta_j}  \left\lfloor \frac{\be_{i,j}
mp^r}{p^\ell}\right\rfloor   
= mp^{r-\ell} \left(\sum_{i=1}^{\mu_j} \al_{i,j} - \sum_{i=1}^{\eta_j}
\be_{i,j} \right)=0. 
$$
Moreover, 
$$
\sum_{i=1}^{\mu_j} 
 \left\lfloor \frac{\al_{i,j}(w+mp^r)}{p^\ell}\right\rfloor -
\sum_{i=1}^{\eta_j}  \left\lfloor
\frac{\be_{i,j}(w+mp^r)}{p^\ell}\right\rfloor \ge0
$$
because of Lemma~\ref{lem:10a}$(iii)$ with $N=N_j$. It therefore
suffices to show
\begin{multline} \label{eq:ells+1}
\sum_{j=1}^k\sum_{\ell=r+1}^{\infty} \Bigg( \bigg(\sum_{i=1}^{\mu_j} 
 \left\lfloor \frac{\al_{i,j}(w+mp^r)}{p^\ell}\right\rfloor -
\sum_{i=1}^{\eta_j}  \left\lfloor
\frac{\be_{i,j}(w+mp^r)}{p^\ell}\right\rfloor 
\bigg)
\\
- 
\bigg( \sum_{i=1}^{\mu_j} \left\lfloor \frac{\al_{i,j}
mp^r}{p^\ell}\right\rfloor - \sum_{i=1}^{\eta_j}  \left\lfloor  
\frac{\be_{i,j} mp^r}{p^\ell}\right\rfloor\bigg)
\Bigg) \ge 0,
\end{multline}
(The reader should note the difference to \eqref{eq:ultime2}
occurring in the summation bounds for $\ell$.)
For $\ell>r$, set $x_{\ell}=\{mp^r/p^{\ell}\}$
and  $y_{\ell}=\{(w+mp^r)/p^{\ell}\}$. Using again 
$\sum _{i=1} ^{\mu_j}\al_{i,j}-\sum _{i=1} ^{\eta_j}\be_{i,j}=0$,
we see that the left-hand side of~\eqref{eq:ells+1} is equal to 
\begin{equation} \label{eq:ultime3} 
\sum_{j=1}^k \sum_{\ell=r+1}^{\infty} \Bigg( \bigg(\sum_{i=1}^{\mu_j} 
 \left\lfloor \al_{i,j} y_{\ell}\right\rfloor - \sum_{i=1}^{\eta_j}
\left\lfloor \be_{i,j} y_{\ell} \right\rfloor \bigg) 
\\
- 
\bigg( \sum_{i=1}^{\mu_j} \left\lfloor \al_{i,j} x_{\ell}\right\rfloor
- \sum_{i=1}^{\eta_j}  \left\lfloor \be_{i,j} 
x_{\ell}\right\rfloor\bigg)\Bigg).
\end{equation}

We now claim that $x_{\ell}\le y_{\ell}$ for $\ell>r$.
To see this, we begin by the observation that, 
since $m$ and $p$ are coprime and $\ell>r$, the rational number
$m/p^{\ell-r}$ is not an integer. It follows that 
$$
x_{\ell} + \frac{1}{p^{\ell-r}}= \left\{\frac{m}{p^{\ell-r}}\right\}
+ \frac{1}{p^{\ell-r}} \le 1. 
$$
Hence, since $w<p^r$, we infer that 
$$
x_{\ell}+ \frac{w}{p^\ell} <1.
$$
On the other hand, we have
$$
y_{\ell}=\left\{\frac{w}{p^\ell} + \left\lfloor \frac{m}{p^{\ell-r}}
\right\rfloor+ x_{\ell}\right\} =  
\left\{\frac{w}{p^\ell} +  x_{\ell}\right\}=\frac{w}{p^\ell} +  x_{\ell}.
$$
Since $w\ge 0$, we obtain indeed $y_{\ell} \ge x_{\ell}$, as we claimed.

Using $x_{\ell}\le y_{\ell}$ together with Lemma~\ref{lem:10a}$(ii)$, 
we see that, for $\ell>r$ and $j=1,2,\ldots k$, 
we have 
$$
 \sum_{i=1}^{\mu_j} 
 \left\lfloor \al_{i,j} y_{\ell}\right\rfloor - \sum_{i=1}^{\eta_j}
\left\lfloor \be_{i,j} y_{\ell} \right\rfloor   
\ge 
 \sum_{i=1}^{\mu_j} \left\lfloor \al_{i,j} x_{\ell}\right\rfloor -
\sum_{i=1}^{\eta_j}  \left\lfloor \be_{i,j} x_{\ell}\right\rfloor,  
$$
which shows that the expression in \eqref{eq:ultime3} is non-negative,
thus establishing~\eqref{eq:ells+1} and also \eqref{eq:ultime2}. 
This finishes the proof of the lemma.
\end{proof}

We conclude this section with a result which was announced near the
end of the Introduction.
It is used nowhere else, but we give it for the sake of completeness. 
It is a generalisation of Lemma~\ref{lem:multinomial/N!}. By the same
techniques used to prove Theorem~\ref{thm:2},  
it enables one to prove that 
${\mathbf q}_{1,\mathbf{N}}(z)^{1/\mathbf{B}_{\mathbf{N}}(1)}\in
\mathbb{Z}[[z]]$ (see Footnote~\ref{foot:1}). 

\begin{lem}\label{lem:diviBB} 
For any vector $\mathbf{N}$ and any integer $m \ge 1$, we have that 
$\mathbf{B}_{\mathbf{N}}(1)$ divides $\mathbf{B}_{\mathbf{N}}(m)$.
\end{lem}
\begin{proof} Obviously, it is sufficient to prove the assertion for
$k=1$ and $\mathbf{N}=(N)$. Let  
$\Delta$ be the function associated to $N$ as defined in Lemma~\ref{lem:10a}.
We want to prove that, for any prime $p$, we have
$v_p(\mathbf{B}_{N}(m))\ge v_p(\mathbf{B}_{N}(1))$. 
We can assume that $m$ and $p$ are coprime because 
$v_p(\mathbf{B}_{N}(mp^t))=v_p(\mathbf{B}_{N}(m))$ for any 
integers $m,t\ge 0$. 

Now, when $\gcd(m,p)=1$, we have that 
\begin{align*}
v_p(\mathbf{B}_{N}(m))& = \sum_{\ell=1}^{\infty} \Delta(m/p^\ell)
= \sum_{\ell=1}^{\infty} \Delta(\{m/p^\ell\})
\\
&\ge \sum_{\ell=1}^{\infty} \Delta(1/p^\ell) = v_p(\mathbf{B}_{N}(1)).
\end{align*}
Here,
we used the $1$-periodicity of $\Delta$ for the second equality.
For the inequality, we used that 
$\{m/p^\ell\} \ge 1/p^\ell$ (because $\gcd(m,p)=1$ 
implies that $m/p^\ell$ is not an integer) and the 
the (partial) monotonicity of $\De$ described in Lemma~\ref{lem:10a}$(ii)$.
\end{proof}

\section{Proof of Lemma~\ref{lem:12a}}
\label{sec:9}

We follow the first part of the proof of Lemma~\ref{lem:12}. We write
$\al_{i,m}$, $\beta_{i,m}$, $\mu_m$, and $\eta_m$ for
the parameters associated to $N_m$, $m=1,2,\dots,k$. We may assume 
that 
$\max(N_1, \ldots, N_k)=N_k$.  
Then, using again Lemma~\ref{lem:10a}$(iii)$,
\begin{align*}
v_p\big(\mathbf{B}_{\mathbf{N}}(a+pj)\big) & =
\sum_{m=1}^k\sum_{\ell=1}^{\infty}\left( \sum_{i=1}^{\mu_m} 
\left\lfloor \frac{\al_{i,m}(a+pj)}{p^\ell}\right\rfloor - \sum_{i=1}^{\eta_m}
\left\lfloor \frac{\be_{i,m}(a+pj)}{p^\ell}\right\rfloor\right)\\
&\ge \sum_{\ell=1}^{\infty}\left( \sum_{i=1}^{\mu_k}
\left\lfloor \frac{\al_{i,k}(a+pj)}{p^\ell}\right\rfloor - \sum_{i=1}^{\eta_k} 
\left\lfloor
\frac{\be_{i,k}(a+pj)}{p^\ell}\right\rfloor\right)=\sum_{\ell=1}^{\infty}
\Delta_k\!\left(\frac{a+jp}{p^\ell}\right) 
\end{align*}
with 
$$
\Delta_k(x) : = \sum_{i=1}^{\mu_k}
\left\lfloor \al_{i,k}x \right\rfloor - \sum_{i=1}^{\eta_k} 
\left\lfloor \be_{i,k}x\right\rfloor.
$$

We want to prove that 
\begin{equation}\label{eq:rajout8}
v_p\big(\mathbf{B}_{\mathbf{N}}(a+pj)\big) \ge
1+ v_p(Lj+\ep)
\end{equation} 
for any integer $\ep$ such that $1\le \ep\le \lfloor La/p\rfloor.$  
We have 
$$
\frac{a+jp}{p^\ell} = \frac{a-p\ep/L}{p^\ell} + \frac{pj+ p\ep/L }{p^\ell}.
$$

\subsection{First step} We claim that 
\begin{equation} \label{eq:claim1}
\Delta_k\!\left(\frac{a+jp}{p^\ell}\right) \ge  \Delta_k\!\left(\frac{pj+
p\ep/L }{p^\ell}\right). 
\end{equation}
To see this, we first observe that 
$$
\Delta_k\!\left(\frac{a+jp}{p^\ell}\right)=\Delta_k\!\left(\frac{a-p\ep/L}{p^\ell}
+ \frac{pj+ p\ep/L }{p^\ell}\right)  
= \Delta_k\!\left(\frac{a-p\ep/L}{p^\ell} + \left\{\frac{pj+ p\ep/L
}{p^\ell}\right\}\right) 
$$ 
because $\Delta_k$ is $1$-periodic. 

We now claim that 
\begin{equation}\label{eq:inequalities}
0\le \frac{a-p\ep/L}{p^\ell} + \left\{\frac{pj+ p\ep/L }{p^\ell}\right\} < 1.
\end{equation}
Indeed, positivity is clear and we now concentrate on the upper bound.  
We write $j=up^{\ell-1}+v$ with $0\le
v<p^{\ell-1}$. Hence,
$$
 \left\{\frac{pj+ p\ep/L }{p^\ell}\right\}  =  \left\{u+ \frac{pv+
p\ep/L }{p^\ell}\right\} =  
\left\{\frac{v}{p^{\ell-1}} + \frac{p\ep/L }{p^\ell}\right\}.
$$
Since $0 \le \ep \le \lfloor La/p\rfloor<L$, we have  $0\le
\frac{p\ep/L }{p^\ell} <1/p^{\ell-1}$ and therefore 
$$
0\le \frac{v}{p^{\ell-1}} + \frac{p\ep/L }{p^\ell} <
\frac{v}{p^{\ell-1}} + \frac{1}{p^{\ell-1}} \le 1 
$$
(where the last inequality holds by definition of $v$), whence
$$
\left\{\frac{pj+ p\ep/L }{p^\ell}\right\} = \frac{pv+ p\ep/L }{p^\ell}.
$$
Therefore, we have 
$$
 \frac{a-p\ep/L}{p^\ell} + \left\{\frac{pj+ p\ep/L }{p^\ell}\right\} =
\frac{a-p\ep/L}{p^\ell} +  \frac{pv+ p\ep/L }{p^\ell}  
= \frac{a}{p^\ell} + \frac{v}{p^{\ell-1}}.
$$
Since $\frac{v}{p^{\ell-1}}<1$ and $a<p$, we necessarily have  
$$
\frac{a}{p^\ell} + \frac{v}{p^{\ell-1}} < 1,
$$
as desired.

Since $\frac{a-p\ep/L}{p^\ell}\ge 0$, it follows from 
Lemma~\ref{lem:10a}$(i)$,$(ii)$ (with $\Delta=\Delta_k$) and
\eqref{eq:inequalities} that  
$$
\Delta_k\!\left(\frac{a-p\ep/L}{p^\ell} + \left\{\frac{pj+ p\ep/L
}{p^\ell}\right\}\right) \ge  
\Delta_k\!\left( \left\{\frac{pj+ p\ep/L }{p^\ell}\right\} \right) =
\Delta_k\!\left( \frac{pj+ p\ep/L }{p^\ell} \right). 
$$
Thus, we have proved the claim~\eqref{eq:claim1} 
made at the beginning of this step.

\subsection{Second step} Let us write $Lj+\ep=\beta p^d$,
where $d=v_p(Lj+\ep)$, so that  
$$
\frac{pj+ p\ep/L }{p^\ell} = \frac{\beta p^{d+1-\ell}}{L}.
$$
We have proved in the first step that 
\begin{equation}\label{eq:rajout7}
v_p\big(\mathbf{B}_{\mathbf{N}}(a+pj)\big) \ge \sum_{\ell=1}^\infty
\Delta_k\!\left(\frac{\beta p^{d+1-\ell}}{L}\right). 
\end{equation}
By the same argument as the one that 
we used in the first part of the proof of Lemma~\ref{lem:12}
in Section~\ref{sec:4}, the rational
number $\frac{\beta p^{d+1-\ell}}{L}$ is not an integer and for  
$\ell\le d+1$, the denominator of  $\frac{\beta p^{d+1-\ell}}{L}$ is
at most $L$. Since $L\le N_k$, it follows then from
Lemma~\ref{lem:10a}$(iv)$, again with $\Delta=\Delta_k$,
that $\Delta_k(\beta p^{d+1-\ell}/L) \ge 1$ for any $\ell$ in
$\{1, 2, \ldots,  d+1\}.$  
Use of this estimation in \eqref{eq:rajout7} gives
$$v_p\big(\mathbf{B}_{\mathbf N}(a+pj)\big)\ge d+1=1+v_p(Lj+\ep).$$
This completes the proof of~\eqref{eq:rajout8} and, hence, of
Lemma~\ref{lem:12a}.

\section{Proof of Lemma~\ref{lem:strat3}}
\label{sec:10}

Again, we want to use Proposition~\ref{prop:dworkcongruence}, this
time with $A_r(m)=g_r(m)=\mathbf B_{\mathbf N}(m)$.
Clearly, the proposition would imply
that $\mathbf S(a,K,s,p,m)\in p^{s+1}
\mathbf B_{\mathbf N}(m)\mathbb Z_p$, and, thus, the claim.
So, we need to
verify the conditions $(i)$--$(iii)$ in the statement of the proposition.

Conditions $(i)$ and $(ii)$ of Dwork's congruences theorem 
clearly hold, and we must check Condition $(iii)$. 
To do the latter, we follow the  
method developed to prove Lemma~\ref{lem:10}
(see Section~\ref{sec:5}). The reader should recall that
$$
\mathbf{B}_{\mathbf{N}}(m) := \prod_{j=1}^k \mathbf{B}_{N_j}(m),
$$
where $\mathbf B_{N_j}(m)$ is given by \eqref{eq:definitionBNgras} 
or \eqref{eq:Bzudilin}.

\subsection{First step} Let us fix $j\in \{1, 2, \ldots, k\}$ 
and set $D_N:=N^{-\varphi(N)}C_N$, which is an integer. 
We claim that
\begin{equation}\label{eq:step1} 
\frac{\mathbf{B}_{N_j}(v+up+np^{s+1})}{\mathbf{B}_{N_j}(up+np^{s+1})} = 
\frac{\mathbf{B}_{N_j}(v+up)}{\mathbf{B}_{N_j}(up)} + \mathcal{O}(p^{s+1}).
\end{equation}
To prove~\eqref{eq:step1}, we 
can use the same arguments as in the first
step of the proof of Lemma~\ref{lem:10}, 
due to  
the identities~(\footnote{They are immediate consequences of 
the definition~\eqref{eq:definitionBNgras} of $\mathbf{B}_{N}$. Zudilin 
used them in his proof of the following stronger version
of~\eqref{eq:step1}: 
$$
\frac{\mathbf{B}_{N_j}(v+up+np^{s+1})}{\mathbf{B}_{N_j}(up+np^{s+1})} = 
\frac{\mathbf{B}_{N_j}(v+up)}{\mathbf{B}_{N_j}(up)}\,\big(1 +
\mathcal{O}(p^{s+1})\big). 
$$
However, for this, he assumes that $p$ divides $N_j$ (see~\cite[Eq.~(35)]{zud}).
Here, we do not assume that $p$ divides $N_j$, and therefore we obtain
the weaker equality~\eqref{eq:step1}, which is fortunately enough for
our purposes.}) 
\begin{align}
\frac{\mathbf{B}_{N_j}(v+up+np^{s+1})}{\mathbf{B}_{N_j}(up+np^{s+1})} 
&=  
\frac{D_{N_j}^v \prod_{\ell=1}^{\varphi(N_j)}\prod_{i=1}^v 
\big(r_\ell+(i-1)N_j+uN_jp + nN_j p^{s+1}\big)}
{\big((v+up+np^{s+1})(v-1+up+np^{s+1})\cdots (1+up+np^{s+1})\big)^{\varphi(N_j)}}, \notag
\\
\label{eq:rajout9}
\frac{\mathbf{B}_{N_j}(v+up)}{\mathbf{B}_{N_j}(up)} &=  
\frac{D_{N_j}^v \prod_{\ell=1}^{\varphi(N_j)}\prod_{i=1}^v 
\big(r_\ell+(i-1)N_j+uN_jp\big)}
{\big((v+up)(v-1+up)\cdots (1+up)\big)^{\varphi(N_j)}},
\end{align}
and 
to the fact that $(v+up)(v-1+up)\cdots (1+up)$ is not
divisible by $p.$ 

A side result of~\eqref{eq:rajout9} (which was actually 
used to prove~\eqref{eq:step1}) is that  
$$
\frac{\mathbf{B}_{N_j}(v+up)}{\mathbf{B}_{N_j}(up)} \in \mathbb{Z}_p.
$$
We deduce from this fact and from~\eqref{eq:step1} that 
$$
\prod_{j=1}^k
\frac{\mathbf{B}_{N_j}(v+up+np^{s+1})}{\mathbf{B}_{N_j}(up+np^{s+1})}
=
\prod_{j=1}^k
\left(\frac{\mathbf{B}_{N_j}(v+up)}{\mathbf{B}_{N_j}(up)} +
\mathcal{O}(p^{s+1})\right)  
= \prod_{j=1}^k \frac{\mathbf{B}_{N_j}(v+up)}{\mathbf{B}_{N_j}(up)} +
\mathcal{O}(p^{s+1}) ,
$$
or, in other words, 
\begin{equation}\label{eq:main1}
\frac{\mathbf{B}_{\mathbf{N}}(v+up+np^{s+1})}
{\mathbf{B}_{\mathbf{N}}(up+np^{s+1})}
=  
\frac{\mathbf{B}_{\mathbf{N}}(v+up)}{\mathbf{B}_{\mathbf{N}}(up)} +
\mathcal{O}(p^{s+1}). 
\end{equation}

\subsection{Second step} Let us fix $j\in \{1, 2, \ldots, k\}$. 
Exactly as in the proof of the second step of Lemma~\ref{lem:10}, 
the properties of $\Gamma_p$ imply that 
\begin{align*}
\frac{\mathbf{B}_{N_j}(up+np^{s+1})}{\mathbf{B}_{N_j}(u+np^{s})} 
&= (-1)^{\mu_j-\eta_j}\frac{\prod_{i=1}^{\mu_j}\Gamma\big(1+\al_{i,j}(up+np^{s+1})\big)}
{\prod_{i=1}^{\eta_j}\Gamma\big(1+\be_{i,j}(up+np^{s+1})\big)} \notag \\
&= (-1)^{\mu_j-\eta_j}\frac{\prod_{i=1}^{\mu_j}\Gamma(1+\al_{i,j}up)+\mathcal{O}(p^{s+1})}
{\prod_{i=1}^{\eta_j}\Gamma(1+\be_{i,j}up)+\mathcal{O}(p^{s+1})} \notag
\\
&= (-1)^{\mu_j-\eta_j}\frac{\prod_{i=1}^{\mu_j}\Gamma(1+\al_{i,j}up)}
{\prod_{i=1}^{\eta_j}\Gamma(1+\be_{i,j}up)}\,\big(1 +\mathcal{O}(p^{s+1})\big)\notag
\\
&= \frac{\mathbf{B}_{N_j}(up)}{\mathbf{B}_{N_j}(u)} \,\big(1 +\mathcal{O}(p^{s+1})\big).
\end{align*}
Here we used again the fact that $\Gamma_p(m)$ is never divisible by
$p$ for any integer $m$. 

Hence, taking the product over $j=1, 2, \ldots, k$, we obtain 
\begin{equation}\label{eq:main2}
\frac{\mathbf{B}_{\mathbf{N}}(up+np^{s+1})}{\mathbf{B}_{\mathbf{N}}(u+np^{s})}
=  
\frac{\mathbf{B}_{\mathbf{N}}(up)}{\mathbf{B}_{\mathbf{N}}(u)}\,\big(1
+ \mathcal{O}(p^{s+1})\big). 
\end{equation}

\subsection{Third step} We follow verbatim the beginning of the
third step in the proof of Lemma~\ref{lem:10} (i.e.,  
we multiply the left-hand and right-hand sides
of~\eqref{eq:main1} and~\eqref{eq:main2}, etc.) to obtain 
$$
\frac{\mathbf{B}_{\mathbf{N}}(v+up+np^{s+1})}{\mathbf{B}_{\mathbf{N}}(v+up)} 
-
\frac{\mathbf{B}_{\mathbf{N}}(u+np^{s})}{\mathbf{B}_{\mathbf{N}}(u)}  
=
\frac{\mathbf{B}_{\mathbf{N}}(u+np^{s})}{\mathbf{B}_{\mathbf{N}}(u)}\,
\mathcal{O}(p^{s+1})  
+
\frac{\mathbf{B}_{\mathbf{N}}(u+np^{s})}{\mathbf{B}_{\mathbf{N}}(v+up)}\,
\mathcal{O}(p^{s+1}). 
$$

It remains to check that 
\begin{equation}\label{eq:check}
\frac{\mathbf{B}_{\mathbf{N}}(u+np^{s})}{\mathbf{B}_{\mathbf{N}}(u)} \in 
\frac{\mathbf{B}_{\mathbf{N}}(n)}{\mathbf{B}_{\mathbf{N}}(v+up)} \,
\mathbb{Z}_p 
\quad
\textup{and}
\quad
\frac{\mathbf{B}_{\mathbf{N}}(u+np^{s})}{\mathbf{B}_{\mathbf{N}}(v+up)}\in 
\frac{\mathbf{B}_{\mathbf{N}}(n)}{\mathbf{B}_{\mathbf{N}}(v+up)} \,
\mathbb{Z}_p. 
\end{equation}
The first assertion in~\eqref{eq:check} can be rewritten as
\begin{equation} \label{eq:ass1} 
\frac{\mathbf{B}_{\mathbf{N}}(u+np^{s})}{\mathbf{B}_{\mathbf{N}}(n)}\cdot
\frac{\mathbf{B}_{\mathbf{N}}(v+up)}{\mathbf{B}_{\mathbf{N}}(u)} \in
\mathbb{Z}_p ,
\end{equation}
while the second assertion can be rewritten as 
\begin{equation} \label{eq:ass2} 
\frac{\mathbf{B}_{\mathbf{N}}(u+np^{s})}{\mathbf{B}_{\mathbf{N}}(n)}\in
\mathbb{Z}_p.
\end{equation}
Now, the assertion \eqref{eq:ass2} is the special case 
$w=u$, $m=n$ and $r=s$ of 
Lemma~\ref{lem:ultime}, while \eqref{eq:ass1} follows from 
\eqref{eq:ass2} combined with the special case 
$w=v$, $m=u$ and $r=1$ of Lemma~\ref{lem:ultime}.

\medskip
This completes the proof of the lemma.

\section*{Acknowledgements}
The authors are extremely grateful to Alessio Corti 
and Catriona Maclean for illuminating discussions 
concerning the geometric side of our work, 
and to David Boyd for helpful information on
the $p$-adic behaviour of the harmonic numbers $H_N$ 
and for communicating to us the value \eqref{eq:boyd} from his files
from 1994.

\def\refname{Bibliography}


\begin{thebibliography}{1}\label{sec:biblio}
\addcontentsline{toc}{section}{Bibliographie}

\bibitem{almkvist} G. Almkvist and W. Zudilin, {\em 
Differential equations, mirror maps and zeta values}, in: Mirror
Symmetry~V,  N. Yui, S.-T. Yau,  
and J.D. Lewis (eds.), AMS/IP Studies in Advanced Mathematics {\bf 38} (2007), 
International Press \& Amer. Math. Soc., 481--515.

\bibitem{andre} Y. Andr\'e, {\em G-fonctions et transcendance},
J. reine angew. Math. {\bf 476} (1996), 95--125. 

\bibitem{AAR} G.~E.~Andrews, R.~A.~Askey et R.~Roy, {\em Special
Functions},
The Encyclopedia of Mathematics and Its Applications, vol.~{\bf 71},
(G.-C.~Rota, ed.), Cambridge University Press, Cambridge (1999).

\bibitem{batstrat} V. V. Batyrev and D. van Straten,
{\em Generalized hypergeometric functions and rational curves on Calabi--Yau
complete intersections in toric varieties},
Comm. Math. Phys. {\bf 168} (1995), 493--533.

\bibitem{boyd} D. W. Boyd, 
{\em A $p$-adic study of the partial sums of the harmonic series},
Experimental Math.\ {\bf 3} (1994), 287--302.

\bibitem{candelas} P. Candelas, X. de la Ossa, P. Green and L. Parkes,
{\em A pair of Calabi--Yau  
manifolds as an exactly soluble superconformal theory}, 
Nucl. Phys. {\bf B359} (1991), 21--74.

\bibitem{corti} A. Corti and V. Golyshev, {\em Hypergeometric 
equations and weighted projective spaces}, preprint 2006; available at
{\tt http://ar$\chi$iv.org/abs/0607.5016}.

\bibitem{cott} E. Cotterill, {\em Rational curves of degree $10$ on a general quintic threefold}, 
 Comm. Algebra  {\bf 33}  (2005),  no. 6, 1833--1872.

\bibitem{doran} C. F. Doran, 
{\em Picard--Fuchs uniformization and modularity of the mirror map}, 
Comm. Math. Phys. {\bf 212} (2000), 625--647.

\bibitem{dworkihes} B. Dwork, {\em p-Adic cycles}, Publ. I.H.E.S. {\bf
37} (1969), 27--115. 

\bibitem{dwork} B. Dwork, {\em On $p$-adic differential equations IV:
generalized  
hypergeometric functions as $p$-adic analytic functions in one
variable}, Ann.\ Sci.\ \'E.N.S. (4) {\bf 6}, no.~3 (1973), 295--316.

\bibitem{givental} A. B. Givental, {\em 
Equivariant Gromov--Witten invariants}, 
Internat. Math. Res. Notices {\bf 13} (1996), 613--663.

\bibitem{hw} G.~H.~Hardy et E.~M.~Wright,
{\it An introduction to the theory of number theory},
5th edition, Oxford University Press, 1979.

\bibitem{HeRSAA}
N. Heninger, E. M. Rains and N. J. A. Sloane, 
{\em On the integrality of $n$th roots of generating functions},
J. Combin.\ Theory Ser.~A {\bf 113} (2006), 1732--1745. 

\bibitem{hosono} S. Hosono, A. Klemm, S. Theisen, S.-T. Yau, 
{\em Mirror symmetry, mirror map and applications to complete intersection 
Calabi--Yau spaces}, 
Nuclear Phys. B {\bf 43}, no.~3 (1995), 501--552.

\bibitem{KoSVAA} M. Kontsevich, A. Schwarz and V. Vologodsky, 
{\em Integrality of instanton numbers and $p$-adic $B$-model},
Phys. Lett. B {\bf 637} (2006), 97--101. 

\bibitem{lang} S. Lang, {\em Cyclotomic fields II}, Graduate Texts in
Mathematics {\bf 69}, Springer--Verlag, 1980.  

\bibitem{lianyau1} B. H. Lian and S.-T. Yau, {\em Arithmetic properties
of mirror map and  
quantum coupling}, Comm. Math. Phys. {\bf 176}, no.~1 (1996),
163--191.

\bibitem{lianyau} B. H. Lian and S.-T. Yau, 
{\em Mirror maps, modular relations and hypergeometric series I},
appeared as {\em Integrality  
of certain exponential series}, in: Lectures in Algebra and Geometry,
Proceedings of the International Conference on Algebra and Geometry,
Taipei, 1995, M.-C.~Kang (ed.), 
Int. Press, Cambridge, MA, 1998, pp.~215--227.

\bibitem{lianliuyau} B. H. Lian, K. Liu and S.-T. Yau, {\em 
Mirror principle. I.}, Asian J. Math. {\bf 1}, no.~4 (1997), 729--763.

\bibitem{lianyau2} B. H. Lian and S.-T. Yau, 
{\em The $n$th root of the mirror map}, in: Calabi--Yau varieties and
mirror symmetry, Proceedings of the Workshop on Arithmetic, Geometry
and Physics around Calabi--Yau Varieties and Mirror Symmetry,
Toronto, ON, 2001, 
N.~Yui and J.~D.~Lewis (eds.),
Fields Inst. Commun., {\bf 38}, Amer. Math. Soc., Providence, RI,
2003, pp.~195--199.

\bibitem{morrisonjams} D. R. Morrison, {\em 
Mirror symmetry and rational curves on quintic threefolds: a guide for
mathematicians}, J. Amer. Math. Soc. {\bf 6} (1993), 223--247. 

\bibitem{oeis} The On-line encyclopedia of integers sequences, {\tt
http://www.research.att.com/\~{}njas/sequences/}.

\bibitem{pandha} R. Pandharipande, {\em Rational curves on 
hypersurfaces (after A. Givental)}, S\'eminaire Bourbaki. Vol.~1997/98.
Ast\'erisque No.~{\bf 252} (1998), Exp. No.~848, 5, 307--340.

\bibitem{stienstra} J. Stienstra, {\em GKZ Hypergeometric Structures},
Proceedings of the Summer School  
``Algebraic Geometry and Hypergeometric Functions'', Istanbul, June
2005, preprint 2005; available at {\tt http://ar$\chi$iv.org/abs/0511.5351}.

\bibitem{villegas} F. Rodriguez--Villegas, {\em 
Hypergeometric families of Calabi--Yau manifolds}, 
 Proceedings of the Workshop on Arithmetic, Geometry
and Physics around Calabi--Yau Varieties and Mirror Symmetry,
Toronto, ON, 2001, 
N.~Yui and J.~D.~Lewis (eds.),
Fields Inst. Commun., {\bf 38}, Amer. Math. Soc., Providence, RI,
2003, 223--231.

\bibitem{VoisAA}
C. Voisin, {\em Mirror symmetry},
SMF/AMS Texts and Monographs, vol.~1, American
Mathematical Society, Providence, RI; Soci\'et\'e Math\'ematique de
France, Paris, 1999.

\bibitem{volog} V. Vologodsky, {\em Integrality of instanton numbers},
preprint 2007; 
{\tt http://ar$\chi$iv.org/abs/0707.4617}.

\bibitem{yoshida} M. Yoshida, {\em Fuchsian differential equations}, 
Aspects of Mathematics {\bf 11}, Vieweg, 1987.

\bibitem{zud} W. Zudilin, {\em Integrality of power expansions related
to hypergeometric series},  
Mathematical Notes {\bf 71}.5 (2002), 604--616.
\end{thebibliography}
\end{document}